\input psfig
\input amstex
\magnification 1200
\documentstyle{amsppt}
\Monograph
\pageheight{9truein}
\pagewidth{6truein}
\NoBlackBoxes
\NoRunningHeads
\input cyracc.def
\catcode`\@=11
\font@\tencyr=wncyr10
\catcode`\@=13
\addto\tenpoint{\def\cyr{\tencyr\cyracc}}

\def\Ug{U_q\frak g}

\def\half{{\frac{1}{2}}}
\def\1{\bold 1}
\def\sltwo{\frak {sl} _2 }
\def\sln{\frak{sl}_n}
\def\a{\alpha}
\def\l{\lambda}
\def\lk{{\lambda^k}}
\def\g{\frak g}
\def\h{\frak h}
\def\Z{\Bbb Z}
\def\C{\Bbb C}
\def\R{\Bbb R}
\def\N{\Bbb N}

\def\V{\bold V}
\def\eps{\Lambda_0}
\def\d{\partial}
\def\i{\text{\rm i}}
\def\ghat{\hat\frak g}
\def\slnhat{\widehat{\frak{sl}_n}}
\def\hhat{\hat\frak h}
\def\gtilde{\tilde{\frak g}}

\def\H{\Cal H}
\def\tJ{\tilde J}
\def\Mhat{{\widehat M}}
\def\What{{\widehat W}}
\def\Phat{\widehat{P}}
\def\Pp{\widehat P'}
\def\Rhat{\widehat R}
\def\ahat{{\hat \alpha}}
\def\dhat{\hat\delta}
\def\rhat{{\hat\rho}}
\def\lhat{{\widehat{\lambda}}}
\def\Tr{\operatorname{Tr}}
\def\Id{\operatorname{Id}}
\def\ch{\operatorname{ch }}
\def\Hom{\operatorname{Hom }}
\def\dim{\operatorname{dim }}
\def\lot{{\text{\it lower order terms}}}
\def\v{^\vee}
\def\<{\langle}
\def\>{\rangle}
\def\o{\otimes}
\redefine\headmark#1{}
\def\tiphi{\tilde\Phi^{\mu, u}_\l}

\hcorrection{0.2truein}
\document

\rightheadtext{}
\leftheadtext{}
\vbox{\vskip 1in}
\specialhead 
\centerline{ ACKNOWLEDGEMENTS} 
\endspecialhead

\vbox{\vskip 0.5in}

I would like to express my deep gratitude to 
my advisor Igor Frenkel for his guidance and encouragement throughout
my study at Yale. Most of what I know  about
representation theory, mathematical physics and special functions comes
from him. 

Also, I want to thank my friend and co-author Pavel Etingof; most of
the results in this dissertation are based on our joint work, and
would never be written without him; the idea of systematic study of
the weighted traces of intertwiners belongs to him. His enthusiasm and
new ideas often helped to find  a new way when I was ready to give up.

My understanding of the subjects discussed in this dissertation
benefited a lot from fruitful discussions with other people, and
among them Ivan Cherednik, Jintai Ding, Howard Garland, Ian Grojnowski, 
 David Kazhdan,  Ian Macdonald,   Masatoshi Noumi,  Gregg Zuckerman 
 and many others. 

Financial support during my first three years of study was provided by
 Yale and during the final year by Alfred P.~Sloan graduate  
dissertation fellowship. 

Finally -- and most important of all -- I want to thank my wife {\cyr
Varya} for her love and patience -- even when I started discussing
mathematics with Pavel while skiing in New Hampshire...

\vfill\newpage
\rightheadtext{INTRODUCTION }
\leftheadtext{INTRODUCTION}
\vbox{\vskip 1in}
\specialhead 
\centerline{INTRODUCTION} 
\endspecialhead

\vbox{\vskip 0.5in}

This dissertation is devoted to study of some structures appearing in
representation theory of simple Lie algebras (quantum groups, affine
Lie algebras) and their relations with the theory of special functions
and mathematical physics. This subject, of course, has a long
history. 
A representation-theoretic approach to special
functions was developed in the 40-s and 50-s in the works
of I.M.Gelfand, M.A.Naimark, N.Ya.Vilenkin, and their collaborators
(see \cite{V}, \cite{VK}). The essence of this approach is the fact 
that most classical special functions can be obtained as suitable
specializations of matrix elements or characters
of representations of groups. In geometrical terms, interesting
special functions appear as spherical functions on symmetric spaces
associated with $G$, which was studied by Harish-Chandra and Helgason
(see \cite{HC}, \cite{W}). 

However, this approach does not cover all interesting cases. For
example, it was shown in recent works on representations of (quantum)
affine Lie algebras that matrix elements of 
intertwining operators between certain representations of these
algebras are interesting special functions --
($q$-)hypergeometric functions and their
generalizations \cite{TK, FR}. 

In this dissertation we present a more general scheme, which includes
both the classical theory of characters and matrix elements and
intertwining operators. This scheme was suggested in the joint paper
of the author and Pavel Etingof \cite{EK1} and developed in papers
\cite{EK2-- EK5, E1, K1, K2}. Briefly it can be described as ``theory
of characters for intertwining operators'', which we call generalized
characters. This approach allows to get in a natural  way many 
interesting special functions (for example,  Lam\'e functions, Macdonald's
polynomials) and easily prove a number of their properties. 

Here are the main ideas of this approach. Let $\g$ be a  simple
finite-dimensional Lie algebra over $\C$, and let $G$ be the
corresponding compact real Lie group. An important role in the representation 
 theory of $G$ play  the characters 
of finite-dimensional representations; in particular, they form a
distinguished basis in the space of  functions on $G$ invariant with respect to
conjugation; in other words, they are zonal spherical functions on the
symmetric space $G\times G/G_{diag}$. We can generalize 
this considering {\it equivariant}
functions, i.e. functions on $G$ with values in a representation $U$
which satisfy 

$$f(gxg^{-1})=g f(x).$$

This is a particular example of so-called $\mu$-spherical functions
which were studied in papers of Harish-Chandra (see
\cite{HC, W}).

It turns out that this space has equally natural description. Namely,
let $V$ be any finite-dimensional representation of $G$, and
$\Phi:V\to V\o U$ be an intertwining operator. Define the
corresponding generalized character $\chi_\Phi$ by

$$\chi_\Phi(g)=\Tr_{V}(\Phi g).$$

This is a natural
generalization of usual characters; however, it takes values not in $\C$
but in the representation $U$. It is easy to see that $\chi_\Phi$ is
equivariant; moreover, every equivariant function is a linear
combination of generalized characters. Note that due to the
equivariance, the 
generalized character is completely determined by its values on the
elements of the form $g=e^h, h$ from Cartan subalgebra of $\g$. For
this reason, from now on we consider $\chi_\Phi$ as a function on
Cartan subalgebra. 

In  the classical case many properties of
representations could be expressed in terms of their characters. The
same holds here: many properties of intertwining operators can be
expressed in terms of corresponding generalized characters.  For
example, the natural inner product on the space of intertwining
operators  gives an inner product on generalized characters and an 
analogue of the orthogonality theorem for them. 
Similarly, the fact that $\Phi$ commutes with  the center of $U\g$ 
yields differential equations satisfied by the
$\chi_\Phi$ -- each element of the center gives a differential
equation. Thus, the generalized characters are common eigenfunctions of
a family of commuting differential operators, obtained from the 
 center of  $U\g$. 

These generalized characters are a rich source of special functions.
In general, these  functions have not been studied before. However,
in some special cases they coincide with well-known objects. For example,
if $U$ is a trivial module then this is nothing but the usual characters
and thus we recover the classical theory.

More generally, the generalized characters (considered as functions on
the Cartan subalgebra) take values in the zero-weight
subspace $U[0]$; thus, if we consider the case of Lie algebra $\sln$ and
take $U$ to be the space of homogeneous polynomials in $n$ variables of
degree $(k-1)n$ then $U[0]$ is one-dimensional and thus the generalized
characters can be considered as complex-valued. In this case we show 
that the ratio $\chi_\Phi/\delta^{k-1}$, where $\delta$ is
the Weyl denominator, is the  Jack symmetric polynomial.
The simplest way to describe Jack  polynomials -- or more generally,
Jacobi polynomials, which are generalization of Jack polynomials to
the case of arbitrary root systems introduced by Heckman
and Opdam  -- is to say that they are  eigenfunctions of  
(conjugated by certain function) Sutherland operator:

$$L_k=\Delta - 2k(k-1)\sum_{\alpha\in R^+}\frac{1}
{(2\sinh \frac{(\a, h)}{2})^2},$$
where $\Delta$ is the Laplace operator on $\h$.
For special values of $k$ this operator
is just the radial part of the Laplace operator on certain symmetric spaces
associated with the group $G$, and the eigenfunctions are zonal spherical
functions. Heckman and Opdam have showed that for any $k$, these
polynomials have a number of
remarkable properties. In particular, they showed that the Sutherland
operator can be included in a commutative family of differential
operators, isomorphic to the algebra $S(\h)^W\simeq \Cal Z(U \g)$ of Weyl
group invariant polynomials on $\h$; also, these polynomials are
orthogonal with respect to a certain inner product. All of these results
are highly non-trivial and required a good deal of ingenuity. We will
show  that our representation-theoretic
approach 
 allows to obtain both these properties immediately as a
corollary of the general results about the generalized characters. 
This new interpretation also suggests new generalizations, such as
vector-valued analogue of these polynomials or affine analogue (see
below). 

In a similar way, the generalized characters can be defined for the
quantum group $U_q 
\g$, corresponding to $\g$. All the constructions above can be
generalized to this case with some minor changes; most importantly,
differential operators should be replaced by difference operators. We
will show  that for the case $\g=\sln$ and
$U$ chosen as above we get 
the famous Macdonald's polynomials -- a family of symmetric
polynomials which was recently introduced by I.~Macdonald and has been
the object of intensive study since that time. Again, this approach
allows us to reprove many properties of Macdonald's polynomials (inner
product identity, symmetry identity) in a very simple way.

In a similar way, we can define generalized characters for affine Lie
algebras. It turns out that if we take $U$ to be tensor product of
evaluation representations then these characters are precisely the
correlation functions on the torus of the Wess-Zumino-Witten (WZW)
model of conformal filed theory. We deduce differential equations,
describing dependence of these generalized characters on the
modular parameter of the torus and the parameters of the evaluation
representations; we call these equations elliptic
Knizhnik-Zamolodchikov equations (this was first done, in a different
language, by Bernard \cite{Be}) and study their monodromies. 

Again, in the case $\g=\sln$, $U$ -- evaluation representation
corresponding to the representation $U$ of $\sln$ described above we
get an affine analogue of the theory of Jacobi polynomials. This is
closely related with the elliptic analogue of the Sutherland
operator. We can
generalize to this case some results from the theory of usual Jacobi
polynomials; there are also arise  new phenomena, such as modular
properties of these polynomials (as functions of the modular parameter
of the torus). 


\vfill\newpage

\vbox{\vskip 1in}
\specialhead\chapter{1} 
\centerline{BASIC DEFINITIONS} 
\endspecialhead

\vbox{\vskip 0.5in}

\head 1.1 Simple Lie algebras and their representations\endhead

In this whole work, all the objects are always defined over the ground
field $\C$ of complex numbers or its extensions.
 In this section we briefly list the
main facts on the simple Lie algebras and their representations which
we are going to use; all of these facts are quite standard and can be
found in any textbook on Lie algebras (see, for example, \cite{Hu, B}). 

Let $\g$ be a simple Lie algebra over $\C$ of rank $r$, $\h\subset \g$
-- its Cartan subalgebra, $U \g$ -- its universal enveloping algebra.
Then $\g$ has the root decomposition: 

$$\g=\h\oplus\bigoplus_{\a\in
R}\g_\a,$$ 
where $R\subset \h^*$ is the corresponding root system. We
fix a polarization of $R$: $R=R^+\sqcup-R^+$, where $R^+$ is the
subset of positive roots. We denote by $\alpha_1,\ldots, \alpha_r\in
R^+$ the basis of simple roots.  The polarization of
roots gives rise to polarization of $\g: \g=\frak n^-\oplus \h\oplus
\frak n^+$, where $\frak n^\pm =\bigoplus_{\a\in R^\pm} \g_\a$, and
corresponding polarization of $U\g: U\g=U\frak n^-\o U\h\o
U\frak n^+$.

We fix an invariant symmetric
bilinear form $(\, , \,)$ on $\g$ by the condition that for the associated
bilinear form on $\h^*$ we have $(\a, \a)=2$ for long roots;
 this form allows us
to identify $\h^*\simeq \h: \lambda\mapsto h_\lambda$. Abusing the
language, we will also use the
notation $(\, , \,)$ for  the associated bilinear form on $\h^*$. 
We denote by $\<\ , \ \>$ the canonical pairing $\h\o
\h^*\to \C$.

As usual,
for every $\a\in R$ we define the dual root $\a\v
=\frac{2h_\a}{(\a,\a)}\in \h$ and  introduce the following notions: 

$Q=\bigoplus \Z \a_i$ -- root lattice;

$Q^+=\bigoplus \Z_+\a_i$; 

$Q\v=\bigoplus \Z \a_i\v\subset \h$ --  coroot lattice;

$P=\{\l\in \h^*|\<\l, \a\v\>\in \Z\}$ -- weight lattice; 

$P^+=\{\l\in \h^*| \<\l, \a_i\v\>\in \Z_+\}$ -- 
cone of integer dominant weights;

$\omega_i$ -- fundamental weights: $\<\omega_i, \a\v_j\>=\delta_{ij}$;

$P\v=\{ h\in \h|\<h, \a\>\in \Z\}$ -- coweight lattice;
    
$\rho=\half \sum_{\a\in R^+}\a=\sum\omega_i$;

$\theta\in R$ -- the highest root: $\theta-\a\in Q^+$ for every $\a\in
R$;  

$h\v=\<\rho, \theta\v\>+1$ -- dual Coxeter number for $\g$;

We have a (partial) order on $P$ defined as follows: $\l\le
\mu$ if $\mu-\l\in Q^+$.

    We denote by $\C[P]$ the group algebra of the abelian
group $P$, i.e. the algebra over $\C$ spanned by formal
exponents $e^\l, \l\in P$ with relations $e^0=1,
e^{\l+\mu}=e^\l e^\mu$.  If an element $f\in \C[P]$ can be
written in the form $f=\sum_{\mu\le \l} a_\mu e^\mu, a_\l\ne
0$ then we say that $a_\l e^\l$ is the highest term of $f$
and write $f=a_\l e^\l+\lot$.

As usual, we denote by $W$ the Weyl group of $R$ and by $l(w)$ length
of an element $w\in W$ with respect to the generators $s_i=s_{\a_i}$.
This group acts naturally on $\h^*$, preserving $R$;
thus, it acts on $P$ and  $\C[P]$. Of special interest for us
will be the algebra $\C[P]^W$ of $W$-invariant elements in
$\C[P]$. We will often use the fact that the orbitsums 

    $$m_\l=\sum_{\mu\in W\l} e^\mu,\quad \l \in P^+\tag 1.1.1$$
form a basis in $\C[P]^W$, which follows from the fact that
for every $\l\in P$ the orbit $W\l$ contains precisely one
element from $P^+$.

We 
choose a basis $e_\a\in \g_\a, f_\a\in \g_{-\a}$ for $\a\in R^+$ such
that $(e_\a, f_\a)=1$; then $[e_\a, f_\a]=h_\a$. 
In particular, if $\a=\a_i$ is a simple root then the
elements 
$$\gathered
e_i=\frac{2}{(\a_i,\a_i)}e_{\a_i},\quad  
f_i=\frac{2}{(\a_i,\a_i)}f_{\a_i},\\
h_i=\a_i\v=\frac{2}{(\a_i,\a_i)}h_{\a_i}\endgathered$$
 satisfy the usual relations of
the Lie algebra $\sltwo$. Moreover, $e_i, f_i, h_i,
i=1\dots r$ generate $\g$.

Let $a_i$ be any orthonormal basis in $\g$. Define the Casimir element
$C\in U\g$ by $C=\sum a_i^2$. This element is central and does not
depend on the choice of orthonormal basis. Therefore, it is easy to
check that it can be written as follows: let $x_l,
l=1\ldots r$ be an orthonormal basis in $\h$; then 
 $C=\sum_{\a\in R^+} e_\a f_\a +f_\a e_\a+\sum x_l^2$. 

We will also use the following 
useful identity: 

$$\Delta(C)=C\o 1 + 1\o C + 2\Omega,\tag 1.1.2$$
where 
$$\Omega=\sum a_i\o a_i 
= \sum_{\a\in R^+} e_\a \o f_\a +f_\a\o  e_\a+\sum x_l\o x_l\tag 1.1.3$$ 
is the $\g$-invariant element in $\g\o \g$, and $\Delta:U\g \to
U\g\o U\g$ is the 
comultiplication: $\Delta(x)=x\o 1+ 1\o x, x\in \g$ and
$\Delta(ab)=\Delta(a) \Delta(b)$.

    In general, there is no explicit construction of the
whole center of $U\g$.  However, its structure as a graded algebra
(and more over, as a module over $U\g$) is known.

\proclaim{Theorem 1.1.1} 

$$\Cal Z (U\g)\simeq (S\h)^W.$$
\endproclaim
We will construct the isomorphism (Harish-Chandra isomorphism) later;
now we only say that under this isomorphism the Casimir element
$C\mapsto p_2: p_2(\l)=(\l, \l)-(\rho, \rho)$. Note that since it is
known that $(S\h)^W$ is a free polynomial algebra, the same is true
for $\Cal Z(U\g)$.

Also, there exists a unique involutive  algebra
automorphism (Chevalley involution) $\omega:\g\to\g$ such that 
$$\gathered
\omega(e_i)=-f_i,\\
 \omega(f_i)=-e_i, \\
\omega (h)=-h, \quad h\in \h.\endgathered\tag 1.1.4$$

Let us now consider the representation theory of $\g$. If $V$ is a
$\g$-module then we denote an action of $x\in \g$ in $V$ by
$\pi_V(x)$. 
We will always
consider modules with a weight decomposition: $V=\bigoplus_{\l\in
\h^*} V[\l]$, where $V[\l]$ are finite-dimensional spaces such that 
$h v=\<\l, h\>v$ if $v\in V[\l], h\in \h$. We will
call $v\in V[\l]$ vectors of weight $\l$. 
For every $\l\in
\h^*$ we define the Verma module $M_\l$ by the following conditions: 

1. $M_\l$ is spanned over $U\g$ by a single  vector $v_\l$ (highest weight
vector) such that $e_\a v_\l=0$, $v_\l$ has weight $\l$. 

2. $M_\l$ is a free module over $U\frak n^-$.

    These modules are simplest examples of modules from
category $\Cal O$. By definition, a $\g$-module $V$ is said
to be from category $\Cal O$, if it has weight
decomposition, is finitely generated and satisfies the
following condition: for every $v\in V$, the space $U\frak
n^+ v $ is finite-dimensional.

Let $I_\l$ is the maximal ideal in $M_\l$ which does not
contain $v_\l$ (it exists and is unique); denote by 
    $L_\l=M_\l/I_\l$ the corresponding irreducible highest-weight
module. For generic $\l$ (that is,
for all $\l$ except a subvariety of codimension 1), $M_\l=L_\l$.
Moreover, it is known that there is a unique bilinear form on $M_\l$
(Shapovalov form) such that $(v_\l, v_\l)=1$ and $(v_1, x v_2)= -
(\omega(x)v_1, v_2)$ for all $x\in \g$. This form is symmetric and its
kernel is $I_\l$; thus, it descends to a non-degenerate form on $L_\l$.

It is known that $L_\l$ is finite-dimensional iff $\l\in P^+$, and
every finite-dimensional irreducible representation of $\g$ is
obtained in this way. Moreover, if $\l\in P^+$ then the dual
representation $L_\l^*$ is also finite-dimensional irreducible:
$L_\l^*= L_{\l^*}$, where 

    $$\l^*=-w_0(\l),$$
 $w_0$ being the longest element in the Weyl group. 

If $c\in \Cal Z(U\g)$ is a central element, then
$c|_{M_\l}=\chi(c)(\l+\rho)\Id_{M_\l}$, where $\chi:\Cal Z(U\g)\to
(S\h)^W$ is the Harish-Chandra isomorphism mentioned before. 
The same holds for any subquotient of $M_\l$, in
particular, for $L_\l$. Note that 
for the Casimir element we have $C|_{M_\l}=(\l,
\l+2\rho)\Id_{M_\l}$; thus, in the adjoint representation $C$
acts by multiplication by $2h\v$.

For a  finite-dimensional representation $V$  of $\g$ let 
$\ch V= \sum_{\l\in \h^*} \dim V[\l] e^\l\in \C[P]$. For an
irreducible representation, the character is given by the Weyl
formula: 

$$\ch L_\l = \frac{\sum_{w\in W} (-1)^{l(w)}
e^{w(\l+\rho)}}{\delta},\quad \l\in P^+,\tag 1.1.5$$
where 
$$\delta=\prod_{\a\in R^+}(e^{\a/2}-e^{-\a/2})
=e^\rho\prod_{\a\in R^+} (1-e^{-\a})
    \tag 1.1.6$$
is the Weyl denominator. Obviously, $\ch L_\l=e^\l+\lot$.

    More generally, if $V$ is a module from category $\Cal
O$ then we can define its character by the same formula 
$\ch V= \sum_{\l\in \h^*} \dim V[\l] e^\l$. However, in
this case it lies not in $\C[P]$ but in its completion
$\overline{\C[P]}$, which is defined as follows: 

    $$\overline{\C[P]}=\{\sum_{n=0}^{\infty} a_n e^{\l_n}|
\lim (\rho, \l_n)=+\infty\}.\tag 1.1.7$$

    In the following we will be especially interested in
the case $\g=\sln$, i.e. the root system of type $A_{n-1}$.
In this case the root system and related objects admit the
following explicit realization: 

    $\h^*=\{(\l_1, \dots, \l_n)|\sum \l_i=0\}\subset \C^n$;

    $R=\{\varepsilon_i-\varepsilon_j\}_ {i\ne j},$
    where $\varepsilon_i$ is the standard basis in $\C^n$: 
$\varepsilon_i=(0,\ldots,0,1,0,\ldots, 0)\in \Z^n $ (1 in the $i$-th
place);

 $(\l, \mu)=\sum \l_i\mu_i $;
    
$R^+=\{\varepsilon_i-\varepsilon_j\}_ {i< j},$;

$\alpha_i=\varepsilon_i-\varepsilon_{i+1}$;

$Q=Q\v=\{(\l_1, \dots, \l_n)\subset \Z^n|\sum \l_i=0\}$;

    $P=\{(\l_1, \dots, \l_n)\subset\frac{1}{n}\Z^n |\sum
\l_i=0, \l_i-\l_j\in \Z\}$;    

    $P^+=\{(\l_1, \dots, \l_n)\subset\frac{1}{n}\Z^n |\sum
\l_i=0, \l_i-\l_{i+1}\in \Z_+\}$;    

    $\rho= (\frac{n-1}{2}, \frac{n-3}{2}, \dots, \frac{1-n}{2})$;

$W=S_n$;

    $\C[P]\simeq \{\text{homogeneous polynomials in
}x_1^{\pm 1}, \dots, x_n^{\pm 1} \text{ of degree zero}\}$
(this isomorphism is given by letting $e^{\varepsilon_i}=x_i$).        

Note also in this case $\l, \mu\in P, \l < \mu$ with respect to the
order on $P$ we defined before implies that $\l \prec \mu$, where
$\prec$ is the lexicographic order on $\R^n$. 

    Many results for the Lie algebra $\sln$ are also valid
for (reductive) Lie algebra $\frak{gl}_n$ with suitable
changes. Namely, the group algebra $\C[P]$ should be 
replaced by the algebra $\C[x_1, \dots, x_n]$, $P^+$ should
be replaced by the set of all partitions in $n$ parts, i.e.
sequences $(\l_1, \dots, \l_n)$ such that $\l_i\in Z_+,
\l_1\ge \l_2\ge \dots\ge \l_n\ge 0$, $\rho$ should be
replaced by $(n-1, n-2, \dots, 0)$ and Weyl denominator
should be replaced by Vandermonde determinant
$\prod_{i<j}(x_i-x_j)$. 
    
\head 1.2 Quantum groups \endhead

In this section we briefly list the main definitions from the theory
of quantum groups (see \cite{Dr1, J1,J2}).
 In this section all the objects are considered over
the ground field $\C_q=\C(q^{1/2N})$, where $q^{1/2N}$ is a formal
variable, and $N=|P/Q|$ (this will be justified below).

Let $R, R^+,..$ be as defined in Section~1.1. Let us renormalize the
inner product in $\h$, introducing $(x,y)'=m(x,y)$, where $m$ is
chosen so that $(\a, \a)'=2$ for {\bf short} roots. Thus, $m$ can 
take values 1 (for simply-laced root systems), 2 (for root
systems of type $B, C, F$) or 3 (for $G_2$). This renormalization is
quite standard in the theory of quantum groups. In this case 
$d_i=\frac{(\a_i, \a_i)'}{2}$ are positive integers 
 such  that  $d_i a_{ij}=d_j a_{ji}$, 
 where $a_{ij}=2(\a_i, \a_j)/(\a_i, \a_i)$ is the Cartan
matrix, and $\text{g.c.d.} (d_i)=1$.

    In this section  we will identify $\h$ and $\h^*$ using the
form $(\, ,\,)'$; (this differs by a factor of $m$ from the
identification of the previous chapter!); then
$\a_i=d_i\a_i\v$, and therefore $Q\subset Q\v, P\subset
P\v$. Since, by definition, for every $\l\in P$ we have
$N\l\in Q$, this implies 
$(\l, \mu)'\in 
\frac{1}{N} \Z$ for every $\l\in P,\mu\in P$, which is
the reason why have chosen the ground field to be
$\C(q^{1/2N})$. 

Define the quantum group $\Ug$ as an associative algebra over $\C_q$
with generators $e_i, f_i, q^{h},h\in \frac{1}{2}P\v\subset \h, 
 i=1\ldots r$ and relations

$$\gathered
q^0=1, \quad q^{a+b} =q^{a}q^b,\\
q^h f_j q^{-h}=q^{-\<h, \a_j\>} f_j,\\
q^h e_j q^{-h}=q^{\<h, \a_j\>} e_j,\\
[e_i, f_j]=\delta_{ij}\frac{q^{d_i h_i}-q^{-d_i h_i}}{q^{d_i}-q^{-d_i}},
\endgathered \tag 1.2.1$$

$$\gathered 
\sum_{n=0}^{1-a_{ij}} 
	\frac{(-1)^n}{[n]_i! [1-a_{ij}-n]_i!}
	e_i^n e_j e_i^{1-a_{ij}-n}=0,\\
\sum_{n=0}^{1-a_{ij}} 
	\frac{(-1)^n}{[n]_i! [1-a_{ij}-n]_i!}
	f_i^n f_j f_i^{1-a_{ij}-n}=0,\endgathered
\tag 1.2.2$$

where, as before,  $h_i=\a_i\v$ and 
$$[n]_i=\frac{q^{nd_i} - q^{-nd_i}}{q^{d_i}-q^{-d_i}}, 
\quad [n]_i! =[1]_i \dots [n]_i.\tag 1.2.3$$

This is a Hopf algebra with the following comultiplication, counit and
antipode: 
$$\gathered 
\Delta e_i =e_i\o q^{-d_i h_i/2} +q^{d_i h_i/2}\o e_i,\ 
\Delta f_i =f_i\o q^{-d_i h_i/2} +q^{d_i h_i/2}\o f_i,\\
\Delta q^h=q^h\o q^h,\\
\epsilon(q^{h})=1, \epsilon(e_i)=\epsilon(f_i)=0,
S(e_i)=-q^{-d_i}e_i, S(f_i)=-q^{d_i}f_i, S(q^h)=q^{-h}.
\endgathered\tag 1.2.4$$

    Often in the literature a smaller algebra is
considered, which only contains $q^h$ for $h\in \half Q\v$,
in which case the definition is given in terms of
generators $K_i=q^{h_i}$. Also, one can get  many different forms
by replacing our $e_i$ by $q^{\pm d_i h_i/2} e_i$ or
replacing $q$ by $q^{-1}$; however, essentially all  these
definitions are equivalent (as long as they are over
rational functions of $q$ and not over formal power
series).

Many of the facts about $U\g$ can be generalized to $\Ug$ without
difficulties. For example, $\Ug$ has a polarization: $\Ug
=U^+\cdot U^0\cdot U^-$. It also admits a ($\C_q$-linear) 
Chevalley involution $\omega$ 
which is defined by $e_i\mapsto -f_i, f_i\mapsto -e_i, q^h\mapsto
q^{-h}$;
 this is an algebra automorphism and coalgebra antiautomorphism. We
will also need the following property of this involution:
$S\omega=\omega S^{-1}$.  

Similarly, the representation theory of $\Ug$ is quite parallel to the
classical case (see \cite{L1}), the only difference being
that all the representations 
must be considered as linear spaces over $\C_q$. Namely, for every
$\l\in \h^*$ we can define Verma module $M_\l$ in a way similar to that
of Section~1.1. If $\l$ is not an integer weight, then this module
should be considered as a linear space over the field $\C_{q,q^\l}$ of
rational functions in $q^{1/2N}, t_i=q^{\<\l, \a_i\v\>/2N}$. Again, each
$M_\l$ has a unique (non-trivial) irreducible quotient $L_\l$, and
$L_\l$ is finite-dimensional iff $\l\in P^+$. Verma modules
and all its subfactors (in particular, $L_\l$) have weight
decomposition, and character (and therefore, dimension) of
$L_\l$ is the same as in the classical case. More
generally, the same is true for the whole category $\Cal
O$: every $\g$-module from category $\Cal O$ has a
$q$-analogue.

    On each $M_\l$ one
can define Shapovalov form $(\ ,\ )$ such that $(xv,v')=(v,\omega
S(x)v')$; this form is defined uniquely up to a factor, it is
symmetric and its
restriction to $L_\l$ is non-degenerate. 

    Though it is not literally true
that every finite-dimensional irreducible module is of the form
$L_\l$, it is almost so: if we restrict ourselves to consideration of
the modules with weight decomposition: $V=\bigoplus V[\l],
q^h|_{V[\l]}=q^{\<h, \l\>}\Id_{V[\l]}$
then the category of such  finite-dimensional
representations of $\Ug$ is equivalent to that of representations of
$\g$ as a tensor category;  in particular, 
every finite-dimensional representation is
completely reducible, the only irreducible
finite-dimensional representations are $L_\l, \l \in P^+$, 
and $L_\l\o L_\mu=\sum N_{\l\mu}^\nu L_\nu$, where the
multiplicities $N_{\l\mu}^\nu$ are the same as in classical case. 
In this whole work we only consider 
representations with weight decomposition. 

However, there are  several very important differences between the
representations of $\Ug$ and of $\g$. First of all, we have a
different notion of dual representation. Namely, 
for every finite-dimensional representation $V$ of $\Ug$ we define the
action of $\Ug$ on the space $V^*$ of linear functionals on $V$ by the
rule $\<xv^*, v\>=\<v^*, S(x)v\>$ for $v\in V, v^*\in V^*, x\in
\Ug$, where $S$ is the antipode in $\Ug$. This endows $V^*$ with the
structure of $\Ug$ representation 
which we will call right dual to $V$.  In a similar way, the left dual
${}^*V$ is the representation of $\Ug$ in the space of linear
functionals on $V$ defined by $\<xv^*, v\>=\<v^*, S^{-1}(x)v\>$. Then
the following natural pairings and embeddings are $\Ug$-homomorphisms
(warning: the order is important here!):

$$\gathered V\o {}^*V\to \C_q, \quad V^*\o V\to\C_q\\
\C_q\to V\o V^*, \quad \C_q\to {}^*V\o V\endgathered\tag 1.2.5$$

Note that $V^*$ and ${}^*V$, considered as two structures of
representation of $\Ug$ on the same vector space, don't
coincide, but they are isomorphic. Namely, $q^{-2\rho}:{}^*V\to V^*$
is an isomorphism of $\Ug$-modules (recall that we
identified $\h$ and $\h^*$ using the form $(\, , \,)'$, so
$2\rho\in Q\subset Q\v$). This follows from the relation
$S^2(x)= q^{-2\rho} x q^{2\rho}$, which can be easily
checked (it suffices to check it for  the generators $e_i, f_i$).

If $V=L_\l$ is an irreducible
finite-dimensional representation, then so is $V^*$: $L_\l^*\simeq
L_{\l^*}$. 

Second, unlike the classical case, if $V,W$ are finite-dimensional 
representations of $\Ug$ then $V\o W$ and
$W\o V$ are isomorphic, but the isomorphism is non-trivial. More
precisely (see \cite{Dr1}), there exists a universal R-matrix
$\Cal R$ which is an element of a certain completion of 
$ \Ug\o \Ug$ 
such that
$$\check R_{V,W}= P\circ \pi_V\o \pi_W(\Cal R) \colon V\o W\to
W\o V
\tag 1.2.6$$
is an isomorphism of representations. Here $P$ is the transposition:
$Pv\o w=w\o v$.
Also, it is known that $\Cal R$ has the following
form:
$$\gather
\Cal R=q^{-\sum a_i\o a_i}\Cal R^*, \quad 
	\Cal R^*\in U^+\hat \o  U^-\tag 1.2.7\\
(\epsilon\o 1)(\Cal R^*)=(1\o\epsilon)(\Cal R^*)=1\o
1,\endgather$$
where $a_i$ is an orthonormal basis in $\h$ with respect to $(\, ,\,)'$. 
As we said, $\Cal R$ does not lie in the tensor square of
$\Ug$ but in its certain completion; however, for any pair
of finite-dimensional representations $V, W$ of $\Ug$ the
operator $\pi_V \o \pi_W(\Cal R)$ is well-defined (this is
where we need fractional powers of $q$ in the definition of
$\C_q$).

Finally, if $V$ is a representation of $\Ug$ let us
consider the action of
$\Ug$ in $V$ given by $\pi_{V^{\omega}}(x)=\pi(\omega x)$, where $\omega$
is the Chevalley involution defined above. We denote $V$ endowed with this
action by $V^\omega$. One can easily check that if $V$ is
finite-dimensional then $V^\omega\simeq {}^*V$
(which is, of course, equivalent to saying that $V^\omega\simeq V^*$):
if $(\ ,\ )$ is Shapovalov form on $V$ then $x\o y\mapsto (x,y)$
is a homomorphism of $\Ug$-modules $V\o V^\omega\to \C_q$.
Note also that $(V\o
W)^\omega=W^\omega\o V^\omega$ and that if $\Phi\colon V\to W$ is an
intertwiner then $\Phi$ is also an intertwiner considered as a map
$V^\omega\to W^\omega$.

    In certain instances, the structure of $\Ug$ is even
simpler than that of $U\g$. A good example is the following
description of the center of $\Ug$, due to Drinfeld
(\cite{Dr2}) (a similar
construction was independently proposed by N.~Reshetikhin \cite{R}).

    \proclaim{Theorem 1.2.1} 
\roster\item Let $V$ be a finite-dimensional
representation of $\Ug$. Define 

    $$c_V=(\text{\rm Id} \otimes \Tr_{V})\left(\Cal R^{21}\Cal R
(1\otimes q^{-2\rho})\right),\tag 1.2.8$$
where $\Cal R^{21}= P(\Cal R), P(x\o y)=y\o x$. Then
$c_V$ is a central element in $\Ug$.

    \item $$ c_V|_{M_\l}=\ch V (q^{-2(\l+\rho)})\Id_{M_\l},$$
 where $f(q^{2\l})$ for $f\in \C[P]$
denotes the polynomial in $q^{1/2N}$ obtained by replacing  
each $e^\mu$ by $q^{2(\l, \mu)'}$.
    
    \item The elements $c_{L_\l},\l\in P^+$ form a basis in
$\Cal Z(\Ug)$.

\endroster\endproclaim

    \demo{Proof} 
(1) This is based on the following statement (see \cite{Dr2}):
if $\theta\colon \Ug\to \C_q$ is such that $\theta(xy)=\theta(y S^2(x))$
then the element $c_\theta=(\Id \otimes \theta)(\Cal R^{21}\Cal R)$ is
central. On the other hand, we know that $S^2(x)=q^{-2\rho }x q^{2\rho}$, so
$\theta(x)=\Tr|_V(xq^{-2\rho} )$, where $V$ is any finite-dimensional
representation of $\Ug$, satisfies $\theta(xy)=\theta(y
S^2(x))$.

(2) Let $v_\l$ be a highest-weight vector in $M_\l$; let us
calculate $c_V v_\l$. Let $w\in V[\mu]$. Then (1.2.7) implies

$$\Cal R^{21}\Cal R (v_\l\otimes
w)=q^{-2(\l,\mu)'}v_\l\otimes w + \sum v'_i\otimes w'_i,$$
where $\text{wt }w'_i <\mu$. Thus, $c_r v_\l =(\sum \limits_\mu
(\text{dim }V[\mu]) q^{-2(\l,\mu)'}q^{-2(\rho,\mu)'})
v_\l$,  where the sum is
taken over all the weights of $V$.

    (3) This is based on the quantum version of
Harish-Chandra isomorphism. Name\-ly, 
if $C$ is a central element then it acts in Verma
module $M_\l$ by some constant; denote it $\gamma(C)
(\l+\rho)$. It is easy to see that $\gamma(C)(\l+\rho)$ can
be uniquely written in the form $f(q^{\l+\rho})$ for
some $f\in \C_q[\frac{1}{2}P\v]^W$. Moreover, it is known  (see
\cite{T})
that  $\gamma$ is an isomorphism  $\Cal Z(\Ug)\simeq\C_q[2P]^W$.
 Since  characters of irreducible
representations form a basis in $\C_q[P]^W$, we obtain the
statement (3).
\qed\enddemo
\remark{Remark} Note that though $\Cal R$ is defined only
in some completion of $\Ug^{\o 2}$, the central elements we
constructed belong to $\Ug$ itself.\endremark


\vfill\newpage

\vbox{\vskip 1in}
\specialhead\chapter{2} 
\centerline{GENERALIZED CHARACTERS} 
\endspecialhead

\vbox{\vskip 0.5in}

This chapter is the heart of the whole dissertation. In this chapter
we define for every $\g$-intertwining operator $\Phi$  a certain
weighted trace $\chi_\Phi$, which we call {\it generalized character},
which is justified by the fact that if $\Phi:V\to V$ is the identity
operator then $\chi_\Phi$ is the usual character of $V$. We also study
basic properties of these characters: namely, we prove an
orthogonality theorem for these characters, and show that they are
eigenfunctions of a large family of commuting differential operators. 
We also describe the $q$-analogue of these generalized characters,
which is defined in a similar way, but $\g$ must be replaced with the
quantum group $\Ug$. 

All of constructions of this chapter belong to the author, P.~Etingof
and in some part I.~Frenkel, see \cite{EK1, EK2, EK4, EFK}.

    \head 2.1 Weighted traces of intertwiners and
orthogonality theorem 
\endhead 

In this section  we define 
 the basic object of our study:
the weighted traces of intertwining operators, which we
call generalized characters. {\bf All the definitions and
theorems of this section apply equally to the quantum group
$\Ug$ and the  usual universal enveloping algebra $U\g$.} To avoid
unnecessary repetitions, we formulate all the results for
$\Ug$; however, the reader should keep in mind that all of
the statements also hold if one replaces $\Ug$ by $U\g$,
and $\C_q$ by $\C$.

\definition{Definition 2.1.1} 
Let $V$ be a $\Ug$-module from category $\Cal O$, 
$U$ -- an arbitrary $\Ug$-module with weight
decomposition (not necessarily finite-dimensional), and let
$\Phi\colon V\to V\otimes U$ be a non-zero intertwining
operator. 
 
A generalized character $\chi_\Phi\in
\overline{\C[P]}\otimes U$ is defined by 

$$\chi_\Phi=\sum_{\l\in P} e^\l \Tr |_{V[\l]} (\Phi).\tag 2.1.1$$

    In particular, if $V$ is finite-dimensional then
$\chi_\Phi\in \C[P]\o U$.
    \enddefinition

    \example{Example} Let $U=\C$ be the trivial representation
of $\Ug$, and let $\Phi: V\to V$ be the identity operator. Then
$\chi_{\Id}$ is the usual character of $V$, which
justifies the name ``generalized character''.\endexample

For finite-dimensional $V$, 
we can consider these traces as functions on $\h$ with
values in $U$ by the rule $e^\l(h)=e^{2\pi\i\<\l, h\>}$. This is equivalent
to letting 
$$\chi_\Phi(h)=\Tr _V (\Phi e^{2\pi\i h}),\tag 2.1.2 $$
where $e^{2\pi\i h}$ is an operator in $V$ defined by 
$e^{2\pi\i h}|_{V[\l]}=e^{2\pi\i\<h,\l\>}\Id_{V[\l]}$. 

    Using the fact that $\Phi$ preserves weight, it is easy
to see that in fact $\chi_\Phi$ takes values in the
zero-weight subspace $U[0]$. Also, 
for $V=L_\l, \l\in P^+$  the highest term  of $\chi_\Phi$ is
$ue^\l$ and the lowest term is $u'e^{-\l^*}$ for some
$u,u'\in U$.

    Note that in the quantum case (i.e., for $\Ug$ rather
than $U\g$)  $\chi_\Phi$ has  no Weyl group symmetry 
if $U$ is a nontrivial representation (see a
counterexample in Theorem~4.2.2).

    Let us recall the well-known results on the
orthogonality of (usual) characters. Let us introduce the
($q$-linear) bar involution in $\C_q[P]$ by
$\overline{e^\l}=e^{-\l}$,  and let $[\ \ ]_0: \C[P]\to \C$ 
be the constant term: $[\sum c_\l e^\l]_0=c_0$. Then it is
known that characters of irreducible finite-dimensional  
representations of $\Ug$ (which are the same as for  $\g$)
 are orthonormal with respect to
the following inner product in $\C[P]$: 

    $$\<f, g\>_1=\frac{1}{|W|}[f\bar g\Delta]_0,\tag 2.1.3$$
where 

    $$\Delta = \delta\bar \delta= 
\prod_{\alpha\in R}(1-e^\alpha).\tag 2.1.4$$
    
    It turns out that a similar statement holds for
generalized characters:

    \proclaim{Theorem 2.1.2}{\rm 
(Orthogonality theorem for generalized characters)}
Let $\l, \mu\in P^+$, and let $U$ be a finite-dimensional
representation of $\Ug$. Let  
$\Phi_\l:L_\l\to L_\l\otimes U, \Phi_\mu\colon L_\mu\to
L_\mu\otimes U$ be $\Ug$-intertwiners, and $\l\ne\mu$. Then 
the generalized characters
$\chi_1=\chi_{\Phi_\l}, \chi_2=\chi_{\Phi_\mu}$ are orthogonal
with respect  to the following  inner
product: $\<f,g\>_1=\frac{1}{|W|}[(f,\bar g)_U \Delta]_0$, 
where
$(\cdot,\cdot)_U$ is the Shapovalov form in $U$   and all the other
notations are as before.
\endproclaim
\demo{Proof}
As was explained above, we can as
well consider $\Phi_\mu$ as an intertwiner $L_\mu^\omega\to
U^\omega\o  L_\mu^\omega$. This implies that 
$(\chi_1,\overline{\chi_2})_U=\chi_{\Psi}$
 (note the change of sign of $h$ in the second
factor!),  where the intertwiner $\Psi\colon
L_\l\o  L_\mu^\omega\to L_\l\o  L_\mu^\omega$ is
defined as the following composition

$$L_\l\o  L_\mu^\omega @>{\Phi_\l\o  \Phi_\mu^\omega}>>
 L_\l\o  U\o 
U^\omega\o  L_\mu^\omega @>{\text{Id}\o 
(\cdot,\cdot)_U\o  \text{Id}}>>
L_\l\o  L_\mu^\omega.\tag 2.1.5$$

Since $L_\l\o  L_\mu^\omega=\bigoplus N_\nu L_\nu$, we see that
$(\chi_1,\overline{\chi_2})_U$ is a linear combination
of usual characters $\ch {L_\nu}$.
But since these characters are the same as
for $\g$, we know that $\frac{1}{|W|}[\chi_\nu\Delta]_0=\delta_{\nu,0}$.
On the other hand, it is known that if $\l\ne \mu$ then the
decomposition of $L_\l\o  L_\mu^\omega$ does not contain
the trivial representation (i.e. $N_0=0$); thus, in this case $\chi_1$
and $\chi_2$ are orthogonal.
\qed
\enddemo

    \remark{Remark} Note that the
calculation of the norm $\<\chi_\Phi, \chi_\Phi\>_1$ (which is trivial
for $U=\C$) in general case is very
complicated, since it requires calculation of restriction of the
intertwiner $\Psi$ defined above to the trivial subrepresentation in
$L_\l\o L_\l^\omega$. In some simplest cases (where the space of
intertwiners is one-dimensional) it can be calculated explicitly (we
will do it in Chapter~5), and the answer is quite non-trivial; it
coincides with so-called Macdonald inner product identities for root
system $A_n$, see \cite{M2}.
\endremark

    \head 2.2 Center of $U\g$ and commuting differential
operators\endhead 

    In this section we consider some properties of
generalized characters for $U\g$, following the paper of
Etingof \cite{E1}. In this section, the words
``representation'' and ``intertwiner'' stand for
representation of $\g$ and $\g$-intertwiner. In this
section we always consider the case when $V, U$ are
finite-dimensional, though in fact all the results can be
easily generalized (see remark in the end of this section).

    Let $\C[P](e^\a-1)^{-1}$ denote the ring of fractions
$\frac{f}{g}$, where $f,g\in \C[P]$ and $g$ is a product of
linear terms of the form $e^\a-1, \a\in R$. 
Let us introduce the following
ring of differential operators:

    $$DO=\{\text{differential operators 
    on $\h$ with coefficients from } \C[P](e^\a-1)^{-1}\}.$$ 

    These operators can be also considered formally, as
derivations of the ring of fractions $\C[P](e^\a-1)^{-1}$. In
particular, we will use the following differential
operators: 

    $\Delta_\h$ -- Laplace operator in $\h$ normalized so that
$\Delta_\h e^\l= (\l, \l)e^\l$; 

    $\d_\a, \a\in \h^*$ -- differentiation along the vector
$\frac{1}{2\pi\i} h_\a$, so that $\d_\a e^\l= (\a,
\l)e^\l$. 
\proclaim{Theorem 2.2.1} Let us denote by $I$ the left
ideal in $U\g$ generated by $\h$. Then 
 for any $u\in U\g$ there exists a unique  differential
operator $D_u\in \text{DO}\o  U\g/I$, 
such that for any intertwiner $\Phi:V\to V\o U$ 
we have 

$$\Tr_V(\Phi u e^{2\pi\i h})=D_u\Tr_V(\Phi e^{2\pi\i h}).\tag 2.2.1$$

    Note that  $D_u f$ is well defined for any function
$f$ with values in $U[0]$ .

\endproclaim

    We do not give the proof of this theorem (which can be
found in \cite{E1}) here, since the proof is quite similar
to the proof in the $q$-case, which we give later (see
Theorem~2.3.1). The main idea of the proof will be  quite clear
from the example of Casimir element, which we discuss
below. 

    \proclaim{Corollary 2.2.2} Let $c$ be a central element:
$c\in \Cal Z(U\g)$, $V, U$ be finite-dimensional
representations  of $\g$ such that $c$ acts in $V$ by
multiplication by a constant $c(V)$, and $\Phi: V\to V\o U$
be a non-zero intertwiner. Then the generalized character
$\chi_\Phi$ satisfies the following differential equation 

    $$D_c \chi_\Phi = c(V) \chi_\Phi, \tag 2.2.2$$
where the differential operator $D_c$ is defined in
Theorem~\rom{2.2.1}. 
\endproclaim

    \demo{Proof} Let us consider $f(h)=\Tr_V (\Phi c e^{2\pi\i h})$.
On one hand, since $c$ acts by multiplication by $c(V)$, we
see that $f(h)=c(V) \chi_\Phi(h)$. On the other hand, by
Theorem~2.2.1, $f(h)=D_c \chi_\Phi$. \qed\enddemo

    In general, $u\mapsto D_u$ is not an algebra
homomorphism. However, it becomes a homomorphism when
restricted to the center:

    \proclaim{Theorem 2.2.3} Let us keep the notations of
Theorem~\rom{2.2.1}.     Let $c,c'\in \Cal Z(U\g)$. Then $D_{cc'}=D_c D_{c'}$;
thus, $c\mapsto D_c$ is a homomorphism $\Cal Z(U\g)\to
\text{DO}\o U\g/I$.
    \endproclaim

    \demo{Proof} By definition, $D_{cc'}\chi_\Phi=\Tr (\Phi c c'
e^{2\pi\i h})$. On the other hand, since $c$ is central, $\Phi c$ is
also an intertwiner, and therefore we can write $\Tr (\Phi
c c' e^{2\pi\i h}) = D_{c'} \chi_{\Phi c} = D_{c'}D_c \chi_\Phi$. \enddemo

Let us  apply the construction of Theorem~2.2.1  to the
Casimir element $C$ defined in Section~1.1. In this case
the answer can be written down explicitly: 

    \proclaim{Proposition 2.2.4} 
\roster\item
    Let $C$ be the Casimir
element. Then 
$$\multline
 D_C=\Delta_\h -2\sum_{\a\in R^+} \frac{e_\a
f_\a }{(e^{\a/2}-e^{-\a/2})^2}-\sum_{\a\in R^+}
\frac{1+e^\a}{1-e^\a}\d_\a \\
=\delta^{-1}\circ
\left(\Delta_\h-2\sum_{\a\in R^+} \frac{e_\a
f_\a }{(e^{\a/2}-e^{-\a/2})^2}-(\rho, \rho)\right)
\circ \delta,\endmultline \tag 2.2.3$$
where $\delta$ is Weyl denominator.

    \item Let $\Phi:L_\l\to L_\l \o U$ be an intertwiner.
Then the corresponding generalized character $\chi_\Phi$  satisfies the
following equation: 

    $$\left(\Delta_\h-2\sum_{\a\in R^+} \frac{e_\a
f_\a }{(e^{\a/2}-e^{-\a/2})^2}\right) (\delta \chi_\Phi)
=(\l+\rho, \l+\rho)\delta\chi_\Phi.\tag 2.2.4$$
\endroster\endproclaim

    \demo{Proof} Recall that the Casimir element was
defined by: $C=\sum_{\a\in R^+} (e_\a f_\a +f_\a e_\a)
+\sum x_l^2$, where $e_\a\in \g_\a, f_\a\in \g_{-\a}$ are
such that $(e_\a, f_\a)=1$, and $x_l$ is an orthonormal
basis in $\h$ with respect to $(\, , \,)$. Therefore, 
(1) follows from the following identities:

 $$\gather
\Tr_V(\Phi (\sum x_l^2) e^{2\pi\i h})= \Delta_\h \chi_\Phi, \tag 2.2.5\\
\Tr_V(\Phi e_\a f_\a e^{2\pi\i h})=-\frac{e_\a
f_\a}{(e^{\a/2}-e^{-\a/2})^2} \chi_\Phi - \frac{e^\a}{1-e^\a}
\d_\a \chi_\Phi,\tag 2.2.6\\
 \Tr_V(\Phi f_\a e_\a e^{2\pi\i h})=-\frac{f_\a
e_\a}{(e^{\a/2}-e^{-\a/2})^2} \chi_\Phi - \frac{1}{1-e^\a}
\d_\a \chi_\Phi.\tag 2.2.7\endgather$$

    The first of these identities is immediate corollary of
the definition of $\chi_\Phi$. Let us prove the second one.
Denote $X=\Tr_V (\Phi e_\a f_\a e^{2\pi\i h})$. Then, using
the identities $\Phi e_\a = (1\o e_\a+ e_\a\o 1) \Phi$,
$e^{2\pi\i h} e_\a = e^{2\pi\i \<\a, h\>} e_\a e^{2\pi\i
h}$ and the cyclic property
of the trace, we can write: 

    $$\aligned
X=&\Tr_V (\Phi e_\a f_\a e^{2\pi\i h})
= e_\a \Tr (\Phi f_\a e^{2\pi\i h})+ \Tr(\Phi f_\a e^{2\pi\i h} e_\a)\\
=& e_\a \Tr (\Phi f_\a e^{2\pi\i h})+ e^{2\pi\i \<\a, h\>}\Tr (\Phi f_\a
e_\a e^{2\pi\i h})\\
=&e_\a \Tr (\Phi f_\a e^{2\pi\i h})+ e^\a\Tr (\Phi (e_\a f_\a
-h_\a)e^{2\pi\i h})\\
=& e_\a \Tr (\Phi f_\a e^{2\pi\i h})+ e^\a (X-\d_\a \chi_\Phi),
\endaligned$$ 
(recall that $e^\a$ denotes both the element of $\C[P]$ and
a function on $\h$ given by $e^\a(h)=e^{2\pi\i \<\a,
h\>}$).   Therefore, 

    $$X=\frac{1}{1-e^\a}\left( e_\a \Tr (\Phi f_\a e^{2\pi\i h})
-e^\a\d_\a \chi_\Phi\right).$$

    Similar arguments show that 

    $$\Tr (\Phi f_\a e^{2\pi\i h})= \frac{f_\a}{1-e^{-\a}}
\chi_\Phi,$$ 
which gives us (2.2.6). Proof of (2.2.7) is quite similar. 

    This proves the first identity in (1); the second
identity can be checked by direct calculation or deduced from
Proposition~3.1.1 (see below).

    (2) obviously follows from (1) and Corollary~2.2.2
along with the fact that $C|_{L_\l}=(\l, \l+2\rho)$.\enddemo

    \demo{Example} Let $U=\C$ be the trivial
representation, $\Phi=\Id_{L_\l}$. Then Proposition~2.2.4
implies that the character $\ch L_\l$ satisfies the
following equation: 

    $$\Delta (\delta\ch L_\l) =(\l+\rho, \l+\rho)\delta\ch
L_\l;\tag 2.2.8$$
in particular, 

    $$\Delta \delta=(\rho, \rho)\delta.\tag 2.2.9$$

    These equations are well-known and have a nice
geometric interpretation -- as the heat equation on the
compact Lie group $G$ associated with $\g$. 
    \enddemo

\remark{Remarks} \roster \item
   Note that the theorems of this section can be easily
generalized to the case when $V$ is from category $\Cal O$.
In this case, we must consider the differential operators
with coefficients from the ring $\overline{\C[P]}$: since 
for $\a\in R^-$, $(1-e^\a)^{-1}=1+e^\a+e^{2\a}+\dots$ converges in
$\overline{\C[P]}$, there is no need to extend this ring. 

\item The equations (2.2.3), (2.2.4) are special cases of the formulas
for radial part of Laplace operator in the theory of $\mu$-special
functions, due to Harish-Chandra (see \cite{W, Chapter 9}). We discuss
the relation between generalized charcters and spherical functions in
Section~2.4 below. \endroster
\endremark

\head 2.3 Center of $\Ug$ and difference operators \endhead

    In this section we discuss the analogue of the
construction of Section~2.2 in the $q$-case. The
construction is quite parallel to the classical case, but
the differential operators should be replaced by difference
operators.  In this section, the words
``representation'' and ``intertwiner'' stand for
representation of $\Ug$ and $\Ug$-intertwiner, and the ground
field is the field $\C_q$ of rational functions in
$q^{1/2N}$ (see Section~1.2)

    Introduce
the following ring of difference operators, acting on the functions
$f\in \C_q[P]$:
$$DO_q=\{ D =\sum\limits_{\a\in\frac{1}{2} P\v} a_\a T_\a|\text{
almost all }a_\a=0\},$$
where $T_\a e^\l=q^{\<\a,\l\>}e^\l$, and $a_\a\in
\C_q[P](q^me^\mu-1)^{-1}$ are fractions of the form
$\frac{f}{g}, f,g \in \C_q[P]$ such that $g$ is a product
of factors  of the form $q^me^\mu-1$  for some $\mu\in Q, m\in
\Z$.

    Note that if we consider elements of $\C_q[P]$ as
functions on $\h$ then $T_\a$ can be rewritten as follows: 
$(T_\a f)(h)= f(h+\a \frac{\log q}{2\pi\i})$; however, to avoid
difficulties with proper definition of $\log q$ when $q$ is
a formal variable, we prefer not to use this language and
by ``function'' we always mean an element of $\C_q[P]$ (or
a its  field of fractions).  
However, to simplify the notations  we
sometimes  will write the generalized character
$\chi_\Phi$ as $\Tr (\Phi e^{2\pi\i h})$ meaning by this the element
of $\C_q[P]$ defined by (2.1.1). One can check that in fact
all our arguments can be carried in this formal language.  

    As in the previous section, let $V, U$ be
finite-dimensional representations of $\Ug$. 
    
\proclaim{Theorem 2.3.1} 
Let $I$ be the left ideal in $\Ug$ generated by $q^h-1$.
Then for any $u\in \Ug$ there exists a unique difference
operator $D_u\in DO_q\o  \Ug/I$
such that for any intertwiner $\Phi:V\to V\o U$ 
we have 

$$\Tr_V(\Phi u e^{2\pi\i h})=D_u\Tr_V(\Phi e^{2\pi\i h}).\tag 2.3.1$$

Again, 
 $D_u f$ is well defined for any function
$f$ with values in $U[0]$ .

\endproclaim

\demo{Proof} Without loss of generality we can assume that
$u$ is a monomial in the generators
$e_i, f_i, q^{h}$ of the form
$u=u^- u^0 u^+$, $u^\pm\in U^\pm, u^0\in U^0$. Define $\text{sdeg
}u=\text{deg }u^+ -\text{deg }u^-$, where $\text{deg }e_i=-\text{deg
}f_i=1$. We prove the theorem by induction in $\text{sdeg }u$.

If $\text{sdeg }u=0$ then $u=u^0=q^{\a}$ for some
$\a\in\frac{1}{2} P\v$. Then it follows immediately from the definition that
$D_u=T_\a$, so the theorem holds.

Let us make the induction step. Assume that $\text{sdeg }u>0$; then either
$\text{deg }u^+\ne 0$ or $\text{deg }u^-\ne 0$. We can  assume
that $u^-\in U[\mu],\mu\in -Q^+, \mu\ne 0$.

Since $\Phi$ is an intertwiner, $\Tr(\Phi
u^- u^0 u^+ e^{2\pi\i h})=\Tr (\Delta(u^-)\Phi u^0 u^+ e^{2\pi\i h})$. 
From the definition of comultiplication one easily sees that
$\Delta (u^-)= u^-\o  q^{\mu/2}+ \sum u_j\o  v_j$
for some   $u_j, v_j\in U^0 U^-$ such that 
$\text{sdeg }(u_ju_0u^+)<\text{sdeg }u$. Thus,

$$\Tr(\Phi u e^{2\pi\i h})= q^{\mu/2} \Tr (\Phi u^0 u^+ e^{2\pi\i h} u^- )+ \sum
v_j\Tr(\Phi u^0 u^+ e^{2\pi\i h} u_j).$$

Since commuting with $e^{2\pi\i h}$ does not change $\text{sdeg }u_j$, by
the induction assumption we can write

$$\align
\Tr(\Phi u e^{2\pi\i h})=& q^{\mu/2} \Tr (\Phi u^0 u^+ e^{2\pi\i h} u^- )+
 D'\Tr(\Phi e^{2\pi\i h})\\
=&q^{(\mu,\mu)'/2}e^\mu \Tr (\Phi u^0 u^+ u^- e^{2\pi\i h} )+
 D'\Tr(\Phi e^{2\pi\i h})\\
=&q^{(\mu,\mu)'/2}e^\mu \Tr (\Phi (u+[u^0 u^+, u^-])  e^{2\pi\i h} )+
D'\Tr(\Phi e^{2\pi\i h})\endalign$$
for some $D'\in DO_q \o \Ug$. Since $\text{sdeg }$ of all
terms in $ [u^0u^+, u^-]$ is less than
$\text{sdeg } u^- u^0 u^+$, we can again apply induction assumption and get

$$\Tr(\Phi u e^{2\pi\i h})=\frac{1}{1- q^{(\mu,\mu)'/2}e^\mu} D''\Tr(\Phi
e^{2\pi\i h}). $$

This proves the existence part of the theorem. Uniqueness follows from the
following lemma:
\proclaim{Lemma 2.3.2} Let us fix a finite-dimensional $\Ug$-module $U$.
Suppose that $D\in DO\o  \operatorname{ Hom  }(U[0], U[\mu])$ is such
that 
for any $\Ug$-intertwiner
$\Phi:V\to V\o  U, V$ --
finite-dimensional $\Ug$-module
we have $D\chi_\Phi=0$. Then $D=0$.
\endproclaim

\demo{Proof of the lemma} Let us assume that $D\ne 0$. Multiplying
$D$ by a suitable element from $\C_q[P]$ we can assume that $D$ has polynomial
coefficients: $D=\sum_{\l\in P} e^\l D_{(\l)}$, $D_{(\l)}$ being
difference operators with constant matrix-valued coefficients.  Let us
take the maximal (with respect to the standard  order in $P$) $\l$
such that $D_{(\l)}\ne 0$. Then if we have a generalized
character  $\chi$ such that $\chi=e^\mu u+\lot$ then,
taking the highest term of $D\chi$,  we see that
$D_{(\l)}(e^\mu u)=0$. On the other hand, if we take $\mu$
such that $\<\mu, \a_i\v\>\gg 0$ then for every $u\in U[0]$ 
there exists a non-zero
intertwiner $\Phi\colon L_\mu\to L_\mu\o  U$ such that the
corresponding generalized character
 has the form $\chi_\Phi=e^\mu u+\text{lower
order terms}$. Thus $D_{(\l)}(e^\mu u)=0$ for all $\mu\gg 0,
 u\in U[0]$. It is easy to show that it is only possible if
$D_{(\l)}=0$,  which contradicts the
assumption $D_{(\l)}\ne 0$.\qed
\enddemo \enddemo

\proclaim{Proposition 2.3.3} 
Let us keep the notations of Theorem~\rom{2.3.1}. Then 
    $c\mapsto D_c$ is an algebra homomorphism
of $\Cal Z(\Ug)$  to  $DO_q\o  \Ug[0]/I$.
\endproclaim

    \demo{Proof} The same as in the classical case (Theorem~2.2.3)\enddemo
    
    Unlike the classical case, in the quantum case we have
an explicit construction of central elements (see Theorem~1.2.1).
Applying the previous construction to those central
elements, we get the following theorem:

    \proclaim{Theorem 2.3.4} Let $V$ be a
finite-dimensional representation of $\Ug$,  $c_V$ --
the corresponding central element in $\Ug$ \rom{(}see
Theorem~\rom{1.2.1)}, and $D_{c_V}$ -- the corresponding
difference operator constructed in Theorem~\rom{2.3.1}. Let 
$\Phi:L_\l\to L_\l \o U$ be a $\Ug$-intertwiner, and
$\chi_\Phi$-- the corresponding generalized character. Then
$\chi_\Phi$ satisfies the following difference equation: 

    $$D_{c_V} \chi_\Phi= \ch V (q^{-2(\l+\rho)})\chi_\Phi.
\tag 2.3.2$$

    \endproclaim
\demo{Proof} The same as in classical case (see
Corollary~2.2.2); the value of $c_V$ in $L_\l$ is taken
from Theorem~1.2.1.\enddemo

    \remark {Remark} As in the classical case, all the
constructions of this section can be easily generalized to
the case where $V$ is an 
arbitrary module from category $\Cal O$. \endremark

\head 2.4 Generalized characters as spherical functions \endhead

In this section we show that generalized characters  for $\g$ can
be interpreted as spherical functions on the space $G\times G/G$. The
constructions of this section are  due to the author, I.~Frenkel and
P.~Etingof (see \cite{EFK}). 

As before, let $\g$ be a simple Lie algebra over $\C$, and let $G$ be
the compact real simply connected Lie group corresponding to $\g$. It is
known that every finite-dimensional complex  representation of $\g$
can be lifted 
to $G$. Also, it is known that every finite-dimensional representation
$V$ of $G$ is unitary: there exists a positive definite hermitian form
$(\, ,\,)_V$ on $V$ such that $(gv_1, gv_2)_V=(v_1, v_2)_V$. 

Let us fix some finite-dimensional representation $U$ of $\g$. Let
$C^\infty(G, U)$ be the space of all smooth functions on $G$ with values
in $U$. We can define a (hermitian) inner product in  $C^\infty(G, U)$
by 

$$(f_1, f_2) =\int_G  (f_1(g), f_2(g))_U\, dg, \tag 2.4.1$$ 
where $(\, , \,)_U$ is the inner product in $U$ and $dg$ is the Haar
measure on $G$. It is easy to see that the inner product defined by
(2.4.1) is positive-definite. We denote by $L^2(G, U)$ the closure of
the space $C^\infty(G, U)$ with respect to the norm $\Vert
f\Vert=\sqrt{(f,f)}$. 

\definition{Definition 2.4.1} A function $f\in C^\infty(G, U)$ is called
{\it equivariant} (notation: $f\in C^\infty(G, U)^G$) 
if for every $g,x\in G$ we have

$$f(gxg^{-1})=g f(x).\tag 2.4.2$$

The same applies to $f\in L^2(G,U)$. 
\enddefinition

This definition can be rewritten as follows. Let us consider the group
$G\times G$, and let $G_d\subset G\times G$ be the diagonal subgroup.
Consider the functions $f:G\times G\to U$ such that for every $k\in G_d,
x\in G\times G$ we have

$$\gathered 
f(xk)=f(x),\\
f(kx)=k f(x).\endgathered\tag 2.4.3$$

This is a special case  of what is called $\mu$-spherical functions
(see \cite{W}) on the group $G\times G$ with respect to the subgroup
$G_d$. On the other hand, it is easy to see that 

$$f(x,y)\mapsto f(xy^{-1})$$ 
establishes isomorphism between the space of spherical functions on
$G\times G/G_d$ defined by (2.4.3) and the space of equivariant
functions in the sense of Definition~2.4.1.

Here is yet another description of the same space. Let $f\in C^\infty(G,
U)^G$, and let $u\in U^*$. Define a complex-valued function on $G$ by 
$f_u(g)=\<u, f(g)\>$. Define the action of $G$ on scalar functions on $G$
by $(T_g f)(x) =f(g^{-1}xg)$. Then it is easy to see that $T_g
f_u=f_{gu}$, and thus $f_u$ 
satisfies the following condition: 

$$\aligned 
&\text{Under the action of $G$ defined above, $f_u$ spans a
finite-dimensional subspace } \\
&\text{in $C^\infty(G)$, and as a representation of
$G$, this space is isomorphic to $U^*$.}\endaligned\tag 2.4.4$$

Vice versa, it is easy to see that every scalar function on $G$
satisfying the condition (2.4.4) can be obtained as $f_u(g)$ for some
$f\in C^\infty(G, U)^G, u\in U^*$.

Now, let $V$ be a finite-dimensional representation of $G$ and let
$\Phi:V\to V\o U$ be a $G$-intertwiner. Define coresponding
generalized character $\chi_\Phi\in C^\infty(G, U)$
by 

$$\chi_\Phi(g)=\Tr_V(\Phi g).\tag 2.4.5$$

Note that for $g=e^{2\pi\i h}$ this coincides with the previously
given  definition of the
generalized character (see (2.1.2)). 

\proclaim{Lemma 2.4.2}\roster
\item For every intertwiner $\Phi:V\to V\o U$, the generalized character
$\chi_\Phi$ defined by \rom{(2.4.5)} is equivariant: $\chi_\Phi\in
C^\infty(G, U)^G$. 

\item Assume that $V$ is irreducible. Then $\chi_\Phi=0$ iff $\Phi=0$. 
\endroster\endproclaim

\demo{Proof} (1) is quite trivial: 

$$\Tr_V(\Phi gx g^{-1})=\Tr ((g\o g)\Phi x g^{-1} )=g\Tr (\Phi x).$$

(2) Let us consider $\chi_\Phi$ on the maximal torus, i.e. on the points
of the form $e^{2\pi\i h}, h\in \h_\R$. As was explained in Section~2.1,
we can as well consider it as an element of $\C[P]\o U$, and
$\chi_\Phi=0$ as a function on $\h$ iff $\chi_\Phi=0$ as an element of $\C[P]\o
U$. Let $v_\l$  be a highest-weight vector in $V$. Then $\Phi
v_\l=v_\l\o u_0+\lot $ for some $u_0\in U[0]$. It is known that $u_0=0$
iff $\Phi=0$. On the other hand, $\chi_\Phi=e^\l\o u_0 +\lot$, which
proves the theorem. \qed\enddemo

Let $\l\in P^+$. Recall that $L_\l$ is the irreducible highest-weight
module over $\g$ with highest weight $\l$, which in this case is
finite-dimensional and thus can be considered as a module over $G$. Let 
$H_\l=\Hom_G (L_\l, L_\l\o U)$ (note that this space is
finite-dimensional). Due to Lemma~2.4.2, we have an injective map

$$H_\l \to C^\infty(G, U)^G: \Phi\mapsto \chi_\Phi.$$

We will denote the image of $H_\l$ in $C^\infty(G, U)^G$ also by $H_\l$.

\proclaim{Theorem 2.4.3} \roster
\item
For $\l\ne \mu$, $H_\l$ and $H_\mu$ are orthogonal with respect to the
inner product in $C^\infty(G, U)$ \rom{(}see \rom{(2.4.1))}. 

\item 
$$L^2(G,U)^G=\bigoplus_{\l\in P^+} H_\l.$$
\rom{(}direct sum should be understood in the sense of Hilbert spaces\rom{)}. 
\endroster\endproclaim

We refer the reader to \cite{EFK} for the proof of this theorem. 

\example{Example} Let $U=\C$ be the trvial representation. Then
Theorem~2.4.3 coincides with well-known Peter-Weyl theorem. \endexample

Finally, since the conjugacy classes in $G$ are the same as $W$-orbits
in the maximal torus 

$$T=\exp(\i \h_\R)\simeq \h_\R/Q\v, \quad \h_\R=\bigoplus \R
\a_i\v,\tag 2.4.6$$
it is easy to see that every equivariant function on $G$ is uniquely
defined by its values on $T$. More precisely, 
we have the following isomorphism: 
$$C^\infty(G, U)^G \simeq C^\infty (T, U[0])^W, \tag 2.4.7$$
i.e. functions on $T$ with values in $U[0]$ satisfying 
$f(w t)=wf(t)$ for every element from the Weyl group. 

Thus, restriction of the inner product defined by (2.4.1) to
the space of equivariant functions can be rewritten in terms of
integral over $T$. The answer is given by the following theorem:

\proclaim{Theorem 2.4.4} If $f_1, f_2\in C^\infty(G, U)^G$ then

$$(f_1, f_2)=\frac{1}{|W|} \int_T (f_1(t), f_2(t))_U \delta\bar\delta\, dt,
\tag 2.4.8$$
where, as before, $dt$ is the Haar measure on $T$, $(\, , \,)_U$ is
the inner product in $U$, and $\delta$ is the Weyl denominator: 

$$\delta(e^{2\pi\i h})=\prod_{\a\in R^+} \bigl( e^{\pi\i \<\a, h\>} -
e^{-\pi\i \<\a, h\>}\bigr),$$
and $\bar \delta$ is the complex conjugate of $\delta$. 
\endproclaim
\remark{Remark} Again, in the case $U=\C$ this theorem is well known;
it was first proved by H.~Weyl.\endremark

\demo{Proof} It is very easy to check that the scalar function $\psi(g)=
(f_1(g), f_2(g))_U$ is conjugation invariant. Thus, the theorem is
reduced to the following statement: for every central function $\psi$
on $G$, we have

$$\int_G \psi (g) \, dg  =\frac{1}{|W|} \int_T \psi(t)\delta\bar\delta\, dt,$$
which is well-known (essentially, this is the result of Weyl we
referred to above).\qed 
 \enddemo

Note that this theorem along with the orthogonality statement 
of Theorem~2.4.3 gives a new proof of the hermitian version of the 
orthogonality theorem for
generalized characters for $\g$ (Theorem 2.1.2); of course, these
arguments wouldn't help in the $q$-case.

Also, using the description of the space of equivariant functions as
scalar-valued functions on $G$ transforming under the conjugation as
$U^*$ (see (2.4.4)), it is easy to see that every
conjugation-invariant scalar differential operator $D$ on $G$ defines
an operator in $C^\infty(G, U)^G$, acting by a scalar in every $H_\l$
(see details in \cite{EFK}). Since the space of such operators is
isomorphic to  $(U\g)^\g\simeq \Cal Z(U\g)$, this gives  a
natural action of $\Cal Z(U\g)$ by differential operators on
$C^\infty(G, U)^G$. Using the isomorphism (2.4.7), we can rewrite this
action in terms of differential operators on $T$ with coefficients
from $\text{End }U[0]$ (this operation is often referred to as
``taking the radial part of the Laplacian''). It can be easily shown
that the differential operators appearing in this way coincide with
those constructed in Theorem~2.2.1. In particular, the operator $D_C$
given by (2.2.3) is exactly the radial part of the second order
Laplace-Beltrami operator $\Delta_G$. We refer the reader to
\cite{EFK} for more detailed information. Note that it would be quite
difficult to calculate this radial part by geometrical arguments,
whereas the algebraic approach makes it very simple. 

Again, we note that the radial part of the Laplace operator has been
calculated in more general situation by Harish-Chandra (see \cite{W,
Chapter 9}), so (2.2.3) can be deduced from his results. 
   

\vfill\newpage

\vbox{\vskip 1in}
\specialhead\chapter{3} 
\centerline{JACOBI POLYNOMIALS} 
\endspecialhead

\vbox{\vskip 0.5in}

    In this chapter we show that in some special  cases
the generalized characters for the Lie algebra $\sln$
coincide with so-called Jack polynomials,  
     which were
studied by Heckman and Opdam (see \cite{HO, O1, H1, H2}). These
results are due to the author and Pavel Etingof (see \cite{EK2, EK4}). 
In fact, Heckman and Opdam defined and studied analogues of these
polynomials for arbitrary root systems; they called them
Jacobi polynomials associated with a root system, 
the reason being  that for $\g=\frak{sl}_2$ these polynomials coincide with
so-called ultraspherical (Gegenbauer) polynomials, which are a
special  case of Jacobi polynomials. In this chapter we
only consider classical case, i.e. representations,
intertwiners and generalized characters for $\g$ rather
than for $\Ug$.

\head 3.1 Jacobi polynomials and corresponding differential operators 
\endhead

    In thie section we give definition and main properties
of Jacobi polynomials associated with an arbitrary reduced
irreducible root system $R$, following papers of Heckman
and Opdam.
    
    Recall (see Chapter 2) that $\C[P]$ denotes 
 the group algebra of the
weight lattice, and $\C[P](e^\a-1)^{-1}$ is the ring 
obtained by adjoining to $\C[P]$ 
the expressions of the form $(e^\alpha -1)^{-1}, \alpha\in
R$. Note that the elements of $\C[P]$ may be considered as
functions on $\h$ by the rule:
$e^\a(h)=e^{2\pi\i \<\a,h\>}$. Under this identification 
 elements of $\C[P](e^\a-1)^{-1}$ become functions 
 on the real torus $T= \h_\R/Q\v, \h_\R=\bigoplus \R\alpha\v_i$ 
with singularities
on the hypersurfaces $e^\a (h)=1, \a\in R$. However, we will use the
formal language as far as possible. Abusing the language,
we will call elements of $\C[P]$ polynomials, and elements
of $\C[P]^W$ symmetric polynomials. Similarly, we will talk
of divisibility of polynomials meaning divisibility in the
ring $\C[P]$. 

    Recall also (see Section~2.2) that we denote by $DO$
the ring of differential operators in $\h$ with
coefficients from $\C[P](e^\a-1)^{-1}$; again, they can be
treated formally as derivations of the ring
$\C[P](e^\a-1)^{-1}$. In particular, 
for every $\a\in \h^*$ we defined $\d_\a$ so that  $\d_\a e^\l =
(\l, \a)e^\l$ and  the Laplace
operator $\Delta_\h$ so that 
$\Delta_\h e^\l = (\l,\l)e^\l$.

    Through this whole chapter, we fix some positive
integer $k$. 

    Consider the following differential operator

$$L= L_k=\Delta_\h - k(k-1)\sum\limits_{\alpha\in
R^+}\frac{(\a,\a)}{(e^{\alpha/2}- e^{-\alpha/2})^2}.\tag 3.1.1$$

 This operator for the root system $A_n$  was introduced 
by Sutherland (\cite{Su}) and for an arbitrary root system by
Olshanetsky  and Perelomov
(\cite{OP}) as a Hamiltonian of an integrable quantum system. We will
call $L$  the  Sutherland operator.

    As before, let $\delta$ be the Weyl denominator defined
by (1.1.6).  

Define the following version of the  Sutherland operator:

$$M_k=\delta^{-k}(L_k-k^2(\rho,\rho))\delta^k.\tag 3.1.2$$

\proclaim{Proposition  3.1.1}{\rm (\cite{HO})} 
\roster
\item 
$$M_k=\Delta_\h-k\sum_{\alpha\in R^+}\frac{1+e^\alpha}
{1-e^\alpha} \d_\alpha
=\Delta_\h -2k \sum_{\a\in R^+} \frac{1}{1-e^\a}\d_\a +2k\d_\rho.
 \tag 3.1.3$$

\item Both $L_k, M_k$ commute with the action of the Weyl group.

\item $M_k$ preserves the algebra of symmetric polynomials
$\C[P]^W\subset \C[P](e^\a-1)^{-1}$.
\endroster\endproclaim

    \demo{Proof} We do not give the proof of (1) here,
referring the reader to \cite{HO}. 
Note that the proof involves some non-trivial statement
about the root systems, which we will discuss later in the
proof of the  affine analogue of this statement (Theorem~7.1.2).

    (2) is obvious, and (3) immediately follows from (1).
\enddemo 
    
    Recall (see 1.1.1) that we denoted by $m_\l$ the basis of orbitsums
in $\C[P]^W$: 
$m_\l=\sum\limits_{\mu\in W\l}e^\mu,\l\in P^+$.  

\proclaim{Lemma 3.1.2}
$$M_km_\l =(\l, \l+2k\rho)m_\l+\sum\Sb \mu<\l\\ \mu\in P^+\endSb
c_{\l\mu}m_\mu.\tag 3.1.4$$
\endproclaim
\demo{Proof} Explicit calculation.\enddemo

Now we can consider the eigenfunction  problem for $M_k$. Let
us  consider the action of $M_k$ in the finite-dimensional space
spanned by $m_\mu$ with $\mu\le \l$. Then the eigenvalue $(\l,
\l+2k\rho)$ has multiplicity one in this space due to the following
trivial but very useful fact:

\proclaim{Lemma 3.1.3} Let $\l,\mu\in P^+$, $\mu<\l$. Then $(\mu+\rho,
\mu+\rho)<(\l+\rho, \l+\rho)$. \endproclaim

Thus, we can give the following definition: 

\proclaim{Definition 3.1.4} Jacobi polynomials $J_\l, \l\in P^+$ are
the elements of $\C[P]^W$ defined by the following conditions:

\roster 
\item $J_\l=m_\l+\sum\limits_{\mu<\l}c_{\l\mu}m_\mu$.

\item $M_kJ_\l=(\l, \l+2k\rho)J_\l$.
\endroster
\endproclaim

Due to Lemma 3.1.3, these properties determine $J_\l$ uniquely. Note
that this definition is valid for any complex $k$, and
components of  $J_\l$ are rational functions of  $k$.

Let us introduce an inner product in $\C[P]^W$. Let

$$\< f,g\>_0=\frac{1}{|W|}[f\bar g]_0,\tag 3.1.5$$ 
where, as in Section~2.1, 
 $[\,\,]_0$ is the constant term of a polynomial, and the bar
involution is defined by $\overline{e^\l}=e^{-\l}$. More generally,
let

$$\<f,g\>_k=\<f\delta^k, g\delta^k\>_0.\tag 3.1.6$$

    Note that for $k=1$ this definition coincides with
previously given (2.1.3).

\proclaim{Lemma 3.1.5} 
$M_k$ is self-adjoint with respect to the inner product
$\<\cdot,\cdot\>_k$. \endproclaim

\demo{Proof} This is equivalent to $L_k$ being self-adjoint
with respect to the inner product $\<\cdot,\cdot\>_0$, which is
obvious.\enddemo 

\proclaim{Corollary 3.1.6} 
$\<J_\l, J_\mu\>_k=0 $ if $\l<\mu$.\endproclaim

    In fact, one has a stronger result:

\proclaim{Theorem  3.1.7}{\rm (\cite{O1})} 
$\<J_\l, J_\mu\>_k=0 $ if $\l\ne \mu$.\endproclaim

    We will prove this theorem for the root system $A_n$ in
the next section  --
see Theorem~3.2.4.

    In fact, the operator  $M_k$ can be included in  a  large
commutative family of differential operators. Define 

    $$\Bbb D=\{D\in DO| D\text{ is $W$-invariant}, [D,
M_k]=0\}.\tag 3.1.7$$

    Then we have the following theorem: 
\proclaim{Theorem 3.1.8} 
\roster\item Every operator $D\in \Bbb D$ preserves the
space $\C[P]^W$, and $J_\l$ is a common eigenbasis for the
action of $\Bbb D$ in $\C[P]^W$: there exists a map
$\gamma: \Bbb D\to (S\h)^W$ such that 

    $$D J_\l=\gamma(D)(\l+k\rho) J_\l.\tag 3.1.8$$

    \item $\gamma$ is an isomorphism $\Bbb D\simeq
(S\h)^W$. 
\endroster\endproclaim

    We do not give the proof of this theorem here,
referring the reader to the above mentioned papers of
Heckman and Opdam (in fact, for the root system $A_n$ this
theorem was known before). However, we note here that part
(1) is relatively easy, and so is injectivity of $\gamma$;
the difficult part is to prove surjectivity, or
to construct $\gamma^{-1}$. Again, we will reprove it by
representation-theoretic methods for the root system $A_n$
in the next section.

    \head 3.2 Jack polynomials as generalized characters\endhead

    In this section we show that 
one can get Jacobi polynomials for the root system
$A_{n-1}$, in which case they are also known under the name ``Jack
polynomials'', as generalized characters. 
This construction is due to the author and Pavel Etingof
(\cite{EK2, EK3}).     {\bf In this section,
we only consider the case $\g=\sln$}.

    Recall that we have fixed a positive integer $k$.
Define the representation $U=U_{k-1}$ (which we will later use to
define the generalized characters) as the irreducible
representation of $\sln$ with the highest weight $(k-1)n
\omega_1$; it can be described as the symmetric power of
the fundamental representation:  

    $$U=S^{(k-1)n} \C^n.\tag 3.2.1$$

    In other words, $U$ can be identified with the space of
homogeneous polynomials in $x_1, \dots, x_n$ of degree
$(k-1)n$. The action of $\sln$ is given by the following
formulas: 

    $$\gathered 
    e_i\mapsto x_i\d_{i+1},\quad f_i\mapsto x_{i+1}\d_i\\
h_i\mapsto x_i\d_i-x_{i+1}\d_{i+1},\endgathered\tag 3.2.2$$
where $\d_i=\frac{\d}{\d x_i}$. 

    It is very important for us that all weight subspaces of
$U$ are one-dimensional. In particular, the same is true
for $U[0]$; we fix an element $u_0=(x_1\dots x_n)^{k-1}\in
U[0]$, which allows us to identify 

    $$U[0]\simeq \C:\quad u_0\mapsto 1. \tag 3.2.3$$

Also, we will use the fact that 

    $$e_\a f_\a|_{U[0]}=k(k-1)\Id_{U[0]}.$$
    
    \proclaim{Lemma 3.2.1} Let $\l\in P^+$. A non-zero $\sln$-homomorphism 
$\Phi\colon  L_\l\to
L_\l\otimes U$ exists iff $\l-(k-1)\rho\in P^+$; if it
exists,  it is unique
up to a factor.\endproclaim

    \demo{Proof} Standard exercise.\enddemo

For brevity, from now on we use the following notation: 

    $$\lk=\l+(k-1)\rho.\tag 3.2.4$$

For every $\l\in P^+$ fix an intertwiner 

    $$\Phi_\l: L_\lk\to L_\lk \o U\tag 3.2.5$$
by the condition $\Phi_\l(v_\lk)=v_\lk\o u_0+\lot$, where
$v_\lk$ is the highest weight vector in $L_\lk$. 

    Define by $\varphi_\l$ the corresponding generalized
character: 

    $$\varphi_\l= \chi_{\Phi_\l}= \Tr_{L_\lk} (\Phi_\l
e^{2\pi\i h}).\tag 3.2.6$$

    It takes values in the space $U[0]$; since this space
is one-dimensional, we can identify it with $\C$ using
(3.2.3) and consider $\varphi_\l$ as scalar-valued function. 
Under this identification, $\varphi_\l$ has the following
form: $\varphi_\l=e^{\l+(k-1)\rho}+\lot$. 

\proclaim{Proposition 3.2.2}
$$\varphi_0=\delta^{k-1}.\tag 3.2.7$$
\endproclaim

    The proof of this proposition will be given later in a
more general case (see Proposition~4.2.2). 

    Now we can prove the main theorem of this section: 

    \proclaim{Theorem 3.2.3} 
Let $\varphi_\l, \l\in P^+$ be the generalized characters
for $\sln$  defined by formulas
\rom{(3.2.6), (3.2.7)}. 
    Then
$\varphi_\l$ is divisible by $\varphi_0$, and 
the ratio $\varphi_\l/\varphi_0$  is the Jack  polynomial
$J_\l$ \rom{(}see Definition~\rom{3.1.4)}.\endproclaim

\demo{Proof} Let us first prove that 
$\varphi_\l$ is divisible by
$\varphi_0$, and the ratio is a  symmetric
polynomial with highest term $e^\l$.
Consider the tensor product $V=L_\l \otimes L_{(k-1)\rho}$.  It 
decomposes as follows: $V=L_{\l+(k-1)\rho}+ \sum_{\mu< \l}
N_{\l\mu}L_{\mu+(k-1)\rho}$. Consider the intertwiner
$\Phi=\Id_{L_\l}\otimes \Phi_0 :V\to V\otimes U$. On one
hand, it follows from the definition that $\chi_\Phi
= \ch {L_\l} \cdot \varphi_0$, where $\ch {L_\l}$ is the (usual)
character of $L_\l$. On the other hand, the decomposition of
$V$ implies that $\chi_\Phi = \varphi_\l +\sum_{\mu<\l}
a_{\l\mu} \varphi_\mu$, and thus 
$\varphi_\l/\varphi_0=\ch {L_\l} +\sum_{\mu<\l}
a_{\l\mu} \varphi_\mu/\varphi_0$. Since $\ch {L_\l}$ is a
symmetric polynomial, it follows by
induction in $\l$ that $\varphi_\l/\varphi_0$ is also  a
symmetric polynomial.

Let us prove that $\varphi_\l/\varphi_0$ satisfies the equation
$M_k (\varphi_\l/\varphi_0)=(\l,
\l+2k\rho)(\varphi_\l/\varphi_0)$. Since
$\varphi_0=\delta^{k-1}$ (see Proposition~3.2.2), this
equation is equivalent to 
$$L_k \delta\varphi_\l= (\l+k\rho, \l+k\rho)\delta\varphi_\l.\tag 3.2.8$$

    On the other hand, recall that we have proved in
Section~2.2 that the generalized characters are
eigenfunctions of a certain differential operator $D_C$
which was obtained from the Casimir element $C$.  Comparing
the expressions (2.2.3)  for $D_C$ and (3.1.1) for $L_k$, we see
that 

    $$L_k=\delta\circ\pi_U (D_C)\circ \delta^{-1},\tag 3.2.9$$
since $e_\a f_\a u_0 =k(k-1)u_0$. Therefore, (3.2.8)
immediately follows from the differential equation for
generalized characters derived in Proposition~2.2.4.  
\qed\enddemo

    This immediately implies the orthogonality of Jack
polynomials: 

    \proclaim{Theorem 3.2.4} If $J_\l$ are Jack polynomials
then 

    $$\<J_\l, J_\mu\>_k=0\quad \text{if } \l\ne \mu.\tag 3.2.9$$ 
\endproclaim

    \demo{Proof} It follows immediately from Theorem~3.2.3
and the orthogonality theorem for generalized characters
(Theorem~2.1.2). 
\enddemo

    Thus we have reproved by representation-theoretic
arguments the orthogonality theorem~3.1.7 for Jacobi
polynomials for the root system $A_{n-1}$. Note that 
the proof given by Heckman and Opdam uses transcendental
(i.e., not algebraic) methods; the
representation-theoretic approach outlined above gives
probably the simplest proof of this fact. On the other hand,
assuming (3.2.9) we could give an alternative proof of
Theorem~3.2.3 -- using the orthogonality theorem for
generalized characters rather than the differential
equation satisfied by them. We will use this approach in
the next chapter, where we discuss $q$-analogue of Jacobi
polynomials -- Macdonald's polynomials.

Our construction also allows one to construct differential operators
commuting with $M_k$. Namely, it follows from Theorem~2.2.3 that the
operators of the form $D=\delta^{-(k-1)} \pi_U (D_c)\delta
^{k-1}$, where $c\in \Cal
Z(U\g)$, commute with each other (and thus, with $M_k$, which can be
obtained from the Casimir element) and are $W$-invariant. 
This means that the map 

    $$\Cal Z(U\g)\simeq (S\h)^W\to \Bbb
D: c\mapsto \delta^{-(k-1)} \pi_U (D_c)\delta
^{k-1}$$
is the inverse map to the map $\gamma$ defined by (3.1.8).
This proves (for the root system $A_{n-1}$) surjectivity of
$\gamma$, i.e. the difficult part of Theorem~3.1.8 on the
structure of $\Bbb D$.


\vfill\newpage

\vbox{\vskip 1in}
\specialhead\chapter{4} 
\centerline{MACDONALD'S POLYNOMIALS}
\centerline{ AND GENERALIZED CHARACTERS} 
\endspecialhead

\vbox{\vskip 0.5in}

    In this chapter we develop $q$-analogue of the
constructions of Chapter~3. We define the $q$-analogue of
Jacobi polynomials, called Macdonald polynomials, and show that in some
special cases the
generalized characters for the quantum group $U_q\sln$
coincide with Macdonald's polynomials for the root system
$A_{n-1}$. 
This construction is due to the author and Pavel Etingof
\cite{EK2}. 

    In this chapter, the words ``representation''
etc. always stand for representation of the quantum group
$\Ug$. 

    \head 4.1 Definition of Macdonald's polynomials\endhead 

In this section we give the definition of Macdonald's polynomials for
reduced irreducible root system, following \cite{M2}. 

We preserve the notations of Chapter~1. 
Consider the algebra of $W$-invariant
polynomials $\C[P]^W$. The main goal of this section is to construct a
basis in $\C_{q,t}[P]^W$, where $\C_{q,t}=\C(q,t)$ is the field of
rational functions in two independent variables $q,t$. Unless
otherwise stated, in this section ``polynomial'' will stand for an
element of $\C_{q,t}[P]$, and divisibility will stand for
divisibility in this ring; similarly, elements of
$\C_{q,t}[P]$ will be called symmetric polynomials. 

    Recall (see (1.1.1)) that we denoted by $m_\l$ the basis of orbitsums
in $\C[P]^W$: 
$m_\l=\sum\limits_{\mu\in W\l}e^\mu, \l\in P^+$.  
    
\proclaim{Theorem 4.1.1} {\rm (Macdonald)} 

There exists a unique family of
polynomials $P_\l(q,t)\in \C_{q,t}[P]^W, \l\in P^+$ such that 

\roster 
\item $P_\l=m_\l +\sum_{\mu<\l} c_{\l\mu}m_\mu$.

\item These polynomials are orthogonal with respect to the following
inner product on $\C_{q,t}[P]$:

$$\<f, g\>_{q,t}= \frac{1}{|W|} [f\bar g\Delta_{q,t}]_0,\tag 4.1.1$$
where, as before, the bar involution is defined by $\overline{e^\l}=e^{-\l}$, 
$[\ \ ]_0$ is the constant term: $[\sum c_\l e^\l]_0=c_0$, and 

$$\Delta_{q,t}=\prod\limits_{\a\in R} \prod\limits_{m=0}^\infty
	\frac{1-q^{2m}e^{\a}}{1-q^{2m}t^2e^{\a}}.\tag 4.1.2$$
\endroster
\endproclaim

These polynomials are called Macdonald's polynomials (our notation
slightly differs
 from that of Macdonald: what we denote by $P_\l(q,t)$ in
the notations of \cite{M2} would be $P_\l(q^2,t^2)$).

    \remark{Remark} In fact, for non simply-laced systems
there  this theorem can be generalized, allowing
different variables $t_\a$ for roots of different  lengths;
see \cite{M2} for details.\endremark

It is often convenient to consider Macdonald's polynomials for
$t=q^k, k\in \Z_+$ (see examples below). 
In this case, Macdonald's polynomials lie in $\C(q)[P]$; 
we will write $\<\ ,\ \>_k$ instead of $\<\ ,\
\>_{q,q^k}$, etc. (note that it agrees with Definition~3.1.6
for $q=1$). 
However, most of the properties of Macdonald's polynomials  obtained for
$t=q^k$ can be generalized to the case when $q,t$ are independent
variables. Abusing the notations, we will sometimes say ``$t=q^k$, $k$
is an independent variable''  instead of saying that $q,t$ are
independent variables. 

    \example{Example} \roster\item For $k=0$, we have
$P_\l=m_\l$ independently of $q$; for $k=1$, $P_\l=\ch L_\l$
-- also independently of $q$.

    \item In the limit $q, t \to 1$ so that $t=q^k$
Macdonald's polynomials tend to Jacobi polynomials
introduced in Section~3.1, which follows from Theorem~3.1.7.
\endroster\endexample

For the case $\g=\sln$ (that is, when $R$ is the root system of type
$A_{n-1}$), one can slightly modify the above definition
and define Macdonald's polynomials for the $\frak{gl}_n$
as a basis in $\C[x_1, \dots, x_n]^{S_n}$ labeled by the
partitions $\l$.
We will keep the
same notation $P_\l(x; q,t)$ for these polynomials. 
In this form they were introduced in \cite{M1}. 

For the root system $A_{n-1}$ the polynomials $P_\l(x; q,t)$ can be
defined in a different way; namely, they can be defined as the
eigenfunctions of a certain family of commuting difference
operators. Recall that we have identified  $\C[P]^W$ for the root
system $A_{n-1}$ the space of  
with symmetric polynomials in $x_1^{\pm 1}, \dots,
x_n^{\pm 1}$ of degree zero (see end of Section~1.1). Let
us define the following operators acting in $\C[P]^W$: 

$$M_r=t^{r(r-n)} \sum\Sb I\subset \{1,\ldots, n\}\\|I|=r\endSb
  \biggl(\prod\Sb i\in I \\j\notin I\endSb
       \frac{t^2x_{i} -x_j}{x_{i}-x_j} \biggr)
	T_{q^2, I},\tag 4.1.3$$ 
where 
$(T_{q^2,i}f)(x_1,\ldots,x_n)=f(x_1,\ldots, q^2x_i,\ldots,
x_n)$, $T_{q^2, I}=\prod_{i\in I} T_{q^2, i}$  
and $r=1,\ldots,n$. It is not too
difficult (though it is not quite obvious) to show that
these operators preserve the space of symmetric
polynomials. 
 Note
that since the total degree is zero, $M_n f=f$ for every
$f\in \C[P]^W$, and any difference operator acting in polynomials in
$x_1, \dots, x_n$  is defined by
its action in $\C[P]^W$ uniquely  up to a multiple of $M_n$. We could
consider the case of $\frak{gl}_n$ rather then $\sln$ thus
replacing $\C[P]^W$ by $\C[x_1, \dots, x_n]^{S_n}$; then
all $M_i$ act non-trivially. Let us, however, stick to the
$\sln$ case. 

    \proclaim{Theorem 4.1.2} {\rm (Macdonald)}
\roster\item $ [M_i, M_j]=0$.

\item  $M_r$ is self-adjoint with respect to the inner product
$\<\cdot,\cdot\>_{q,t}$.

    \item $M_rP_\l(x;q,t)=c_\l^r P_\l(x;q,t)$, and 

    $$c_\l^r =\sum_{|I|=r} \prod_{i\in I} 
q^{2\l_i}t^{2\rho_i} ,$$
where $\rho_i=\frac{n+1-2i}{2}$
    \rom{(}see end of Section~\rom{1.1)}. 
\endroster
\endproclaim

    This characterization of Macdonald's polynomials 
is analogous to the definition of Jacobi polynomials as
eigenfunctions of the  commuting family  $\Bbb D$ 
of differential
operators given in Chapter~3. In fact, one can show that any differential
operator $D\in \Bbb D$ can be obtained as a 
certain linear combinations of coefficients of expansion of
Macdonald's difference operators $M_i$ in powers of
$(q-1)$. However, these expressions are rather messy; we
will return to this relation later (see remark at the end
of the next section). 

    We also note (though it is not relevant for our
purposes) that the same holds for an arbitrary root system:
Macdonald's polynomials can be defined as eigenfunctions of
a certain family of commuting difference operators (see
\cite{C2}); unfortunately, for root systems other than
$A_n$ it is very difficult to write these difference
operators explicitly.

\head 4.2 Macdonald's polynomials of type $A$ as generalized characters
\endhead

Through this section, we assume $t=q^k, k\in \N$ and show how one gets
Macdonald's polynomials $P_\l(x;q,q^k)$ for the root system
$A_{n-1}$ as generalized characters. The construction is
quite parallel to that of Section~3.2. In this section we only
consider $\g=\sln$; unless otherwise specified, the words
``representation'', ``intertwiner'', ``generalized
character'' stand for representation of $U_q\sln$, etc.
Also, all the objects are defined over the field
$\C_q=\C(q^{1/2n})$.

Let $U$ be the finite-dimensional representation
of $\Ug$ with the highest
weight $(k-1)n\omega_1$; this is a $q$-analogue of the representation
$U=S^{(k-1)n}\C^n$ considered in Section~3.2.  
As before, this 
representation can be realized explicitly in homogeneous 
polynomials of degree
$(k-1)n$ of $n$ variables $x_1,\ldots, x_n$ with the action
of $U_q\sln$ given by 

$$\gathered
h_i\mapsto x_i\frac{\d}{\d x_i}- x_{i+1}\frac{\d}{\d x_{i+1}},\\ 
e_i\mapsto x_i D_{i+1},\quad
f_i\mapsto x_{i+1} D_i,\\
(D_i f)(x_1,\ldots, x_n)=\frac{f(x_1,\ldots, qx_i,\ldots ,x_n)-
f(x_1,\ldots, q^{-1}x_i,\ldots ,x_n)}{(q-q^{-1})x_i}.
\endgathered\tag 4.2.1$$

    Since the multiplicities in tensor products for $\Ug$
are the same as for $\g$, we have the following proposition:
    
    \proclaim{Proposition 4.2.1}Let $\l\in P^+$.  A non-zero
$U_q\sln$-homomorphism $\Phi\colon L_\l\to 
L_\l\otimes U$ exists iff $\l-(k-1)\rho\in P^+$; if it
exists,  it is unique
up to a factor.\endproclaim

As in Chapter~3, we denote $\lk=\l+(k-1)\rho$ and 
define the intertwiner $\Phi_\l$ and the corresponding generalized
character $\varphi_\l$: 

$$\gathered
\Phi_\l:L_\lk\to L_\lk\o U,\\
\varphi_\l=\chi_{\Phi_\l}.\endgathered\tag 4.2.2$$

As before, we consider $\varphi_\l$ as
a scalar-valued function, and choose the identification
$U[0]\simeq \C$ so that $\varphi_\l=e^{\l+(k-1)\rho}
+\lot$.

\proclaim{Proposition 4.2.2}
$$\varphi_0=\prod_{i=1}^{k-1}\prod_{\alpha\in R^+}
(e^{\alpha/2}-q^{2i}e^{-\alpha/2}).\tag 4.2.3$$\endproclaim

\demo{Proof} First, we prove the following statement:

\proclaim{Lemma 1} $\varphi_\l$ is divisible by
$(1-q^{2j}e^{-\alpha})$ for any  positive root $\alpha$ and $1\le
j\le k-1$.\endproclaim

The proof is done in several steps. 
Let us introduce $F_i=f_iq^{-d_i h_i/2}$; then 

$$\Delta(F_i)= F_i\o q^{-d_i h_i}+1\o F_i.$$

Let $F$ be a (non-commutative) polynomial in $F_1,\ldots, F_{n-1}$ of
weight $-\alpha, \alpha\in Q^+$.
Let  $\varphi_\l ^F=\Tr_{L_\l} (\Phi_\l Fe^{2\pi\i h})$. 
Also, let us fix a basis
in $U$: $U[\alpha]=\C u_\alpha, \alpha\in Q$. Then

\proclaim{Lemma 2}  There exists a polynomial $P_F\in \C(q)[P]$ such
that 

$$\varphi_\l ^F=
\frac{P_F\varphi_\l}{\prod\Sb \beta\le\alpha\\ \beta\in
Q^+\endSb (1-q^{(\beta,\beta)}e^{-\beta})}u_{-\alpha}.\tag 4.2.4$$
\endproclaim

Proof is by induction in $\alpha\in Q^+$. For $\alpha=0$ the statement
is obvious. Now, let $\alpha=\sum m_i\alpha_i, \sum m_i=m$ and assume that the
statement is proved for all $\alpha'<\alpha$. Take $F=F_{j_1}\ldots
F_{j_m}$. Then $\Delta F=\Delta(F_{j_1})\ldots \Delta(F_{j_m})=
\sum_i q^{\sigma_i} F(\gamma_i)\otimes \tilde
F(\alpha-\gamma_i)q^{-h_{\gamma_i}} +F\otimes q^{-h_\alpha}$, 
where $\gamma_i\in Q^+, \gamma_i<\alpha, \sigma_i\in\Z$, and
$F(\gamma)$ has weight $-\gamma$. Therefore, using the intertwining
property of $\Phi_\l$ and the cyclic property of the trace, we
get 

$$\varphi_\l^F=\Tr(\Delta(F)\Phi_\l e^{2\pi\i h})=
q^{-h_\alpha}\Tr (F\Phi_\l e^{2\pi\i h})+ A=q^{(\alpha,
\alpha)}e^{-\alpha}\varphi_\l^F  +A,$$
where $A=\sum q^{\sigma_i}\Tilde
F(\alpha-\gamma_i)q^{-h_{\gamma_i}}\Tr(F(\gamma_i)\Phi_\l
e^{2\pi\i h})$. Thus, 

$$\varphi_\l^F=\frac{1}{1-q^{(\alpha,\alpha)}e^{-\alpha}}A.$$

On the other hand, it follows from the induction assumption that $A$ is an
expression of the form (4.2.4) containing only the factors
$1-q^{(\beta, \beta)}e^{-\beta}$ with $\beta<\alpha$ in the denominator.
This completes the proof of Lemma~2. 

\proclaim{Lemma 3} Let $\alpha\in R^+$, and let $F_\alpha$ be a
\rom{(}non-commutative\rom{)} polynomial in $F_i$ which in
the limit $q=1$ becomes 
a root element of $\frak s\frak l _n$. Then 
$P_{F_\alpha^{k-1}}$ is a non-zero polynomial  relatively prime to
$\prod\limits_{j=1}^{k-1}(1-q^{2j}e^{-\alpha})$. \endproclaim

It suffices to prove this lemma for $q=1$. But for $q=1$,
$\Delta(F_\alpha) =F_\alpha\otimes 1+1\otimes F_\alpha$, and therefore

$$\split
\varphi_\l^{F_\alpha^{k-1}}=\Tr(\Phi_\l F_\alpha^{k-1}e^{2\pi\i h})=
F_\alpha \Tr (\Phi_\l F_\alpha^{k-2}e^{2\pi\i h}) + \Tr(\Phi
F_\alpha^{k-2}e^{2\pi\i h} F_\alpha)\\
=F_\alpha\varphi_\l^{F_\alpha^{k-2}}+
e^{-\alpha}\varphi_\l^{F_\alpha^{k-1}},\endsplit$$
so 

    $$\varphi_\l^{F_\alpha^{k-1}}=(1-e^{-\alpha})^{-1}
		F_\alpha\varphi_\l^{F_\alpha^{k-2}}
=\ldots=(1-e^{-\alpha})^{1-k}\varphi_\l
F_\alpha^{k-1}u_0.$$  

Since $F_\alpha^{k-1}u_0=c_\alpha u_{(1-k)\alpha}$ for some
$c_\alpha\ne 0$, we see that 

$$P_{F_\alpha^{k-1}}=c_\alpha \frac {\prod\limits_{\beta\le
(k-1)\alpha} (1-e^{-\beta}) }{(1-e^{-\alpha})^{k-1}}=
c_\alpha \prod\Sb \beta<(k-1)\alpha\\\beta\ne s\alpha\endSb
(1-e^{-\beta}) 
\prod\limits_{s=1}^{k-1}(1+e^{-\alpha}+\ldots+e^{(1-s)\alpha}).$$

One can easily see that this polynomial is relatively prime to
$1-e^{-\alpha}$. Thus, we have proved Lemma~3.

Now, let us return to the proof of Lemma~1. Let us write 

$$\Tr (\Phi_\l
F_\alpha^{k-1}e^{2\pi\i h})=\frac{P_{F_\alpha^{k-1}}\varphi_\l} 
	{\prod\limits_{j=1}^{k-1}(1-q^{2j}e^{-\alpha})}.$$
Since the left-hand side is a non-zero polynomial,
and  $P_{F_\alpha^{k-1}}$
is relatively prime to
$\prod\limits_{j=1}^{k-1}(1-q^{2j}e^{-\alpha})$, we see that
$\varphi_\l$  must be divisible by
$\prod\limits_{j=1}^{k-1}(1-q^{2j}e^{-\alpha})$. So, Lemma~1 is
proved. 

Now it is easy to prove Proposition 4.2.2: Lemma~1 implies that 
we have the following identity: 
$\varphi_0=f\prod\limits_{j=1}^{k-1}\prod\limits_{\alpha\in R^+}
(1-q^{2j}e^{-\alpha})$ for some polynomial $f$; comparing
the highest and the lowest terms on both sides we see that 
$f=e^{(k-1)\rho}$.
This completes the proof of Proposition 4.2.2.
\qed
\enddemo

Now we can prove the main theorem of this section: 
\proclaim{Theorem 4.2.3} If $\l\in P^+$ 
and $\varphi_\l, \varphi_0$ are the generalized characters 
for $U_q\sln$ 
defined by \rom{(4.2.2), (4.2.3)}     then
$\varphi_\l$ is divisible by $\varphi_0$, and 
the ratio $\varphi_\l/\varphi_0$  is the Macdonald's polynomial
$P_\l(q,q^k)$ for the root system $A_{n-1}$.\endproclaim

    \demo{Proof} First, the same arguments as in the proof
of Theorem~3.2.3 -- no changes are needed -- show that
$\varphi_\l/\varphi_0$ is a symmetric polynomial with
highest term $e^\l$. Therefore, 
to prove the theorem it suffices to prove that these
ratios are orthogonal with respect to the inner product
$\<\cdot,\cdot\>_k$. This immediately follows from the 
orthogonality theorem for generalized characters 
and formula for $\varphi_0$
(Proposition~4.2.2).   Indeed, we know from
the orthogonality theorem that $[\varphi_\l \bar
\varphi_\mu\Delta]_0=0$ if $\l\ne \mu$. Therefore, 
$[(\varphi_\l/\varphi_0)\overline{(\varphi_\mu/\varphi_0)}
\varphi_0\bar\varphi_0\Delta]_0=0$. Due to Proposition~4.2.2,
$\varphi_0\bar\varphi_0\Delta=\prod\limits_{\alpha\in\
R}\prod\limits_{i=0}^{k-1} (1-q^{2i}e^{\alpha})=\Delta_{q,t}$,
which proves the  orthogonality of $\{\varphi_\l/\varphi_0\}$ 
with respect to the inner product $\<\cdot, \cdot\>_{k}$.
\qed\enddemo

\remark{Remark 4.2.4} Note that the proof only uses the following
properties of $U$: 
\roster \item All weight subspaces are one dimensional.

\item $e_\a^k u_0=0, e_\a^{k-1}u_0\ne 0$ for every $\a\in R^+$.
\endroster

(we also used Proposition~4.2.1, which ensures existence of
intertwiners; however, it can be deduced from properties (1), (2)
above.)

Thus, the same theorem must be true if we replace $U$ by any other
representation satisfying these properties. However, one can check
that for $U_q\sln$ we have only two possibilities:
$U=L_{(k-1)n\omega_1}$ (which we used) or $U=L_{(k-1)n\omega_{n-1}}$,
which can be obtained from the previous one by an outer automorphism
of $\g$, i.e. by flip of the Dynkin diagram. For other Lie algebras,
such representations do not exist at all (except some small number of
exceptional cases, where such representations exist only for finite
number of values of $k$), which explains why this theory can not be
generalized to arbitrary root systems in a trivial way. \endremark

\remark{Remark} This theorem can be generalized for the case of
generic $k$, i.e., the case where $q,t$ are independent variables; see
details in \cite{EK2}.  
\endremark 

\head 4.3  The center of $U_q\sln$ and Macdonald's operators\endhead
In this section we show how one can get Macdonald's operators $M_r$
introduced in Section~4.1 from the quantum group $U_q\sln$. 
This construction is
parallel to the one for $q=1$ (see Section~3.2).
As before, in this section we only consider $\g=\sln$ and
 $t=q^k, k\in \N$.

    Recall (see Section~2.3) that we have denoted by $DO_q$
the ring of difference operators, i.e. operators of the
form $\sum\limits_{\a\in \half P\v} a_\a T_\a$, 
where $T_\a e^\l=q^{\<\a,\l\>}e^\l$, and $a_\a\in
\C_q[P](q^me^\mu-1)^{-1}$. We have also constructed for
every element $c$ of the center of $U_q\sln$ a difference
operator $D_c\in DO_q$ such that for every intertwiner
$\Phi$ we have 
$\Tr (\Phi c e^{2\pi\i h})=D_c \Tr (\Phi e^{2\pi\i h})$.
These operators commute, and generalized characters are
their common eigenfunctions. 

    Let us show that applying this construction to the
above case (i.e., $\g=\sln$ and $U$ is chosen as in the
beginning of Section~4.2) we can get Macdonald's difference
operators $M_r$ defined by (4.1.3).

\proclaim{Theorem 4.3.1}  Define the central elements 
$c_r\in \Cal Z(U_q\sln), r=1\,\ldots, n-1$ by
$$c_r= c_{\Lambda^{n-r}_q},\tag 4.3.1$$
\rom{(}cf. Theorem~\rom{1.2.1)}, 
where $\Lambda^i_q$ is the $q$-deformation of the representation of $\sln$
in the $i$-th exterior power $\Lambda^i\C^n$ of the
fundamental representation, and let $D_{c_r}\in DO_q$ be
the corresponding difference operator \rom{(}see Theorem~\rom{2.3.1)}. Then

$$M_r=\varphi_0^{-1}\circ D_{c_r}\circ \varphi_0,\tag 4.3.2$$
where $M_r$ is Macdonald's operator introduced in
\rom{(4.1.3)}, and 
$\varphi_0$ is the
operator of multiplication by the function $\varphi_0$  defined by
\rom{(4.2.3)}. \endproclaim

\demo{Proof} Let us first prove that 

    $$\varphi_0^{-1} D_{c_r} (\varphi_0 P_\l) =
\sum_{I}q^{2\sum_{i\in I} (\lambda+k\rho)_i}P_\l,\tag 4.3.3$$
where the sum is taken over all subsets $I\subset\{1,
\dots, n\}$ of cardinality $r$. Indeed, we know from
Theorem~2.3.4 that 

    $$D_r \varphi_\l = \ch \Lambda^{n-r}
(q^{-2(\l+k\rho)})\varphi_\l.$$ 

    Since $\ch \Lambda^{n-r} =e^{\frac{n-r}{r}(1, \dots,
1)} \sum e^\mu$, where the sum is taken over all
$\mu=(\mu_1, \dots, \mu_n)$ such that $\mu_i =0$
or $1$, $\sum \mu_i=n-r$, we get 

    $$\ch
\Lambda^{n-r}(q^{-2(\l+k\rho)})=\sum_{I}q^{2\sum_{i\in I}
(\lambda+k\rho)_i} .$$ 

Since $P_\l=\varphi_\l/\varphi_0$, we get (4.3.3).

    Comparing (4.3.3) with the formula (4.1.4) for eigenvalues of
Macdonald's operators, we see that $M_r$ and 
$\varphi_0^{-1}\circ D_{c_r}\circ \varphi_0$ coincide on
Macdonald's polynomials, and thus, on all symmetric
polynomials.  Repeating the uniqueness arguments
outlined in the proof of Lemma~2.3.2, we see that it
is  only possible if they are equal.
\qed\enddemo

    \remark{Remarks} 
\roster \item Recently a straightforward proof of
Theorem~4.3.1 was found by Mimachi (\cite{Mi}).

    \item
The central elements $c_r$  are closely related to those
constructed in \cite{FRT}. Essentially, the central elements
constructed in \cite{FRT} are traces of the powers of $L$-matrix,
whereas $c_r$ are coefficients of the characteristic
polynomial of $L$.\endroster\endremark

This theorem  allows one to see the relation between the
differential operators from $\Bbb D$ (see (3.1.7)) and Macdonald's difference
operators. This relation is nothing but the relation between the
center of $U_q\sln$ (which, as we have seen, is spanned by the elements
$c_V= (1\otimes \Tr_V)(\Cal R^{21}\Cal R(1\otimes q^{-2\rho}))$) and
the center of $U\sln$, which does not have a nice explicit
description, but 
can be described by means of
Harish-Chandra isomorphism $\Cal Z(U\g)\simeq (S\h)^W$.


\vfill\newpage

\vbox{\vskip 1in}
\specialhead\chapter{5} 
\centerline{INNER PRODUCT AND SYMMETRY IDENTITIES}
\centerline{ FOR  MACDONALD'S POLYNOMIALS } 
\endspecialhead

\vbox{\vskip 0.5in}

In this chapter we use the technique of generalized characters
developed in the previous sections to prove so-called inner product
(norm) and symmetry  identities for
Macdonald's polynomials of type $A_n$. These results are due to the
author and Pavel Etingof (\cite{EK5}). 

 Inner product identities
express the norm $\<P_\l, P_\l\>$ as a certain product over the
positive roots. They have been conjectured by Macdonald for arbitrary
root systems. He also gave a proof for $A_n$ (unpublished); the proof
for  arbitrary root systems was given in a recent paper of
Cherednik \cite{C2}. 

Symmetry identity relates the values of
$P_\l(q^{2(\mu+k\rho)})$ and $P_\mu (q^{2(\l+k\rho)})$. For $A_n$ it
was first proved by Koornwinder (unpublished), so the first published
proof is in \cite{EK5}. Again, recently Cherednik proved this identity
for arbitrary root systems (\cite{C3}). 

In this chapter we only consider root system of type $A_{n-1}$,
i.e. $\g=\sln$. Our proofs use the realization of Macdonald's
polynomials as generalized characters for $U_q\sln$ 
(see Chapter~4) and the  technique of
representing identities in the category of representations of a
quantum group by ribbon graphs, developed by Reshetikhin and Turaev. 
We refer the reader to their papers \cite{RT1, RT2} or to recent books
of Turaev \cite{T} and Kassel \cite{Kas}  for description of
this technique; a very brief introduction can be found in the Appendix
to \cite{EK5}.

\head 5.1 Inner product identities \endhead

Let us fix a positive integer $k$. 
Recall (see Chapter~4) that we have defined Macdonald's polynomials
$P_\l\in \C_q[P], \l\in P^+$ and an inner product $\<\, ,\,\>_k$ in
$\C_q[P]$ such that $\<P_\l, P_\mu\>_k=0$ if $\l\ne\mu$. The goal of
this section is to calculate $\<P_\l, P_\l\>_k$; our proof is based on
Theorem~4.2.3, which shows that $P_\l$ can be expressed in terms of
generalized characters for $U_q\sln$. 

Recall the notations
$U=U_{k-1}=S^{(k-1)n}\C^n, u_0=u_0^{k-1}\in U[0], 
 \lk=\l+(k-1)\rho, \Phi_\l: L_\lk\to L_\lk \o U,
\varphi_\l=\chi_{\Phi_\l}= e^\lk \o u_0 +\lot $
introduced in Chapter~4; also, recall the Chevalley involution 
 $\omega$  and Shapovalov form $(\, ,\,)_V: V\o V^\omega\to \C$, 
discussed in Section~1.2. We assume that  Shapovalov form in $U$ is
normalized so that $(u_0, u_0)_U=1$. 

\proclaim{Lemma 5.1.1} The  inner product $\<P_\l,
P_\l\>_k$ can be  calculated from 

$$\Psi \1_\l= \<P_\l, P_\l\>_k\1_\l,\tag 5.1.1$$
where  $\1_\l$ is an invariant
vector in $L_\lk\otimes L_\lk^\omega$ and $\Psi:L_\lk\otimes
L_\lk^\omega\to L_\lk\otimes L_\lk^\omega$   is the following
operator:

$$L_\lk\o L_\lk^\omega @>{\Phi_\l\o \Phi_\l^\omega}>>
 L_\lk\o U\o
U^\omega\o L_\lk^\omega @>{\text{Id}\o
(\cdot,\cdot)_U\o \text{Id}}>>
L_\lk\o L_\lk^\omega.\tag 5.1.2$$
\endproclaim

\demo{Proof} It follows from Theorem~4.2.3 that $\<P_\l,
P_\l\>_k=\<\varphi_\l, \varphi_\l\>_1$, where the inner product on
generalized characters is introduced in Theorem~2.1.2. Now we can
repeat the same arguments we used in the proof of
Theorem~2.1.2.\qed\enddemo

It will be convenient to rewrite this in a slightly different way as
follows: 

\proclaim{Theorem 5.1.2} In the notations of previous Lemma, we have:

$$\<P_\l, P_\l\>_k=(\<v_\lk^*,
\Phi_\l\Phi_\l^{\circ}v_\lk\>)_{U},\tag 5.2.3$$
where the intertwiner 
$\Phi_\l^{\circ}:L_\lk\to L_\lk \otimes U^\omega$ is defined by
the condition $\Phi_\l^{\circ}(v_\lk)=v_\lk\otimes
(u_0)^\omega+\lot$.\endproclaim

\demo{Proof} 

The proof is obvious if we use the technique of ribbon graphs. 
Namely, it follows from Lemma~5.1.1 that the inner product 
$\<P_\l, P_\l\>_k=A_\l$ can be defined from the following identity
 of ribbon graphs: 

\psfig{figure=ident1.ps}
where $\text{dim}_qL=\Tr_L(q^{-2\rho})$, 
$\phi:L_\lk^*\to L_{\lk^*}, \psi:L_{n(k-1)\omega_1^*}\to
L_{n(k-1)\omega_1}^*=U_{k-1}^*$ are isomorphisms and $\psi$ is chosen so that 
$\<u_0, \psi(u_0^\omega)\>=1$. It is easy to check
that 

\psfig{figure=ident2.ps}
and thus, 

\psfig{figure=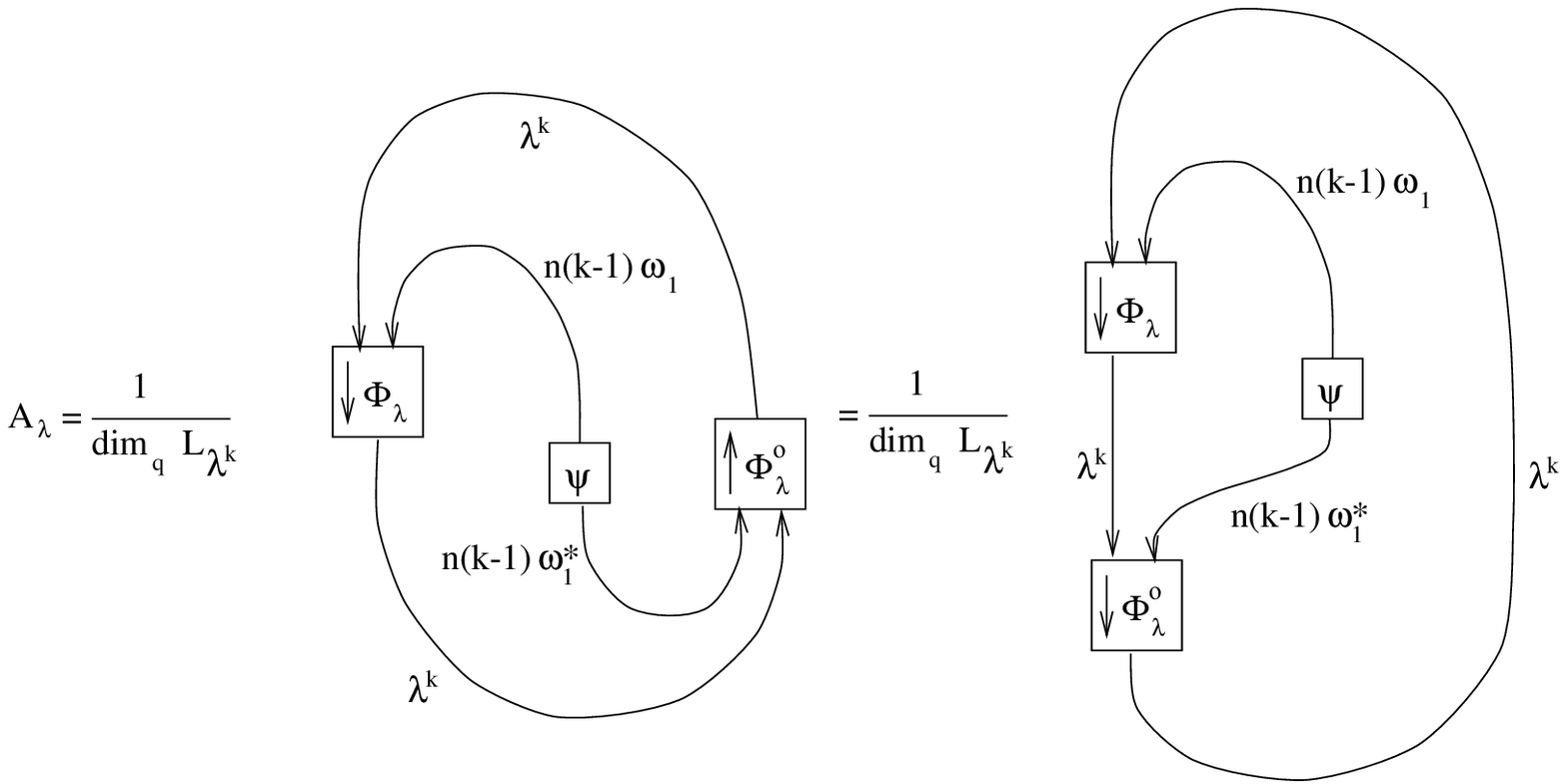}

so

\psfig{figure=ident3a.ps}

\qed\enddemo

\proclaim{Theorem 5.1.3} In the notations of the previous theorem, we
have: 

$$
(\Phi_\l\Phi_\l^{\circ})_U= \prod_{\a\in R^+}\prod_{i=1}^{k-1} 
\frac{1-q^{2(\a, \l+k\rho)+2i}}
     {1-q^{2(\a, \l+k\rho)-2i}}\Id_{L_{\lk}}.\tag 5.1.4
$$
\endproclaim

This theorem will be proved in the next section, using  detailed
analysis of the  poles of the intertwining operators. 

This theorem, along with Theorem~5.1.2 immediately gives the main
result of this section: 

\proclaim{Theorem 5.1.4}{\rm (Macdonald)}
$$\<P_\l, P_\l\>_k=\prod_{\a\in R^+}\prod_{i=1}^{k-1} 
\frac{1-q^{2(\a, \l+k\rho)+2i}}
     {1-q^{2(\a, \l+k\rho)-2i}}.\tag 5.1.5
$$\endproclaim

This is precisely the Macdonald's inner product identity for the  root
system $A_{n-1}$. 

\remark{Remark} As before, these identities can be easily generalized
to the case of arbitrary $k$ -- see \cite{EK5}.\endremark

\head 5.2 Algebra of intertwiners. \endhead 

This section is devoted to the proof of Theorem~5.1.3, which expresses
product of two intertwining operators in terms of a single
intertwiner. The technique we use here is considering the intertwining
operators for arbitrary value of $\l$ (which, of course, requires use
of Verma modules) and study their behavior (in particular, poles) as
functions of $\l$, which can be done  using the Shapovalov determinant
formula. 

Let us start with  more general situation. Consider intertwiners
of the form:

$$\Phi_\l^\mu:M_\l\to M_\l \otimes L_\mu,\tag 5.2.1$$
where
 $M_\l$ is  Verma module over $\Ug$, $L_\mu$ is the finite-dimensional
irreducible module; we assume that $\mu\in P^+\cap Q$, so $L_\mu[0]\ne 0$.
 Let $u\in L_\mu[0]$. We consider all the modules over the
field of rational functions $\C_q=\C(q^{1/2N})$, where $N=|P/Q|$;
 if $\l$ is not an integral
weight then we also have to add $q^{\<\l, \a_i\v\>/2N}$ to this field. 

It is known that if
$M_\l$ is irreducible then there exists a unique intertwiner of the
form (5.2.1) such that $\Phi^\mu_\l(v_\l)=v_\l\otimes u+\lot$.
 We will denote this intertwiner by $\Phi^{\mu, u}_\l$. 
The same is true if we consider the weight $\l$ as indeterminate, i.e.
if we consider $t_i=q^{\<\l, \a_i\v\>/2N}$ as algebraically independent
variables over $\C_q$.

Let us identify $M_\l$ with $U^-$ in a standard way. Then
we can say that we have a family of actions of $\Ug$ in the same space
$M\simeq U^-$,
and thus we have a family of intertwiners $\Phi^{\mu, u}_\l: M\to
M\otimes L_\mu$, defined for generic values of $\l$. This identification
also gives rise to Shapovalov form in $U^-$. We fix a homogeneous basis
$g_k$ in $U^-$ and denote the matrix elements of Shapovalov form in this
basis by $F_{kl}$: $F_{kl}=(g_k v_\l, g_l v_\l)_{M_{\l}}= (v_\l,
(\omega S (g_k))g_l v_\l)_{M_\l}$. 

For $\lambda\in
\h^*$, let us call a trigonometric rational function of
$\l$ a  rational function in $q^{1/2N}, q^{\l/2N}$ (that is, in $q^{1/2N}$
 and $t_i=q^{\<\a_i\v, \l\>/2N},i=1,\ldots, n$) and call  a
trigonometric polynomial in $\l$  a polynomial in $q^{\pm \l/2N}$
with coefficients from $\C_q$. Note that the ring of
trigonometric polynomials is a unique factorization ring, and
invertible elements in this ring are of the form $c(q)q^{\<\l, \a\>}$,
$\a\in\frac{1}{2N}Q^{\vee}$.

\proclaim{Lemma 5.2.1} For fixed $\mu\in P^+, u\in L_\mu[0], u\ne 0$ 
we have

\roster\item
 Let $\l\in \h^*$ be such that $M_\l$ is irreducible. Then
 there exists a unique 
intertwining operator $\Phi^{\mu, u}_\l:M_\l\to M_\l\o L_\mu$ such
that $\Phi^{\mu, u}_\l(v_\l)=v_\l\o u+\lot$.
Its matrix elements are trigonometric rational functions of
$\l$. Moreover, we have the following formula for $\Phi v_\l$: 

$$\Phi^{\mu, u}_\l (v_\l)= \sum_{k,l} (F^{-1})_{kl}g_k v_\l\o
q^{\l+2\rho}\omega (g_l)u,\tag 5.2.2$$
where $g_k$ is a homogeneous basis in $U^-$, $F^{-1}$ is the 
inverse matrix to the Shapovalov form in $M_\l$, and as in
Section~\rom{1.2}, $q^\l|_{V[\nu]}=q^{(\l, \nu)'}\Id_{V[\nu]}$. 

\item Define the operator $\tiphi$ by

$$\tiphi=
 d_\mu(\l)\Phi^{\mu, u}_\l,\tag 5.2.3$$
where
$$ \gathered
d_\mu(\lambda)=\prod_{\a\in R^+}\prod_{i=1}^{n^\a_\mu} 
\bigl(1-q^{(2(\a, \l+\rho)'-i (\a,\a)')}\bigr),\\
n_\mu^\a=\max\{i\in \Z_+| L_\mu[i\a]\ne 0\}.\endgathered
\tag 5.2.4$$

Then matrix elements of $\tiphi$ are trigonometric polynomials, i.e. have
no poles; thus, $\tiphi$ is well defined for all $\l$. 

\item  Consider the special case $\g=\sln, \mu=kn\omega_1$.  Recall
that in this case $L_\mu[0]$ is one-dimensional. Let $u$ be a non-zero
vector in $L_\mu[0]$. Then 
$d_\mu(\l)$ is the least common denominator of matrix elements of
$\Phi^{\mu, u}_\l$; in other words, in this case matrix elements of
$\tiphi$ do not have non-trivial common divisors. 
\endroster
\endproclaim

\demo{Proof} The proof is essentially the same as in the classical
case, which is given in \cite{ES}; however, we repeat it here marking
necessary changes. To prove (1), it suffices to check that the vector in
the right-hand side is the unique highest-weight vector of weight $\l$
in $M_\l\o L_\mu$ of the form $v=v_\l\o u+\dots$.
 Suppose $v\in M_\l\o L_\mu$ is of the form above. Define $E_i= e_i
q^{d_ih_i/2}$; then $v$ is highest-weight iff $\Delta E_i v=0$. On the
other hand, explicit calculation shows that $\Delta E_i=E_i\o 1 -
q^{2d_i} (1\o S E_i) q^{d_ih_i}\o q^{d_ih_i}$. Thus, 

$$\Delta E_i v =
(E_i\o 1 - q^{d_i(2+\<\l, \a_i\v\>)}1\o SE_i)v= 
(E_i \o 1 -q^{(\a_i, \l+2\rho)'}1\o SE_i)v.$$  

For $x\in M_\l, w\in M_\l\o L_\mu$ define 
$(x, w)$ by $(x, w_1\o w_2)= (x, w_1)_{M_\l} w_2 \in L_\mu$. 
If $M_\l$ is irreducible, Shapovalov form in $M_\l$ is non-degenerate,
and therefore, $w=0\iff (x, w)=0$ for all $x\in M_\l$. Therefore, we can
rewrite the condition that $v$ is a highest-weight vector as follows: 

$$\aligned 
&\Delta E_i v=0\iff (x, \Delta E_i v)=0\quad\text{for all }x\in M_\l
\iff\\
&(x, (E_i \o 1 -q^{(\a_i, \l+2\rho)'}1\o SE_i)v)=0\iff\\
&(\omega S(E_i)x, v)= q^{(\a_i, \l+2\rho)'} SE_i (x, v).
\endaligned$$

It is easy to see that the last condition is equivalent to the
following: for any homogeneous $F\in U^-$, we have

$$(Fv_\l, v)=q^{-(\text{wt }F, \l+2\rho)'} \omega(F)(v_\l, v)=q^{-(\text{wt }F,
\l+2\rho)'} \omega(F) u.$$ 

This proves that the highest-weight vector of the desired form 
exists and is unique. It is easy to check that the vector given by
(5.2.2) satisfies the condition above. 

To prove (2), note that  it follows from (5.2.2) that 
matrix coefficients of $\Phi$ may 
have poles only  at the points where the determinant of Shapovalov form
vanishes. The formula for the determinant of the Shapovalov form in
the quantum case can
be found in \cite{CK}, and the factors occuring there
are precisely the factors in formula (5.2.4) (up to invertible
factors). One can check in the same way as it is done in \cite{ES} for
$q=1$ -- that is, by comparing the order of pole of the matrix of
Shapovalov form and its minors -- that in fact all the poles of
$F^{-1}$ are simple. 

The restriction on $i$ in (5.2.4) appears because the coefficients
$(F^{-1})_{kl}$ which may have poles of the form (5.2.4) with $i>n_\mu^\a$
appear with zero coefficient. 

To prove (3), it suffices to check  this statement for $q=1$, which is
also done in the paper \cite{ES} (note that it is quite non-trivial!). 
\qed 
\enddemo

\remark{Remark} It is seen from this proof that $\tilde\Phi^{\mu,u}_\l$ is 
actually a trigonometric polynomial in $\l$ with operator coefficients,
i.e. the degrees of its matrix coefficients, as trigonometric polynomials, 
are uniformly bounded (under a suitable definition of degree). 
\endremark

We will discuss  what happens to these intertwiners when $\l\in
P^+$ and how they are related with intertwiners $L_\l\to L_\l\o L_\mu$
later.

We will also need one more technical lemma. 
\proclaim{Lemma 5.2.2}  Let us write $\tiphi=d_\mu(\l)\Phi^{\mu,
u}_\l$ \rom{(}see \rom{(5.2.4))} in the
following form:

$$\tiphi v_\l= d_\mu(\l)v_\l\otimes u +\ldots + a(\l) v_\l \otimes u_\mu,$$
where $u_\mu$ is the highest-weight vector in $L_\mu$, and 
$a(\l)\in \C_q[q^{\pm \l/2N}]\otimes U^-[-\mu]$ is a trigonometric
polynomial of $\l$ with values in the 
 universal enveloping algebra. Then the greatest common divisor of the
components of $a(\l)$ is 1.
\endproclaim

\demo{Proof} 
It is easy to see, using the irreducibility of $L_\mu$, that if
$a(\l)=0$ then $\tiphi=0$. On the other hand, we have shown before
that the coefficients of $\tiphi$ have no nontrivial common divisors,
and thus $\tiphi$ could  only vanish on a subvariety of
codimension more than one. Thus, the same must be true for $a(\l)$.
\enddemo

Now we want to define a structure of algebra on these intertwiners.
Let $\Phi_1:M_\l\to M_\l\otimes L_{\mu_1}, \Phi_2:M_\l\to
M_\l\otimes L_{\mu_2}$ be non-zero intertwiners. Let us define their
product $\Phi_1 *\Phi_2:M_\l\to M_\l \otimes L_{\mu_1+\mu_2}$ as the
composition 

$$M_\l @>{\Phi_2}>> M_\l\otimes L_{\mu_2}
	@>{\Phi_1\otimes 1}>>M_\l\otimes L_{\mu_1}\otimes L_{\mu_2}
	@>{1\otimes \pi}>>M_\l\otimes L_{\mu_1+\mu_2}, \tag 5.2.5$$
where $\pi$ is a fixed  projection $\pi:L_{\mu_1}\otimes L_{\mu_2}\to
L_{\mu_1+\mu_2}$.

Now, let us apply these notions to the situation considered in the
previous sections in connection with Macdonald's theory. That is,
assume that $\g=\sln$, $\mu=kn\omega_1$ for some $k\in \Z_+$, so that 
$L_\mu=U_k= q-\text{analogue of }S^{kn}\C^n$. 
As before, we identify $U_k$ with the space of
homogeneous polynomials of degree $kn$ in variables $x_1, \dots, x_n$
(see (4.2.1)), and let $u_0=(x_1\ldots
x_n)^k\in L_\mu[0]$; then $L_\mu[0]=\C_q u_0^k$. For brevity, we will
write $U_k$  for $L_{kn\omega_1}$, 
$\Phi^k_\l$ for $\Phi^{\mu=kn\omega_1, u_0^k}_\l$, etc. Let us fix
the projection $\pi:U_{k}\otimes U_{l} \to
U_{k+l}$ by $\pi(u_0^k\otimes u_0^l)=u_0^{k+l}$.  
In this case, $n_\mu^\a=k$ for all $\a\in R^+$ and

$$ d_{k}(\l)=\prod_{\a\in R^+}\prod_{i=1}^k 
	(1-q^{2(\a, \l+\rho)-2i}).\tag 5.2.6$$

Here comes the main result of this section: 

\proclaim{Theorem 5.2.3} 

$$\tilde\Phi^k_\l*\tilde\Phi^l_\l
=\tilde\Phi^{k+l}_\l.\tag 5.2.7$$
\endproclaim

\demo{Proof}
Let us denote  the left-hand side of (5.2.7) by $\Psi$. Then $\Psi$ is
an intertwiner $M_\l\to M_\l\otimes U_{k+l}$, whose matrix
coefficients are trigonometric polynomials in $\l$. In particular, we
can write $\Psi(v_\l)=f(\l)v_\l\otimes u^{k+l}_0+\text{l.o.t}$. 
On the other hand 
$\tilde\Phi^{k+l}_\l(v_\l) =  d_{k+l}(\l)v_\l\otimes
u^{k+l}_0+\text{l.o.t}$. Since the intertwining operator is unique for
generic $\l$, this implies
$\Psi(\l)=\frac{f(\l)}{ d_{k+l}(\l)}\tilde\Phi^{k+l}_\l$.
Since the greatest common divisor of the matrix elements of
$\tilde\Phi^{k+l}$ is 1, this implies that $f(\l)$ is divisible by
$ d_{k+l}(\l)$.

Let us now consider the lowest term of $\Psi$. If we write the lowest
term of $\tilde\Phi^k_\l$ as $a_k(\l)v_\l\otimes u_k$ (cf. Lemma~5.2.2) and
lowest term of $\Phi^l_\l$ as $a_l(\l)v_\l\otimes u_l$ then the
lowest term of $\Psi$ will be $a_l(\l)a_k(\l)v_\l\otimes u_{k+l}$ (up
to some power of $q$). Since
we know that components of $a_k$ have no common divisors, and the same
is true for $a_l$, it follows that the greatest common divisor
of components of $a_k(\l)a_l(\l)$ is 1. Indeed, suppose that $p(\l)$
is a common divisor of components of $a_k(\l)a_l(\l)$. Passing if
necessary to a certain algebraic extension of $\C_q$ we get that
$a_k(\l)a_l(\l)$ vanishes on a certain subvariety of codimension 1. On
the other hand, this contradicts  the fact that both $a_k, a_l$ could
only vanish on subvarieties of codimension more than one, since $U^-$
has no zero divisors. Thus, the greatest common 
divisor of coefficients of $\Psi(\l)$ is one, which implies that
$ d_{k+l}(\l)$ is divisible by $f(\l)$.

This proves that $\Psi(\l)=c(q)q^{(\l, \a)}
\tilde\Phi^{k+l}_\l$ for some $\a\in \frac{1}{2N}Q$ and  
rational function $c(q)$, independent of $\l$. To calculate $\a, c$, let
us consider the limit of both sides of (5.2.7) as $\l\to +\rho\infty$,
 i.e. letting $t_i=q^{\<\l, \a_i\v\>/2N}=0$.

\proclaim{Lemma 5.2.4} 
$$\lim_{\l\to\rho\infty}\Phi_k(v_\l) = v_\l\otimes u^k_0.$$
\endproclaim
\demo\nofrills{}To prove the lemma, note first that 
 $\lim \tilde \Phi^k=\lim \Phi^k$.
 Due to Lemma~5.2.1, we can write 
$\Phi(v_\l)=\sum_{k,l} (F^{-1})_{kl}g_k
v_\l\otimes (\omega g_l)u$. It is known that if we choose a basis
$g_k$ in such a way that $g_0=1, g_k$ has strictly negative weight for
$k>0$, then  $\lim_{\l\to\rho\infty}(F^{-1})_{kl}=
\delta_{k,0}\delta_{l,0}$ (this follows, for example,  from \cite{L2,
Proposition 19.3.7}, which gives much more detailed information about the
asymptotic behavior of $F$; it states that under a suitable normalization
the Shapovalov form in the Verma module 
$M_\l$ formally converges to the Drinfeld's form on 
$U^-$ as $\l\to +\infty\rho$). This proves the lemma.  
\enddemo

Using this lemma and the fact that in the identification $M_\l\simeq
M\simeq U^-$ the action   of
$U^-$ does not depend on $\l$, one can show that 

$$\lim \Psi(v_\l)=v_\l\otimes u_{k+l}.$$ 

Comparing it with the the expression for
$\lim \Phi_{k+l}(v_\l)$, we get the statement of the theorem.
\qed
\enddemo

\proclaim{Corollary 5.2.5}
$$\Phi^{k+l}_\l=\frac{ d_k(\l) d_l(\l)}{ d_{k+l}(\l)}
		\Phi^k_\l*\Phi^l_\l.\tag 5.2.8$$

\endproclaim

So far, we have proved Theorem 5.2.3 only for the case when
$k,l\in \Z_+$. However, it can be generalized. Let us 
consider the space $\tilde U_k=\{(x_1\ldots x_n)^k p(x), p(x)\in
\C_q[x_1^{\pm 1},\ldots x_n^{\pm 1}], \- p(x)\text{ is a homogeneous
polynomial of degree 0}\}$ where $k$ is an arbitrary complex number.
Formula (4.2.1) defines an action of  
$U_q\sln$ in $\tilde U_k$. Also, define $u_0^k=(x_1\ldots x_n)^k\in
\tilde U_k$. 

\proclaim{Lemma 5.2.6}

\roster\item  The set of weights of $\tilde U_k$ coincides with the
weight lattice $Q$, 
and each weight subspace is one-dimensional. In particular,
$\tilde U_k[0]\simeq \C_q u_0^k$. 

\item  For generic $k$, the mapping 

$$x^\l\mapsto 
\frac{\Gamma_q(\l_1+1)\ldots \Gamma_q(\l_n+1)}
     {(\Gamma_q(k+1))^n} 
	x^{-1-\l},\tag 5.2.9  $$ 
defines an isomorphism $\tilde U_k^\omega \simeq \tilde U_{-1-k}$.
The normalization is chosen so that $u_0^k\mapsto u_0^{-1-k}$. 
Here $\Gamma_q(\l)$ is $q$-gamma function:

$$\Gamma_q(x)= 
\frac{1}{(1-q^2)^{x-1}}\prod_{n=0}^\infty
\frac{1-q^{2(n+1)}}{1-q^{2(n+x)}},$$
so $\Gamma(x+1)=q^{x-1}[x] \Gamma(x)$, where

$$[x]=\frac{q^x-q^{-x}}{q-q^{-1}}.$$

Note that the factors of the form $(1-q^2)^\l$ in the product in
\rom{(5.2.9)} cancel, and thus we can 
consider this product   as a formal power  series in $q$ with
coefficients which
are rational functions in $q^\l, q^k$ \rom{(}which we consider as
independent variables\rom{)}.

\item  If $k\in \Z_+$ then $\tilde U_k$ contains a finite-dimensional
submodule, isomorphic to the module $U_k$ defined above: 
$U_k=\tilde U_k\cap \C_q[x_1, \ldots, x_n]$. Also, in this case
$\tilde U_{-1-k}$ has a finite-dimensional quotient
$U^{k}=\tilde U_{-1-k}/( x^\l \text{ such that at least}$
$\text{ one }\l_i\in
\Z_+)$. In particular, $U_{-1}$ can be projected onto $U^0\simeq \C_q$. 
Moreover, formula \rom{(5.2.9)}  above defines an isomorphism 
$U_k^\omega\simeq U^k$ for $k\in \Z_+$. 
\endroster
\endproclaim

\demo\nofrills{Proof }{} of this  lemma is straightforward.\enddemo

Now, let us assume that $\l$ is generic and 
 consider an intertwiner $\Phi^k_\l:M_\l\to M_\l\hat \otimes
\tilde U_k$ such that $\Phi^k_\l(v_\l)=v_\l\otimes u_0^k+\ldots$, and $\hat
\otimes$ is a tensor product completed with respect to $\rho$-grading in 
$M_\l$.
 Note that
if $k\in \Z_+$ then image of $\Phi(v_\l)$ lies in the submodule
$M_\l\otimes U_k$ (which follows from the explicit
formula (5.2.2) for $\Phi$), so this is consistent with our
previous notations. Also, for $k\in\C$ we define 

$$ d_k(\l)=\prod_{\a\in R^+}\prod_{i=0}^\infty\frac
	{1-q^{2(\a, \l+\rho)+2i}q^{-2k}}
	{1-q^{2(\a, \l+\rho)+2i}},\tag 5.2.10$$
which we again consider as a formal power  series in $q$ with
coefficients which
are rational functions in $q^\l, q^k$.
 Note that if $k\in \Z_+$, this coincides with
previously given definition. 

\proclaim{Theorem 5.2.7} For any $k\in \Z_+, l\in \C$ we have 
$$\Phi^k_\l*\Phi^l_\l=\frac{ d_{k+l}(\l)}
			{ d_k(\l) d_l(\l)}
\Phi^{k+l}_\l=\Phi^l_\l* \Phi^k_\l,\tag 5.2.11$$
where $ d_k(\l)$ is given by formula \rom{(5.2.10)}.
\endproclaim

\demo{Proof} Let us fix $k$. Then the matrix elements of the operators
on both sides of (5.2.11) are rational functions in $q, q^l, q^\l$, which
follows from the fact that $ \frac{ d_{k+l}(\l)}
			{ d_k(\l) d_l(\l)}$ is a rational function in
$q, q^l, q^\l$. 
Now the statement of the theorem
follows from the fact that this is
true for $l\in \Z_+$ and the following trivial statement:
\proclaim\nofrills{}
If $F(q,t)\in \C(q,t)$ is such that $F(q, q^l)=0$ for all $l\in\Z_+$ then
$F=0$. \endproclaim
\qed\enddemo

Let us apply this to case when $l=-1-k$. In this case explicit
calculation gives the following answer: 

\proclaim{Corollary 5.2.8} 
For $k\in \Z_+$, 
$$\Phi^k_\l*\Phi^{-1-k}_\l=
\Phi^{-1}_\l\prod_{\a\in R^+}\prod_{i=1}^k 
\frac{1-q^{2(\a, \l+\rho)+2i}}
     {1-q^{2(\a, \l+\rho)-2i}}.
$$
\endproclaim

Now, let us relate these intertwiners with intertwiners for
finite-dimensional representations. 

\proclaim{Theorem 5.2.9} 
Let $\l, \mu\in P^+, u\in L_\mu[0], u\ne 0$ \rom{(}note that  in this
case $L_\l$ is
finite-dimensional\rom{)}. Assume that $\l$ is such
that the intertwiner $\Phi^{\mu, u}_\l:M_\l\to M_\l\o L_\mu$ defined in
Lemma~\rom{5.2.1} is well-defined at this point, i.e. does not have a
pole. Let $I_\l$ be the  maximal submodule in $M_\l$:
$L_\l=M_\l/I_\l$. Then $\Phi^{\mu, u}_\l(I_\l)\subset I_\l\o L_\mu$, and thus, 
$\Phi^{\mu, u}_\l$ can be considered as an intertwiner 
$L_\l\to L_\l\o  L_\mu$. Moreover, this is the unique 
 intertwiner $L_\l\to L_\l\o L_\mu$ such
that $v_\l\mapsto v_\l\o u +\lot$. 

 \endproclaim

\remark{Remark} Note that for $\l\in P^+$ 
the intertwiner $M_\l\to M_\l\o L_\mu$ such
that $v_\l\mapsto v_\l\o u+\lot$ is not unique, so $\Phi^{\mu, u}_\l$ is
a very special intertwiner of this form.\endremark

\demo{Proof} Composing $\Phi^{\mu, u}_\l$ with the projection
$M_\l\to L_\l$, we get an intertwiner $M_\l\to L_\l\o L_\mu$. Since
the tensor product is finite-dimensional, this intertwiner must
annihilate $I_\l$. Uniqueness can be proved in the same way 
as for $M_\l$ for generic $\l$ (see the proof of Lemma~5.2.1).
\qed \enddemo

This shows that we can now apply the results of this section to study of 
products  of intertwiners for finite-dimensional representations. 
In particular, we can now prove Theorem~5.1.3 form the previous section.
Recall that it was formulated as follows:

\proclaim{Theorem 5.1.3} Let $U=U_{k-1}$ and let $\Phi_\l: L_\lk\to L_\lk\o
U, \Phi^{\circ}_\l:L_\lk\to L_\lk \o U^\omega$ be such that 
$\Phi_\l(v_\lk)=v_\lk\o u_0^{k-1}, \Phi_\l^{\circ}(v_\lk)= v_\lk\o
(u_0^{k-1})^\omega$, where $\lk=\l+(k-1)\rho, u_0\in
U[0],(u_0^{k-1})^\omega
\in U^\omega[0]$. Let $(\,
,\,)_U: U\o U^\omega \to \C_q$ be Shapovalov form normalized so that
$(u_0, u_0^\omega)_U=1$. Then 

$$
(\Phi_\l\Phi_\l^{\circ})_U= \prod_{\a\in R^+}\prod_{i=1}^{k-1} 
\frac{1-q^{2(\a, \l+k\rho)+2i}}
     {1-q^{2(\a, \l+k\rho)-2i}}\Id_{L_{\lk}}.
$$
\endproclaim
\demo{Proof} Recall the notation $\tilde U_k$ (see Lemma 5.2.6).
 Then $U=U_{k-1}$ is a
submodule in $\tilde U_{k-1}$, and $U^\omega$ is a factormodule of
$\tilde U_{-k}$; moreover, the Shapovalov form coincides with the
restriction on $U\o U^\omega$ of the map $\tilde U_{k-1} \o \tilde
U_{-k}\to U_{-1}\to \C_q$. This together with Corollary~5.2.8 and 
Theorem~5.2.9 proves the desired statement. 

\enddemo

\head 5.3 Symmetry identities\endhead

In this section we only consider the case $\g=\sln$.

The main goal of this section is to prove Theorem 5.3.3, which
establishes certain symmetry between the values $P_\l(q^{2(\mu+k\rho)})$
and $P_\mu(q^{2(\l+k\rho)})$ (notations will be explained later). The proof
of this theorem is based on the technique of ribbon graphs. 

As before, let $\Phi_\l: L_\lk\to L_\lk\otimes U_{k-1}$ be such that
$\Phi(v_\lk)=v_\lk\otimes u_0^{k-1}+\ldots$, $\lk=\l+(k-1)\rho$, and
let $\varphi_\l= \chi_{\Phi_\l}$ be the corresponding generalized
character. 

\proclaim{Lemma 5.3.1} 

$$\psfig{figure=ident4.ps}
\vbox{$= \varphi_\mu(q^{2(\l+k\rho)})\Phi_\l.$
\vskip 1cm}\tag 5.3.1$$   
where $\varphi_\mu(q^\l)$ stands for polynomial in $q,q^{-1}$
which is obtained by replacing each formal exponent $e^\a$ in the
expression for $\varphi_\mu$ by $q^{(\a,\l)}$. 
\endproclaim

\demo{Proof} Let us consider the operator $F:L_\lk\to L_\lk\otimes
U_{k-1}$ corresponding to the ribbon graph on the left hand side of
(5.3.1). It is some $\Ug$-homomorphism.
 Since we know that such a homomorphism is unique up to a
constant, it follows that $F=a\Phi_\l$ for some constant $a$. To find
$a$, let us find the image of the highest-weight vector. 
First, consider the following part of this picture: 

\psfig{figure=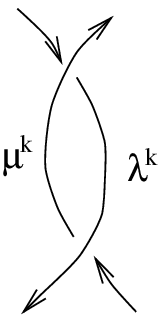}

The corresponding operator is the product $R^{21}R: L_\lk\o 
(L_{\mu^k})^*\to L_\lk\o (L_{\mu^k})^*$. It follows from the 
explicit form of R-matrix (1.2.7) that if $x\in L_{\mu^k}^*[\a]$
then $R^{21}R (v_\lk\otimes x)=q^{-2(\lk, \a)}v_\lk\otimes x+\ldots$. Thus,
if $x_i$ is basis in $L_{\mu^k}$, $x^i$ -- dual basis in $L_{\mu^k}^*$, 
 $x_i$ has weight $\a_i$ then explicit calculation shows that

$$\psfig{figure=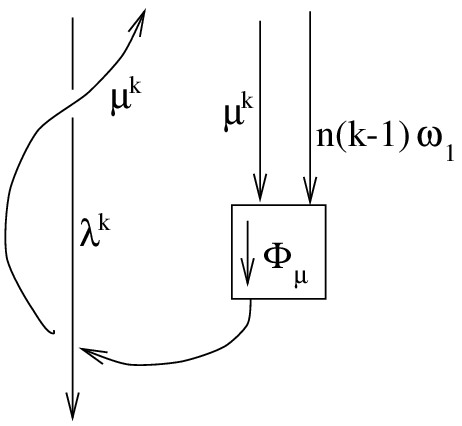}
\vbox{
$:v_\lk\mapsto\sum_i q^{2(\a_i, \l+k\rho)} 
   v_\lk\otimes x^i\otimes \Phi(x_i)+ l.o.t.$\
\vskip 1cm}$$
and thus, $F(v_\lk)=\varphi_\mu(q^{2(\l+k\rho)})v_\lk\otimes
u_0^{k-1}+\ldots$, which completes the proof.
\qed\enddemo

\remark{Remark} In the case $k=1$, i.e. $U=\C$, (5.3.1) reduces to the
formula for the value of the central element $c_{L_{\mu^k}^*}$ in
$L_\lk$ (see Theorem~1.2.1)\endremark

\proclaim{Corollary 5.3.2} Let $\Phi_\l:L_\lk\to L_\lk\o U_{k-1},
\Phi_\l^\circ:L_\lk\to L_\lk\o U_{k-1}^\omega$ be as in 
Theorem~\rom{5.1.2}. Then 

$$\psfig{figure=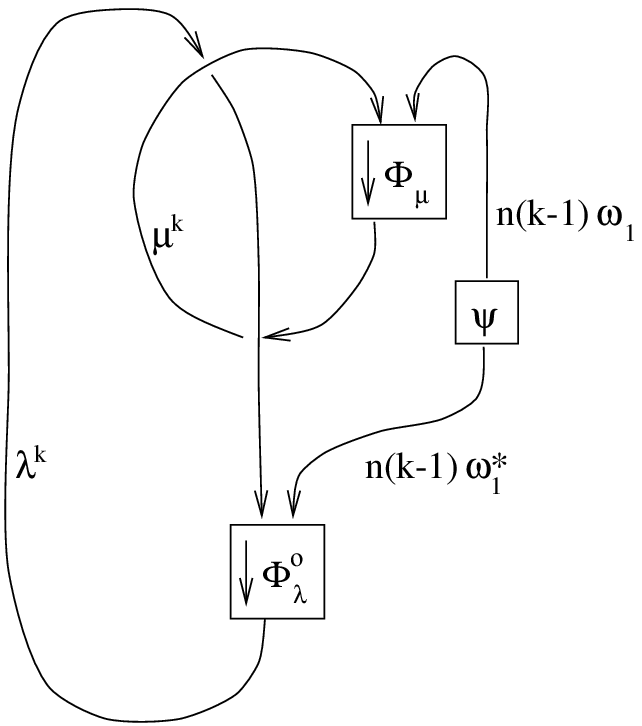}
\vbox{
$=\varphi_\mu(q^{2(\l+k\rho)})\<P_\l, P_\l\>_k\dim_q L_\lk .$\
\vskip 3cm} 
\tag 5.3.2$$
\endproclaim

\demo{Proof} This follows from the previous lemma and the arguments
used in the proof of Theorem~5.1.2.\qed\enddemo

In a similar way, we can replace everywhere $U$ by $U^*$ (see Remark
4.2.4), and 
 repeating with necessary changes all the steps of
previous sections prove the following theorem:

\proclaim{Theorem  5.3.3} Formula \rom{(5.3.2)} remains valid if we replace in
the graph on the left hand side $\Phi_\mu$ by $\Phi_\mu^\circ$,
$\Phi_\l^\circ$ by $\Phi_\l$ and interchange $\omega_1$ and
$\omega_1^*$. \endproclaim

\proclaim{Theorem 5.3.4}
$$\frac
	{P_\mu(q^{2(\l+k\rho)})}
	{P_\l(q^{2(\mu+k\rho)})}
= q^{2k(\rho, \l-\mu)}\prod_{\a\in R^+} \prod_{i=0}^{k-1}
\frac 
	{1-q^{2(\a,\mu+k\rho)+2i}}
	{1-q^{2(\a, \l+k\rho)+2i}}
=\prod_{\a\in R^+}\prod _{i=0}^{k-1}
\frac{ [(\a, \mu+k\rho)+i]}
     { [(\a, \l+k\rho)+i]}.	
				\tag 5.3.3$$
\endproclaim

\demo{Proof} The proof is based on the following identity of the ribbon
graphs: 

$$\psfig{figure=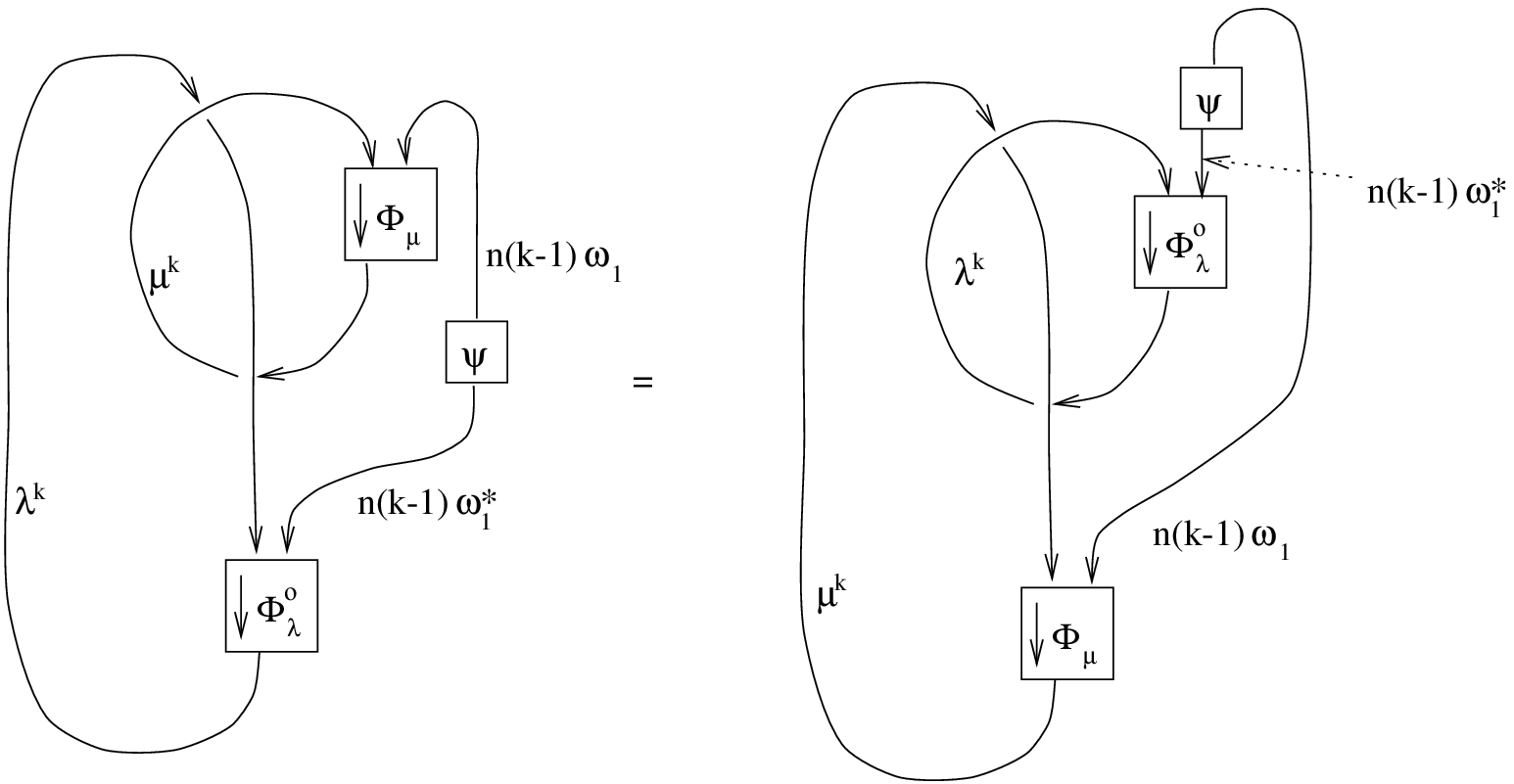} 
\tag 5.3.4$$

Due to Corollary 5.3.2 and Theorem 5.3.3, this implies:

$$\varphi_\mu(q^{2(\l+k\rho)})\<P_\l, P_\l\>_k\dim_q L_\lk=
	\varphi_\l(q^{2(\mu+k\rho)})\<P_\mu, P_\mu\>_k
\dim_q L_{\mu^k}.\tag 5.3.5$$

Substituting in this formula explicit expression for $\<P_\l,
P_\l\>_k$ (Theorem~5.1.4) and using the following   expression for $\dim_q
L_\l$: 

$$\dim_q L_\l=q^{-2(\l, \rho)}
\prod_{\a\in R^+} \frac{1-q^{2(\a, \l+\rho)}}
					{1-q^{2(\a, \rho)}}
=\prod_{\a\in R^+}\frac{[(\a, \l+\rho)]}
			{[(\a, \rho)]}
,$$
which can be easily deduced from the Weyl character formula, 
 we get the statement of the theorem. 
\qed\enddemo

\proclaim{Corollary 5.3.5} {\rm (Macdonald's special value identity,
\cite{M1,M2})}
$$P_\l(q^{2k\rho})=q^{-2k(\rho, \l)}
	\prod_{\a\in R^+}\prod_{i=0}^{k-1}
	\frac
		{1-q^{2(\a, \l+k\rho)+2i}}
		{1-q^{2(\a, k\rho)+2i}}
=\prod_{\a\in R^+}\prod_{i=0}^{k-1}
	\frac{[(\a, \l+k\rho)+i]}
	     {[(\a, k\rho)+i]}
, \tag 5.3.6$$
\endproclaim

\demo{Proof} Let $\mu=0$. Then $P_\mu=1$, and formula (5.3.3) reduces to 
(5.3.6).\qed\enddemo

We refer the reader to \cite{EK5} for the generalization of this formula
for arbitrary  $k$.


\vfill\newpage

\vbox{\vskip 1in}
\specialhead\chapter{6} 
\centerline{GENERALIZED CHARACTERS FOR AFFINE LIE ALGEBRAS} 
\endspecialhead

\vbox{\vskip 0.5in}

    In this chapter we start the study of the generalized
characters for the affine Lie algebras. Mostly, the
constructions are parallel to the finite-dimensional case,
but since the objects we work with are infinite-dimensional,
extra care should be taken.

\head 6.1 Affine Lie algebras\endhead

Here  we review  the notations and facts about affine
Lie algebras and root systems. All of them can be found in \cite{Ka1}.
We keep the notations of Chapter~1; as a rule, we will use 
hat (\ $\hat{}$ ) in the notations of affine analogues of
finite-dimensional objects.  

Let $\ghat$ be the affine Lie algebra corresponding to $\g$:

$$\ghat=\g\o \C[t,t^{-1}]\oplus\C c \oplus\C d,$$
with the commutation rule given by 

$$\gathered
[x\o t^n,y\o t^m]=[x,y]\o t^{m+n}+n\delta_{m,-n}(x,y)c,\\
c \text{ is central},\\
[d,x\o t^n]=nx\o t^n.\endgathered\tag 6.1.1$$

Sometimes we will use a smaller algebra $\gtilde=\g\o
\C[t,t^{-1}]\oplus\C c $. For brevity, we will use the
notation 

    $$x\o t^n =x[n], \quad x\in \g, n\in \Z.$$

Similarly to the finite-dimensional case, we define Cartan subalgebra 
$\hhat=\h\oplus\C c\oplus \C d$, $\hhat^*=\hhat^*\oplus\C
\delta\oplus\C\eps$, where $\<\eps, \h\oplus\C d\>=\<\delta, \h\oplus
\C c\>=0$, $\<\delta,d\>=1,\<\eps, c\>=1$. 
It will be convenient to consider affine hyperplanes
$\hhat^*_K=\h^*\oplus\C \delta +K\eps, K\in \C$; we will refer to the
elements of $\hhat^*_K$ as having level $K$.  

Again, we have a bilinear non-degenerate symmetric form $(\cdot,
\cdot)$ on $\hhat^*$ which coincides with previously defined  on
$\h^*$  and $(\eps, \delta)=1, (\eps,\h^*)=(\delta,
\h^*)=(\eps,\eps)=(\delta,\delta)=0$. 
This gives an identification $ \hhat^*\simeq \hhat: \l\mapsto h_\l$,
and a bilinear form on $\hhat$ such that $(c,d)=1$. Note that under
this identification, $Q\v\subset Q$. 

We define the affine root system $\Rhat=\{\hat\alpha=\alpha+n\delta|\alpha\in
R,n\in\Z \text{ or }\alpha=0, n\in\Z\setminus\{0\}\}$. Again, we have
the notion of positive roots: $\Rhat^+=\{\hat\alpha=\alpha+n\delta\in
\Rhat|n>0 \text{ or } n=0, \alpha\in R^+\}$ and the basis of simple
roots $\alpha_0=-\theta+\delta,\alpha_1,\ldots,\alpha_r$.
For $\ahat=\a+n\delta\in \Rhat^+$ we define
$e_\ahat=e_\a[n]$ if $\a\in R^+$, and $e_\a=f_{-\a}[n]$ if
$\a\in R^-$; $f_\ahat$ are defined similarly. As before, we
use the notations $e_i=2e_{\a_i}/(\a_i, \a_i),
f_i=2f_{\a_i}/(\a_i, \a_i), h_i=\a_i\v =2h_{\a_i}/(\a_i, \a_i),
i=0,\dots, r$. 
This gives the polarization of $\ghat$: $\ghat=\hat \frak
n^+\oplus \hhat \oplus\hat\frak n^-$.

As usual, we denote $\hat Q=\bigoplus \Z\a_i, \hat Q^+=\bigoplus \Z_+
\a_i$.

We define the affine Weyl group $\What$ as the group of
transformations of $\hhat^*$ generated by the reflections with respect
to $\alpha_i, i=0\ldots r$.  We have
notion of sign of an element of $\What$: $\varepsilon(w)=(-1)^l$ if $w$ is a
product of $l$ reflections. This group  preserves the bilinear form;
also, it preserves each of the affine hyperplanes $\hhat^*_K$.

\proclaim{Theorem 6.1.1}{\rm (see \cite{Ka1})}
 $\What\simeq W\ltimes Q\v$, where the action of $W$ is the same as in
the classical case, and the action of $Q\v$ in $\hhat^*_K$ is given by

    $$\alpha\v\colon \lhat \mapsto \lhat+ K\alpha\v-
\bigl(\<\lhat,\alpha\v\> 
	+\frac{1}{2} K(\alpha\v,\alpha\v)\bigr)\delta.
\tag 6.1.2$$
\endproclaim

Now we can define the root lattice $\Phat=P\oplus\Z\delta\oplus
\Z\eps\subset \hhat^*$ and the cone of dominant weights $\Phat^+= \{
\lhat\in\hhat^*| \<\lhat,\alpha\v_i\>\in\Z_+, i=0,\ldots,r\}$. We will
also use the notation $\Phat^+_K=\Phat^+\cap
\hhat^*_K=  \{\lambda+n\delta+K\eps| \lambda\in P^+,
(\l,\theta)\le K\}$. Note 
the cone of dominant weights is invariant with respect to the translations
along $\delta$ direction, but if one factors this out then there is
only a finite number of dominant weights for every level $K$. Abusing
the notations, we will write $P^+_K= \{\lambda\in P^+|
(\l,\theta)\le K\}$. We introduce an order on $\Phat$ as usual: 
$\lhat\le \hat\mu\iff\hat\mu-\lhat\in \hat Q^+$.

  We define the following affine
analogue of $\rho$: 

    $$\hat\rho=\rho+h\v\eps;\tag 6.1.3$$ 
then $\<\hat\rho,
\alpha\v_i\>=1, i=0,\ldots,r$ and thus $\hat\rho\in
\Phat^+$.

\proclaim{Lemma 6.1.2} 
\roster\item $\What$ preserves each $\Phat_K=\Phat\cap\hhat^*_K$.

\item $\Phat^+_K$ is a fundamental domain for the action of $\What$ in
$\Phat_K$ for $K>0$.
\endroster
    \endproclaim

\demo{Proof} The proof is based on the fact that when we restrict the
bilinear form $(\ ,\ )$ to $Q\v$ we get an integer even lattice  --
see \cite{Lo}.
\enddemo


Let us consider representations of $\ghat$. Unless otherwise
specified, we only consider representations with weight
decomposition. We say that a module $V$ over $\ghat$ is of level $K$
if $c|_V=K\Id_V$.

Since we have polarization of $\ghat$, we can define Verma modules
$M_{\lhat}$ over $\ghat$ and their irreducible quotients $L_\lhat$.
Note that $L_\lhat$ is infinite-dimensional unless $\lhat=a\delta$.  
Both of these modules admit weight decomposition. 

    As usual, we
define the category $\Cal O$ to be the
category of finitely-generated  $\ghat$-modules 
with weight decomposition such that for every $v\in V$ the
space $U\hat\frak n^+ v$ is finite-dimensional. Obviously,
for every $\lhat$ the modules $M_\lhat, L_\lhat$ satisfy this
condition. 

    We say that a  module $V$ from category $\Cal O$ is integrable
if for every simple root $\a_i, i=0,\dots, r$ restriction of $V$ to the
Lie algebra $\sltwo$ generated by $e_i, f_i, h_i$ is a
direct sum of finite dimensional modules. Integrable
modules are natural analogues of finite-dimensional modules
for $\g$. It is known that $L_\lhat$ is integrable iff
$\lhat\in \Phat^+$, and that every irreducible integrable
module has the form $L_\lhat$ for some $\lhat\in \Phat^+$.

We will also need another type of modules,
which are called evaluation representations. Let $V$ be a
finite-dimensional module over $\g$ and let  $z$ be a non-zero complex
number. Then we can construct an evaluation representation of $\gtilde$
(not $\ghat$!) in $V$ by 

    $$\gathered
\pi_{V(z)}(a[n])=z^n\pi_{V}(a),\\
\pi_{V(z)}(c)=0.\endgathered\tag 6.1.4$$

Note that $V(z)$ has no $\Phat$-grading but has a natural 
$P$-grading and $V(z)$ is a module of level zero.

    We will be interested in intertwining operators of the
following form, which
are sometimes called vertex operators: 

$$\Phi\colon L_{\lhat}\to \widehat{L_{\lhat}}\otimes V(z),\tag 6.1.5$$
where $\widehat{L_\lhat}$ is the completion of the integrable highest-weight
module $L_\lhat, \lhat\in \Phat^+$ with respect to the $d$-grading. 
To prove the existence of such intertwiners, we use the following
well-known result (see, for example, the arguments in \cite{TK}, which
work for general Lie algebra in the same way as for $\frak s \frak l_2$):

\proclaim{Lemma 6.1.3}
The mapping $\Phi\mapsto (v_{\lhat}, \Phi v_{\lhat})$ establishes
one-to-one correspondence between the space of all intertwiners of the
form \rom{(6.1.5)} and the subspace in $V[0]$ formed by the vectors $v$ such
that 
$x v=0$ for every  $x\in U\hat\frak n^-$ such that 
$xv_\lhat=0$ in $L_\lhat$.
\endproclaim

\head 6.2 Group algebra of the weight lattice \endhead

    One of the most important objects  of study for us will be  the algebra of
$\What$-invariants of the (suitably completed) group algebra of
$\Phat$. In the affine case its definition is much more subtle than in
the finite-dimensional case; our exposition  follows the paper of
Looijenga (\cite{Lo}).

 Let us consider the group algebra $\C[\Phat]$, i.e. the algebra
spanned by the formal exponentials $e^{\lhat}, \lhat\in \Phat$. It is
naturally $\Z$-graded: $\C[\Phat]=\bigoplus\limits_{K\in\Z}
\C[\Phat_K]$. Consider the following completion:

$$\aligned
\overline{\C[\Phat_K]}=\biggl\{\sum_{\lhat\in \Phat_K}
a_\lhat e^\lhat \biggm| &\text{for every $N\in \Z$, there exists only finitely
many $\lhat\in \Phat_K$}\\
&\text{such that $a_\lhat\ne 0$ and
$(\lhat,\hat\rho)\ge N$}\biggr\}.\endaligned\tag 6.2.1$$

Then $\overline {\C[\Phat]}=\bigoplus\limits_K\overline{\C[\Phat_K]}$
is again a $\Z$-graded algebra. This completion is chosen so to
include the characters of Verma modules (and more generally,
modules from category $\Cal O$) over $\ghat$. This algebra has a
natural topology with the basis of neighborhoods of zero given by 

$$X_N=\left\{\sum_{\lhat\in \Phat_K}
a_\lhat e^\lhat\in \overline{\C[\Phat_K]}\biggm| (\lhat, \rhat)<-N\text{ if
} a_\lhat\ne 0\right\}.$$

For our purposes this algebra is too big. 
We will use a smaller algebra: 

    $$A=\bigcap\limits_{w\in
\What} w\biggl(\overline{\C[\Phat]}\biggr), \tag 6.2.2$$
 which is a natural
analogue of the group algebra $\C[P]$ for finite-dimensional case. In
particular, we have a natural action of $\What$ in $A$. It is also
$\Z$-graded: $A_K=\bigcap\limits_{w\in
\What} w\biggl(\overline{\C[\Phat_K]}\biggr)$. Again, this algebra has
a topology, which we describe in terms of convergence: we say that
$f_n\to 0$ in $A$ if for every $w\in \What$, $w(f_n)\to 0$ in
$\overline{\C[\Phat]}$. 

One of the main objects of our study will be
 the algebra of $\What$-invariants  
$A^{\What}\subset A$. We will call elements of
$A^\What_K$ (formal) symmetric theta-functions of level
$K$; indeed, we will show
below that for $\g=\sltwo$ the condition of
$\What$-invariance coincides with definition of even
theta-function.

    For a module $V$ with the weight decomposition define
its character by the formula: $\ch V =\sum_{\lhat\in \Phat}
e^\lhat \dim V[\lhat]$ provided that this sum converges in
the sense of the completion $\overline{\C[\Phat]}$.
    
\demo{Example 1} For any $\lhat\in \Phat_K, K\ge 0$ the orbitsum 
$$m_{\lhat}=\sum\limits_{\hat\mu\in
\What\lhat}e^{\hat\mu}\tag 6.2.3$$ 
belongs to $A^{\What}_K$.\enddemo 

    \demo{Example 2} For any module $V$ from the category $\Cal
O$ (in particular, for $M_\lhat, L_\lhat$) we have $\ch
V\in \overline{\C[\Phat]}$. \enddemo

    \demo{Example 3} If $V$ is a module from the category
$\Cal O$ then $\ch V\in A$ iff $V$ is integrable, in which
case $\ch V\in A^W$. In particular, 
$\ch L_\lhat\in A^\What$ iff  $\lhat\in \Phat^+$.
\enddemo

\demo{Example 4} Define the affine Weyl denominator $\dhat$ by 
  $$\dhat=
e^{\rhat}\prod\limits_{\hat\a\in\Rhat^+}(1-e^{-\hat\a})
=\frac{1}{\ch M_{-\rhat}}.\tag 6.2.4$$

Then $\dhat\in A$. It is easy to see that $\dhat$ is
$\What$-antiinvariant; moreover, it is known (see \cite{Ka1, Lo}) that
every 
$\What$-antiinvariant element
of $A$ has the form $f\dhat, f\in A^{\What}$. \enddemo

Obviously, $A^{\What}$ is $\Z$-graded. Moreover, the following is
well-known: 

\proclaim{Lemma 6.2.1} $A^{\What}_K=0$ for $K<0$, and

    $$A^{\What}_0=\left\{\sum\limits_{n\le n_0}a_ne^{n\delta} ,
a_n\in\C\right\}.\tag 6.2.5$$ \endproclaim

\proclaim{Theorem 6.2.2}{\rm (cf. \cite{Lo})} For every $K\in\Z_+$, the
orbitsums $m_{\lambda+K\eps}, \lambda\in P^+_K$ form a basis of
$A^{\What}_K$ over the field $A^{\What}_0$.\endproclaim  

This theorem follows from the fact that $\Phat^+_K$ is a fundamental
domain for the action of $\What$ in $\Phat_K$ for $K>0$ and from
Lemma~6.2.1. 

It will be convenient to introduce formal variable $p=e^{-\delta}$;
then every element of $A$ can be written as a formal Laurent series in
$p$ with coefficients from $\C[P]$ (note that this is not
so for $\overline{\C[\Phat]}$). It is also compatible with
convergence: if $f_n\to 0$ in $A$ then $f_n\to 0$ in the usual
$p$-adic topology in $\C((p))[P]$.
 In particular, in these notations $A^{\What}_0\simeq \C((p))$. 
(Usually, $e^{-\delta}$ is denoted by $q$; we use the letter $p$ to
avoid confusion with the parameter $q$ of the quantum group).

    We will also need the analogues of the ring
$\C[P](e^\a-1)^{-1}$, considered in Section~2.2. First,
note that if $\ahat\in \Rhat^-$ then the series
$(1-e^\ahat)^{-1}= 1+e^\ahat+e^{2\ahat}+\dots$ converges in
$\overline{\C[\Phat]}$, so there is no need to extend
$\overline{\C[\Phat]}$ any further.

    As for the algebra $A$ defined above, it does not
contain $(1-e^\ahat)^{-1}$. For this reason, we define the
algebra $\hat \Cal R$ as follows. Consider the algebra
$\C[\Phat](1-e^\ahat)^{-1}$, obtained by adjoining to $\C[\Phat]$ the 
inverses of $(1-e^\ahat)$ (no completion so far). Then we have a morphism

    $$\tau\colon \C[\Phat](1-e^\ahat)^{-1}\to
\overline{\C[\Phat]},\tag 6.2.6$$
given by expanding $(1-e^{-\ahat})^{-1}=1+e^{-\ahat}+e^{-2\ahat}+\dots$ for
$\ahat\in \Rhat^+$. Note that the image is not  in $A$. 

Similarly, for every $w\in\What$ we have 

$$\tau_w\colon \C[\Phat](1-e^\a)^{-1}\to
 w\biggl( \overline{\C[\Phat]}\biggr),$$
given by expanding $(1-e^{-\a})^{-1}=1+e^{-\a}+\ldots$ for $\a\in w\Rhat^+$. 

Define 

    $$\multline
\hat\Cal R=\biggl\{\sum a_n \bigm| a_n\in
\C[\Phat](1-e^\a)^{-1},\\
 \sum
\tau_w(a_n) \text{ converges in } w\biggl(\overline{\C[\Phat]}\biggr)
\text{ for every }w\in \What\biggr\}.\endmultline\tag 6.2.7 $$

This algebra obviously contains $A$.   This is the right analogue of the
ring of fractions $\C[P](e^\a-1)^{-1}$, introduced in Section~2.2; 
for example, $\sum_{\ahat\in
\Rhat^+} \frac{1}{1-e^\ahat}\in \hat\Cal R$. Note that there is a natural
action of the Weyl group $\What$ in $\hat\Cal R$; also note that this
algebra has a natural $\Z$-grading given by level.

\head 6.3 Generalized characters for affine Lie algebras\endhead

    In this section, we define and study some properties of
generalized characters for affine Lie algebras. Results of
this section are due to the author and P.~Etingof
(\cite{E1, EK3, EK4}); in less general situation (and different
language) some of the results were found earlier by Bernard
(\cite{Be}).

    \definition{Definition 6.3.1} Let $\Phi:V\to V\o U$ be a
$\gtilde$-intertwiner, where $V$ is a
module from the category $\Cal O$ and  $U$ is a
finite-dimensional $\gtilde$-module of level zero.
Define the corresponding generalized character $\chi_\Phi$
by 

    $$\chi_\Phi=\sum_{\lhat\in \Phat}e^\lhat
\Tr_{V[\lhat]}\Phi \in \overline{\C[\Phat]}\o U[0].\tag 6.3.1$$
    
 \enddefinition

    \demo{Remark} It is easy to see that if $V$ is
integrable then $\chi_\Phi\in A$.\enddemo

    As before, we can interpret generalized characters as
functions of variables $h\in \h, p\in \C$ by the rule 

    $$\chi_\Phi(h,p)=\Tr_V(\Phi p^{-d} e^{2\pi\i h}),\tag 6.3.2$$
provided that this sum converges. In the cases we will be
interested in this sum does converge in a certain region 
(typically, $0<|p|<1$);
however, this requires proof, which we will give later (see Section~8.1). For
this reason, at the moment we consider generalized
characters formally, i.e. as elements of the completed
group algebra of the weight lattice.

    Our next goal is to deduce some differential equations
satisfied by these generalized characters. Unfortunately,
we can not use central elements of $U\ghat$ for this
purpose, since there are no non-trivial central elements in $U\ghat$
except $c$.
However, we can introduce the following central element in
the completion of $U\ghat$:

$$ \widehat C=2(c+h\v)d+\sum_{n\in \Z} 
\left(\sum_{\a\in R^+}(:e_\a[n] f_\a[-n]: +
:f_\a[n]e_\a[-n]:)  + 
\sum :x_l[n]x_l[-n]:\right),\tag
6.3.3$$ 
where $x_l$ is an orthonormal basis in $\h$ with respect to
$(\, ,\,)$ and the normal ordering is defined by 

    $$:x[n]y[m]:=\cases x[n]y[m]\quad m\ge n,\\
		  y[m]x[n]\quad m< n.\endcases\tag 6.3.4$$
    
    This is the affine analogue of the Casimir element.
Though it is an infinite sum, it is well defined in every
module from category $\Cal O$, and $\hat C|_{M_\lhat}=(\lhat,
\lhat+2\rhat)\Id_{M_\lhat}$. 

    Now we can formulate the main result of this section,
which is the affine analogue of Theorem~2.2.4. Similarly to
Section~2.2, let us introduce the following differential
operators, which we consider formally, i.e. as derivations
of the algebra $\C[\Phat]$:

   $$\gathered
\hat\Delta e^\lhat=(\lhat, \lhat)e^\lhat, \\
\d_\ahat e^\lhat=(\ahat, \lhat)e^\lhat, \quad \ahat\in
\hhat^*.\endgathered \tag 6.3.5$$

    If we use the notation $e^{-\delta}=p$ and write 
elements of $\C[\Phat_K]$ as functions of $p$ with
coefficients from $\C[P]$: 
$\C[\Phat_K]=e^{K\eps}\C[p, p^{-1}][P]$ then  

    $$\gathered
\hat\Delta=-2K p\frac{\d}{\d p} +\Delta_\h, \\
\d_{\a+n\delta}=\d_\a+nK,\endgathered
 \tag 6.3.6$$ 
where $\Delta_\h$ is the Laplace operator in $\h$ defined
in Section~2.2.

    Note also that both $\hat\Delta$ and $\d_\ahat$ can be
extended to the completions $\overline{\C[\Phat]}$ and $A$.

    \proclaim{Theorem 6.3.2}{\rm(\cite{E1})}
Let $\lhat\in \Phat^+$, and let  
$\Phi:L_\lhat\to L_\lhat\o
U$ be a $\gtilde$-intertwiner, where $U$ is as in
Definition~\rom{6.3.1}. Denote  by $\chi_\Phi\in A\o U[0]$  the
corresponding generalized character. Then $\chi_\Phi$
satisfies the following differential equation: 

    $$\multline
\left(\hat \Delta - \sum_{\ahat\in \Rhat^+}
\frac{e_\ahat f_\ahat +f_\ahat e_\ahat
}{(e^{\ahat/2}-e^{-\ahat/2})^2}+2\d_{\rhat}-2\sum_{\ahat\in
\Rhat^+} 
\frac{1}{1-e^{\ahat}}\d_\ahat \right)\chi_\Phi\\
 =(\lhat, \lhat+2\rhat)\chi_\Phi,\endmultline \tag 6.3.7$$ 
where we use the following convention: if $\a =n\delta\in \Rhat^+$
then $e_\ahat f_\ahat+ f_\ahat e_\ahat$ should be replaced by $\sum_l
x_l[n] x_l[-n]$, where $x_l$ is an orthonormal basis in $\h$.
Similarly, in the second sum each root should be taken with
multiplicity: for $\ahat=n\delta$ the term
$\frac{1}{1-e^\ahat}\d_\ahat$ must be multiplied by $r$. 

This equation can be rewritten in the following form:
 if $\chi_\Phi$ has level $K$ then 

$$\aligned
&\left(\hat \Delta - \sum_{\ahat\in \Rhat^+}
\frac{e_\ahat f_\ahat +f_\ahat e_\ahat
}{(e^{\ahat/2}-e^{-\ahat/2})^2}\right)(\dhat \chi_\Phi)\\
&=\biggl(\Delta_\h - 2 (K+h\v)p\frac{\d}{\d p} \\
  &\qquad\qquad -\sum_{n\in \Z}
      \biggl( \sum_{\a\in R^+} \frac{p^n e^\a}{(1-p^ne^\a)^2}
			(f_\a[n]e_\a[-n]+ e_\a[-n]f_\a[n])\\ 
&\qquad\qquad\qquad\qquad
	     +\sum _l \frac{p^n}{(1-p^n)^2}x_l[n]x_l[-n]\biggr)
  \biggr) 
  (\dhat \chi_\Phi)\\
&=(\lhat+\rhat, \lhat+\rhat)\dhat\chi_\Phi,\endaligned
\tag 6.3.8$$
where $\dhat$ is the affine Weyl denominator \rom{(6.2.4)}. 

\endproclaim

    In less general  form this equation has first
appeared in the papers of Bernard (\cite{Be}).

\remark{Remark} It is easy to check that the differential operators
written in the right-hand side of (6.3.7), (6.3.8) are well defined as
operators in $\hat \Cal R$.\endremark

\demo{Proof} The proof is analogous to the finite-dimensional
case (Proposition~2.2.4). It is based on the following identities,
which can be easily proved by the same methods as their
finite-dimensional analogues (2.2.5--2.2.7):

$$\gather
\Tr (\Phi e_\ahat f_\ahat p^{-d} e^{2\pi\i h})= 
\left( - \frac{e^\ahat}{(1-e^\ahat)^2}e_\ahat f_\ahat -
\frac{e_\ahat}{1-e^\ahat} \d_\ahat \right)\chi_\Phi,\tag 6.3.9\\ 
\Tr (\Phi f_\ahat e_\ahat p^{-d} e^{2\pi\i h})= 
\left( - \frac{e^\ahat}{(1-e^\ahat)^2}f_\ahat e_\ahat -
\frac{1}{1-e^\ahat} \d_\ahat \right)\chi_\Phi,\tag 6.3.10\\ 
\Tr (\Phi(2 cd +\sum x_l^2)p^{-d} e^{2\pi\i
h})=\hat\Delta\chi_\Phi.\tag 6.3.11\endgather$$

Since $\hat C|_{L_\lhat}=(\lhat, \lhat+2\rhat)\Id_{L_\lhat}$, we have 
 $\Tr (\Phi \hat C p^{-d}e^{2\pi\i h})= (\lhat,
\lhat+2\rhat)\chi_\Phi$. Substituting in this formula the expression
(6.3.3) for $\hat C$ and using identities (6.3.9--6.3.11), we obtain
the desired equation (6.3.7). Formula (6.3.8) can be deduced straightforwardly; it
also follows from more general statement, which we will prove in the
next chapter (Theorem~7.1.2). \enddemo

\example{Example 6.3.3} Assume that $U=U(z)$ is an evaluation
representation. Then $\chi_\Phi\in U[0]$, and thus, independently of
$z$ we have $x_l[n]\chi_\Phi=0, e_\ahat  f_\ahat\chi_\Phi=f_\ahat
e_\ahat \chi_\Phi= e_\a f_\a\chi_\Phi$. Therefore, we can rewrite
(6.3.8) as follows: 

$$\biggl(\Delta_h - 2 (K+h\v)p\frac{\d}{\d p} 
  -2\sum\Sb n\in \Z\\\a\in R^+\endSb 
 \frac{p^n e^\a}{(1-p^ne^\a)^2} e_\a f_\a
  \biggr) 
  (\dhat \chi_\Phi)
=(\lhat+\rhat, \lhat+\rhat)\dhat\chi_\Phi.\tag 6.3.12$$
    \endexample

\remark{Remark} In fact, both Theorem~6.3.2 and Example~6.3.3 are also
true if one replaces the integrable module $L_\lhat$ by Verma module
$M_\lhat$ or any its subfactor, and lets $\lhat$ be an arbitrary (not
necessarily dominant) weight. The only difference is that in this case the
generalized character should be understood as an element of
$\overline{\C[\Phat]}$ rather then $A$. \endremark
    

\vfill\newpage

\vbox{\vskip 1in}
\specialhead\chapter{7} 
\centerline{AFFINE JACOBI POLYNOMIALS AND GENERALIZED CHARACTERS} 
\endspecialhead

\vbox{\vskip 0.5in}

    In this chapter we define and study the affine analogue
of Jacobi polynomials discussed in Chapter~3. This
definitions is due to the author and Etingof \cite{EK4}. 
(The name ``polynomials'' is not quite good, since they
look rather like theta-functions; still, we use it to
stress the analogy with finite-dimensional case).  

    We
also prove that for the root system of type $\hat A_{n-1}$
these polynomials, which in this case we call affine Jack
polynomials, can be obtained as ratio of generalized
characters for the affine Lie algebra
$\widehat{\frak{sl}_n}$.

\head 7.1 Definition of affine Jacobi polynomials. \endhead

    As before, let us fix some positive integer $k$. Recall
that we have defined in the previous chapter 
the algebra $A$, which is a certain
completion of the group algebra $\C[\Phat]$, and the ring
of fractions $\hat \Cal R$, which is a certain completion
of the ring $\C[\Phat](e^\ahat-1)^{-1}$. Recall also the Laplace
operator $\Delta_\h$ defined in Section~2.2.

    \definition{Definition 7.1.1} 
The Calogero-Sutherland operator for (affine) root system $\Rhat$ is
the differential operator which acts in $\hat\Cal R_K$ by
the following  formula (all the notations as before): 

$$\hat L_k= \Delta_\h - 2K p\frac{\d }{\d p}
-k(k-1)\sum\Sb \a \in R^+\\n\in\Z\endSb 
\frac{p^n e^{\a}}{(1-p^ne^{\a})^2}(\a,\a).\tag 7.1.1$$
\enddefinition

    \remark{Remarks} 
\roster \item 
    It is easy to see that the sum does
belong to the ring $\hat \Cal R$, so $\hat L_k$ is a
well defined operator in $\hat \Cal R$. 
\item This operator can be rewritten in the following form:

    $$\hat L_k=\hat \Delta - k(k-1) \sum_{\ahat\in \Rhat^+}
\frac{(\ahat, \ahat)}{(e^{\ahat/2}-e^{-\ahat/2})^2}, \tag 7.1.2$$
which shows that it   is a complete analogue of the operator
$L_k$ introduced in Section~3.1.
This also shows that $\hat L_k$ 
commutes with the action of $\What$ in $\hat\Cal R$.

\item Note that for $K=0$ this operator becomes so-called elliptic
Calogero-Sutherland operator (it becomes clear if one rewrites
$\hat L_k$ in terms of Weierstrass function -- see formula (8.1.6)
below). 
    \endroster
\endremark

Define 

    $$\Mhat_k=\dhat^{-k}\circ(\hat L_k-k^2(\rhat,\rhat))
\circ\dhat^k,\tag 7.1.3$$
where $\dhat$ is the affine Weyl denominator (6.2.4). 

\proclaim{Theorem 7.1.2}
\roster
\item $\Mhat_k$ is a well defined operator in $\hat\Cal R$.

\item 
$$\Mhat_k=\hat\Delta-2k\sum_{\ahat\in\Rhat^+}
\frac{1}{1-e^\ahat}\d_{\ahat} +2k \d_{\rhat}.\tag 7.1.4$$

\item $\Mhat_k$ commutes with the action of $\What$.
\endroster
\endproclaim

\demo{Proof} It will be convenient to use vector fields in $\hhat$, i.e.
the elements of $\hat\Cal R\otimes \hhat$. If $v=\sum v_i\l_i,
v_i\in \hat\Cal R, \l_i\in \hhat$ then we will denote by 
$\d_v=\sum v_i\d_{\l_i}$ the corresponding differential operator. 
Also, for $f\in \hat\Cal R$ denote by $\text{grad }f\in \hat\Cal R
\o \hhat$ its gradient,
which is defined in the usual way, so that $\text{grad } e^\lhat=
e^\lhat \cdot \lhat$ (as before, we identify $\hhat^*$ and $\hhat$). 
Then we have the following obvious formula: 
$$\hat\Delta (fg)= (\hat\Delta f)g + f(\hat\Delta g) + 2 (\text{grad }f, \text{grad
} g)= (\hat\Delta f)g + f(\hat\Delta g) + 2\d_{\text{grad }f }g,$$
  where $(\,,\,)$ is
the inner product of the vector fields, i.e., the inner product on
$\hhat$ extended by $\hat\Cal R$-linearity to $\hat\Cal R\otimes
\hhat$ (it has nothing to do with the inner product on  polynomials!). 
Thus,

$$\hat\delta^{-1}\circ \hat \Delta \circ \hat \delta =
\hat\Delta + 2\d_v + (\hat\delta^{-1}\hat\Delta(\hat\delta)),\tag 7.1.5$$
where 

$$v=\hat\delta^{-1}\text{grad }\hat\delta=
\rhat+\sum_{\ahat\in\Rhat^+}\frac{e^{-\ahat}}{1-e^{-\ahat}} \ahat.\tag 7.1.6$$
Note that $\What$-antiinvariance of $\hat\delta$ implies
$\What$-invariance of $v$. Since $\hat\Delta
(\hat\delta)=(\rhat,\rhat)\hat\delta$, which follows from the
denominator identity for affine root systems,  this proves that
$\hat\delta^{-1}\hat\Delta \hat\delta$ is a well defined operator in
$\hat\Cal R$. 

It is easy to prove by induction that 

$$\hat\delta^{-k}\circ \hat\Delta\circ \hat\delta^k = 
\hat\Delta + k(\rhat, \rhat) + 2k\d_v + k(k-1) \hat\delta^{-1}
(\d_v\hat\delta).$$ 

Obviously, $\hat\delta^{-1}(\d_v\hat\delta)=(v,v)$. 

\proclaim{Lemma 7.1.3} 
$$(v,v)= (\rhat, \rhat)+\sum_{\ahat\in \Rhat^+}
(\ahat,\ahat)\frac{1}{(e^{\ahat/2} -e^{-\ahat/2})^2}.$$
\endproclaim

\demo{Proof} Let us consider $X=(v,v)-(\rhat,\rhat)-\sum_{\ahat\in \Rhat^+}
(\ahat,\ahat)\frac{1}{(e^{\ahat/2} -e^{-\ahat/2})^2}\in\hat\Cal
R$. Obviously, it is $\What$-invariant. Also, from the explicit
expression for $v$ it follows that $X$ has only simple poles.
Consider $X\hat\delta$. It is an element of $\hat\Cal R$
with no poles, thus it is an element of $A$. Also, it is
antiinvariant. Thus, $X\in A$. Obviously, $X\in A_0$; since $X$ is
$\What$-invariant, Lemma~6.2.1 implies that $X\in \C((p))$. 

To complete the calculation, let us write the explicit expression for
$X$; then, let us expand it in a series in $e^{-\ahat},\ahat\in
\Rhat^+$ (i.e., apply the map $\tau$ we used in Section~6.2
to define $\hat \Cal R$) and keep only
the terms of the form $e^{n\delta}$ in the expansion. This gives:

$$X=2\biggl(rh\v - \sum\limits_{\a\in
R^+}(\a,\a)\biggr)\sum\limits_{n=1}^\infty \frac{np^n}{1-p^n}, $$
where $r$ is the rank of $\g$. On the other hand,
for any simple Lie algebra $\g$ we have 

    $$\sum_{\a\in R^+}(\a,\a)= rh\v.$$

    The simplest way to prove it is to consider the action of the Casimir
element $C\in U\g$ in the adjoint representation. On one hand,
$C|_{\g} =2h\v\Id _{\g}$, and thus $\Tr|_{\h}C= 2rh\v$. On
the other hand, it is easy to deduce from the formula
$C=\sum_{\a\in R^+}e_\a f_\a +f_\a e_\a +\sum x_l^2$ that
$\Tr_{\h} C = 2\sum_{\a\in R^+} (\a,\a)$.  
 Thus,  $X=0$. \enddemo

This lemma together with previous results immediately implies statements
1 and 2 of the theorem. Statement 3 follows from $\What$-invariance of
the operator $L_k$. 
\qed\enddemo

Note that Theorem 7.1.2 is a complete analogue of Proposition~3.1.1.
\remark{Remarks}
\roster\item
    This technique is borrowed from \cite{Ma}. 
\item
    In the simply-laced case the identity $\sum(\a,\a)=rh\v$
becomes $\text{dim }\g=r(h+1)$, where $h$ is the Coxeter number for
$\g$. This latter identity is well known and has
a beautiful interpretation in terms of the Coxeter
automorphism (see \cite{Kos}). 
\endroster\endremark

\proclaim{Theorem 7.1.4} 
$\Mhat$ preserves the algebra of $\What$- invariant polynomials: $\Mhat
A^{\What}_K\subset A^{\What}_K$. Moreover, in the basis of orbitsums
 its action is triangular:

$$\Mhat m_\lhat=(\lhat, \lhat+2k\rhat)m_\lhat+\sum\Sb \hat\mu<\lhat\\
\hat\mu\in \hat P^+_K \endSb c_{\lhat\hat\mu}m_{\hat\mu}.\tag 7.1.7$$
\endproclaim

\demo{Proof} Proof is based on the following lemma:

\proclaim {Lemma} Let $f\in A^{\What}, f=e^\lhat +\text{ lower
terms}, \hat\a\in \Rhat^+$. Then 
$\frac{1}{1-e^{\hat\a}}\d_{\hat\a} f$ is 
a well-defined element of $A$ with highest term 
$-(\lhat,\hat\a)e^{\lhat-\hat\a}$. \endproclaim

\demo{Proof of the Lemma} The proof is based on the fact that due to
$\What$-invariance, $f$ contains terms $e^{\hat\mu}$ and $e^{\hat\mu-
\<\hat\mu,\hat\a\v\>\hat\a}$ with equal coefficients, and on explicit
calculation. \enddemo

Now the statement of the theorem follows from the explicit formula
(7.1.4) for $M_k$. 
\qed\enddemo

Our main objective will be the study of the eigenfunctions of action
of this operator in $A^{\What}$, which we will call affine Jacobi
polynomials.  More precisely, let us consider the action of $\Mhat$
in the linear space spanned by $m_{\hat\mu}$ with $\hat\mu\le\lhat$. This
space is not finite-dimensional; however, one can still check that
$\Mhat$ has a unique eigenvector with eigenvalue $(\lhat,\lhat+2k\rhat)$ in
this space (this is based on the affine analogue of
Lemma~3.1.3, and on some standard arguments on convergence). 
Thus, we adopt the following definition:

\definition{Definition 7.1.5} 
Affine Jacobi  polynomials $\hat J_\lhat, \lhat\in\Phat^+$  are
the elements of $A^{\What}$ defined by the following conditions:

\roster 
\item
$\hat J_\lhat=m_\lhat+
\sum\limits_{\hat\mu<\lhat}c_{\lhat\hat\mu}m_{\hat\mu}$.

\item $\Mhat_k \hat J_\lhat=(\lhat, \lhat+2k\rhat)\hat J_\lhat$.
\endroster
\enddefinition

As was said above, these conditions determine $\hat J_\lhat$ uniquely. 
Note that if $\hat\mu=\lhat+n\delta$ then $\hat J_{\hat\mu}=p^{-n}\hat
J_\lhat$. 
Thus, it suffices to consider only the polynomials $J_\lhat$ for
$\lhat=\l+K \eps, \l\in P^+_K$.

\example{Examples} \roster
\item Let $k=0$; then $\hat J_\lhat=m_\lhat$.
\item Let $k=1$; then $\hat J_\lhat=\ch L_\lhat$.
\endroster\endexample

    \proclaim{Theorem 7.1.6} For every $K\in \Z_+$, affine
Jacobi polynomials $\{\hat J_{\l+K\eps}\}_{\l\in P^+_K}$ form a
basis over $\C((p))$ in $A^\What_K$.\endproclaim

    \demo{Proof} It follows from the fact that the
orbitsums form a basis of $A^\What_K$ (Theorem~6.2.2) and
the following elementary lemma the proof of which is left
to the reader: 

    \proclaim{Lemma 7.1.7} Fix $K\in \Z_+$. Let
$A=\{a_{\lhat\hat\mu}\}_{\lhat, \hat\mu\in \Phat^+_K}$ be a matrix  
with complex entries satisfying the following conditions: 
\roster\item 
$a_{\lhat\hat\mu}=0$ unless $\hat\mu\le \lhat$.
\item $a_{\lhat\lhat}=1$.
\endroster

    Then this matrix has an inverse: there exists a unique
matrix $B=\{b_{\lhat\hat\mu}\}$ satisfying the same
conditions such that $AB=\delta_{\lhat\hat\mu}$.
\endproclaim 
\enddemo

\head 7.2 Affine Jack polynomials as generalized characters\endhead

    In this section we consider affine Jacobi polynomials
(see Definition~7.1.5) for the root system $\hat A_{n-1}$,
which we call affine Jack polynomials, and show that they
can be obtained from the generalized characters for the Lie
algebra $\slnhat$. In this section, we only consider
$\ghat=\slnhat$. 

    As before, we fix a positive integer $k$. As in
Section~3.1, let $U$ be the finite-dimensional irreducible
representation of $\sln$ with highest weight
$n(k-1)\omega_1$. Recall that $U[0]$ is one dimensional,
and we have fixed an element $u_0\in U[0]$, thus
identifying $U[0]\simeq \C: u_0\mapsto 1$. Let $U(z)$ be
the corresponding evaluation representation of $\widetilde
\sln$.

\proclaim{Proposition 7.2.1}Let $\mu\in \Phat^+$.  A non-zero intertwiner 
$$\Phi\colon L_{\hat\mu}\to \widehat{L_{\hat\mu}}\otimes U(z)$$
exists iff $\hat\mu=(k-1)\rhat +\lhat, \lhat\in \Phat^+$; if it
exists, it is unique up to a scalar. We will denote such an
intertwiner by $\Phi_\lhat$.
\endproclaim

\demo{Proof} The proof is based on Lemma~6.1.3.\enddemo

Let us consider the corresponding generalized characters
$$\varphi_\lhat=\chi_{\Phi_\lhat}, \quad \lhat\in \Phat^+.\tag 7.2.1$$

    They take
values in $U[0]$, which is one-dimensional, and thus can be considered
as scalar-valued so that 
$$\varphi_\lhat= e^{\lhat+(k-1)\rhat} +\lot.$$

\proclaim{Proposition 7.2.2} For every $\ahat\in \Rhat^+$,
$\varphi_\lhat$ is divisible by $(1-e^\ahat)^{k-1}$ \rom{(}divisibility is
to be understood in the algebra $A$\rom{)}. 
\endproclaim

\demo{Proof}  
This statement is trivially true for $\ahat=n\delta$, since
$(1-p^{-n})$ is invertible in $\C((p))$. Therefore, we can assume that
$\ahat=\a+n\delta, \a\in R$, in which case we can repeat the arguments
used in finite-dimensional situation (see proof of Proposition~4.2.2).
Let us consider  the traces of the form

    $$\varphi_\lhat^F=\sum\limits_{\hat\mu\in \Phat}e^{\hat\mu}
\Tr|_{L_{\lhat}[\hat\mu]}(\Phi_\lhat F) , \tag 7.2.2$$
where $F$ is an arbitrary element of $U\slnhat$.  Let us take $F=f_{\hat\a}
^{k-1}$. Then, using the intertwining property of $\Phi_\lhat$ and the
identity $\Delta (f_\ahat)=f_\ahat\otimes 1+ 1\otimes f_\ahat$, we can
prove by induction that 

$$\varphi_\lhat^F=\frac{f_\ahat
^{k-1}\varphi_\lhat}{(1-e^{-\ahat})^{k-1}}.$$

Since $U[0], U[-(k-1)\a]$ are one-dimensional, we can identify both of
them with $\C$; then 
$f_\ahat^{k-1}\colon U[0]\to U[-(k-1)\a]$ becomes a non-zero constant.
On the other hand, it is easy to see that $\varphi_\lhat^F\in A$, which
proves the proposition.\qed\enddemo

    \proclaim{Theorem 7.2.3} Let $\varphi_\lhat$ be the
generalized character for $\slnhat$, defined by
\rom{(7.2.1)}. Then  

$$\varphi_0=\hat\delta^{k-1},\tag 7.2.3$$
where $\dhat$ is the affine Weyl denominator
\rom{(6.2.4)}. 
\endproclaim

\demo{Proof} Let us consider the ratio $f=\varphi_0/\hat\delta^{k-1}$. It
follows from Proposition~7.2.2 that 
$f\in A$. It has level zero and highest term $1$. Moreover,
similar arguments show that if we twist the order on $\Phat$ by the
action of the Weyl group: $\lhat\ge_w\hat\mu$ if $\lhat-\hat\mu\in
w(\hat Q^+), w\in \What$ then the highest term of $f$ with respect to any such
twisted ordering is still $1$. This is only possible if $f\in \C((p))$. 
To complete the proof,  we have to use the differential equation
for the characters. Indeed, Example~6.3.3 implies that $\varphi_\lhat$
satisfies the following equation:

$$\hat L (\varphi_\lhat \hat\delta)
=(\lhat+k\rhat,\lhat+k\rhat)(\varphi_\lhat\hat\delta),$$
since 
$e_\a f_\a |_{U[0]}=k(k-1) \Id_{U[0]}$.
    
Substituting in this equation $\varphi_0=f(p)\hat\delta^{k-1}$, we see that
$f$ satisfies $\Mhat f=0$. Using   formula (7.1.4)  for
$\Mhat$ we get  $2kh\v p\frac{\d}{\d p}f=0$,
which  is possible only if $f$ is a constant. Comparing the highest terms of
$\varphi_0$ and $\hat\delta^{k-1}$, we get the statement of the theorem.\qed
\enddemo

    Now we can prove the main theorem of this section: 

    \proclaim{Theorem 7.2.4} Let $\varphi_\lhat, \lhat\in
\Phat^+$ be the generalized character for $\slnhat$ defined
by formula \rom{(7.2.1)}, and let $\hat J_\lhat$ be Jack 
\rom{(}Jacobi\rom{)}
polynomial for the root system $\hat A_{n-1}$. Then 

$$\frac{\varphi_\lhat}{\varphi_0}=\hat J_\lhat.$$

\endproclaim
\demo{Proof}
The proof is again quite similar to the finite-dimensional
case (Theorem~3.2.3). 
First, we prove that $\varphi_\lhat/\varphi_0\in
A^{\What}$. Consider the module 
$L=L_\lhat\otimes L_{(k-1)\rhat}$. This  module is
unitary (since both factors are unitary); thus, it is
completely reducible and can be decomposed in a direct sum of the
modules $L_{\hat\mu}$:

$$L=L_{(k-1)\rhat+\lhat}+\sum\Sb \hat\mu\in\Phat^+\\
\hat\mu<\lhat+(k-1)\rhat\endSb N_{\hat\mu}L_{\hat\mu}.$$

This sum is, of course, infinite; however, all  multiplicities
$N_{\hat\mu}$ are
finite. In particular, this implies that the character of this module
belongs to the algebra $A^{\What}$ (that is, $\What$-invariants of completed
group algebra of $\Phat$), so in a certain sense this sum converges. 

Let us construct an intertwiner $\Psi\colon L\to L\otimes U(z)$ as
$\Psi= \Id_{L_\lhat}\otimes \Phi_0$. Consider the corresponding
generalized character $\chi_\Psi$. Then it follows from the
decomposition of $L$  that 

$$\chi_\Psi=\varphi_\lhat+\sum\Sb \hat\mu\in\Phat^+\\ \hat\mu<\lhat\endSb 
a_{\lhat\hat\mu}\varphi_{\hat\mu}.$$

On the other hand, $\chi_\Psi=\varphi_0\ch L_\lhat$.
Dividing both sides by $\varphi_0$ we get 

$$\ch L_\lhat=\sum_{ \hat\mu\in\Phat^+} a_{\lhat\hat\mu}
\frac{\varphi_{\hat\mu}}{\varphi_0},$$
where $a_{\lhat\hat\mu}\ne 0$ only if $\hat\mu\le \lhat$, and
$a_{\lhat\lhat}=1$. 
It follows from Lemma~7.1.7 that we can invert this matrix,
writing
$$\frac{\varphi_\lhat}{\varphi_0}=\sum_{\hat\mu\in \Phat^+}
b_{\lhat\hat\mu} \ch L_{\hat\mu},$$
and the coefficients $b_{\lhat\hat\mu}$ satisfy the same conditions as
$a_{\lhat\hat\mu}$. Thus, $\varphi_{\lhat}/\varphi_0\in A^{\What}$ and has
highest term $e^\lhat$.

We have proved that $\varphi_\lhat/\varphi_0$ satisfies the first 
condition in the definition of  Jack polynomials. Now,
using the differential equation for generalized characters
(Example~6.3.3) along with the fact that 
$e_\a f_\a |_{U[0]}=k(k-1)\Id_{U[0]}$, we get 
$$\hat L (\varphi_\lhat
\hat\delta)=(\lhat+k\rhat,\lhat+k\rhat)
(\varphi_\lhat\hat\delta).$$

 Since
$\varphi_0=\hat\delta^{k-1}$ we can rewrite this as follows:

$$\hat L \biggl(\frac{\varphi_\lhat}{\varphi_0}\hat\delta^k\biggr)
=(\lhat+k\rhat, \lhat+k\rhat)
\biggl(\frac{\varphi_\lhat}{\varphi_0}\hat\delta^k\biggr),$$
which is precisely the definition of affine Jack's polynomials. \qed\enddemo

    
\vfill\newpage

\vbox{\vskip 1in}
\specialhead\chapter{8} 
\centerline{ MODULAR PROPERTIES OF AFFINE JACOBI POLYNOMIALS} 
\endspecialhead

\vbox{\vskip 0.5in}

    In this chapter we consider the affine Jacobi
polynomials defined in the previous chapter from the
analytical point of view. We prove that they define an
analytic function in a certain domain, show their connection
with theta-functions and study their modular properties, following the
paper \cite{EK4}. 

All the constructions of this section are valid for
arbitrary Lie algebra $\g$.

\head 8.1 Functional interpretation of $\C[\Phat]$.\endhead

    So far, we have considered elements of $\C[\Phat]$ and
its various completions and extensions formally. In this
section we discuss the analytic approach. 

    Let us define the following domain

$$Y=\h\times\C\times\H,\tag 8.1.1$$ 
where $\H$ is the upper half-plane: 
$\H=\{\tau\in\C|\text{ Im }\tau >0\}$. Then every element 
$e^\lhat\in \C[\Phat]$ can be considered as a function on
$Y$ as follows: if $\lhat=\lambda+a\delta+K\eps$ then put 

    $$e^\lhat(h,u,\tau)=e^{2\pi\i[\<\lambda,h\>+Ku-a\tau]}.\tag
8.1.2$$

    Note that this in particular implies that
$p=e^{-\delta}$ is given by 

    $$p=e^{2\pi\i\tau}.\tag 8.1.3$$

It is easy to see that if $f\in \C[\Phat]$ then $f(h, u, \tau+1)=f(h,
u, \tau)$, so we can as well consider $f$ as a function of $h, u, p$,
where $p\in \C$ is such that $0<|p|<1$ (the reason for this choice of
possible values for $p$ will be clear later). Note that $f\in
\C[\Phat_K]\iff f(h, u,
\tau)=e^{2\pi\i K u}f(h,0, \tau)$; in this case we say that $f$ is a
function of level $K$. 

Of course, we can't extend this rule to the completion $A$. However,
it turns out that we can extend it to certain elements of
$A$. 

\definition{Definition 8.1.1} Let $f=\sum a_\lhat e^\lhat\in
\overline{\C[\Phat]}$, and let $U\subset Y$. We say that $f$ is an
analytic function on $U$ (or that $f$ converges on $U$) if this sum,
considered as a sum of functions on $Y$ by the above rule absolutely
converges in $U$, and this convergence is uniform on 
 compact subsets in $U$. Similarly, if $f\in \hat \Cal R$ then we say
that $f$ is an
analytic function on $U$ if $\tau(f)$ converges on $U$, where $\tau:
\hat \Cal R\to \overline{\C[\Phat]}$ is the homomorphism obtained by
expanding each $(1-e^\ahat)^{-1}=1+e^\ahat+e^{2\ahat}+\dots, \ahat\in
\Rhat^-$ (see the definition of $\hat\Cal R$ in Section~\rom{6.2}). 
\enddefinition

The following theorem is well known (see, for example, \cite{Ka1,
Chapters 11,13})
\proclaim{Theorem 8.1.2} 
\roster\item 
For every $\lhat\in \Phat^+$, the orbitsum
$m_\lhat$ defined by \rom{(6.2.3)} is an analytic function on $Y$.
\item
The affine Weyl denominator $\hat\delta$ defined by \rom{(6.2.4)}
 is an analytic function in
$Y$, and 

$$\hat\delta(h,u, \tau)= e^{2\pi\i h\v u}
\i  ^{|R^+|}p^{-\frac{\text{dim }\g}{24}}\eta (p)^{r-|R^+|} 
\prod_{\a\in R^+} \theta_1(\<\a, h\>), \tag 8.1.4$$
where the theta \rom{(}$\theta_1$\rom{)} and eta \rom{(}$\eta$\rom{)}
 functions are defined in the Appendix. 
\endroster
\endproclaim

It is easy to rewrite the definition of $\What$-invariance in
analytic terms. Indeed,  
define the action of the affine Weyl group $\What$ on $Y$ so that  
$e^{w\lhat}(h,u,\tau)=e^\lhat(w^{-1}(h,u,\tau))$. One easily
checks that when restricted to the finite Weyl group $W$ this action
coincides with the usual action of $W$ on $\h$ (leaving $u,\tau$
invariant), and the action of $\a\v\in Q\v$ is given by

$$\a\v(h,u,\tau)=(h-\a\v\tau, u+\frac 1 2
(\a\v,\a\v)\tau-\<\a\v,h\>,\tau).$$

This immediately implies the following proposition:

\proclaim{Theorem 8.1.3} Let $f\in A_K$ be such that it converges in
$Y$. Define  $\tilde f(h, \tau)=e^{-2\pi\i K u}
 f (h,u, \tau)$. Let us consider the ring $A^{Q\v}$ of elements
from $A$ invariant with respect to the action of $Q\v\subset \What$.
Then $f\in A^{Q\v}$ iff $\tilde f$ satisfies the
following conditions: for every $\a\v\in Q\v$,  
$$\aligned
\tilde f(h+\a\v,\tau)&=\tilde f(h,\tau),\\
\tilde f(h+\a\v\tau,\tau)&=e^{-2\pi \i K (\half(\a\v,\a\v)\tau +
\<\a\v,h\>)}\tilde f(h,\tau),
\endaligned \tag 8.1.5$$
and $f\in A^\What$ iff $\tilde f$ satisfies \rom{(8.1.5)} and is
$W$-invariant. 
\endproclaim

Functions on $\h\times \H$ satisfying (8.1.5) are usually called
theta-functions of level $K$; thus we see that $A^\What$ is the formal
analogue of the ring of $W$-symmetric theta-functions (cf. \cite{Lo}).

Finally, let us consider the questions of convergence of affine Jacobi
polynomials and generalized characters. Recall  (see Section~2.2) 
that we have introduced the differential
operators $\d_x, x\in \h$ and $\Delta_\h$ acting in functions on $\h$ so
that 

$$\gathered
\d_x e^{2\pi\i \<\a, h\>}=\<\a, x\>e^{2\pi\i \<\a, h\>},\\
\Delta_\h e^{2\pi\i \<\a, h\>}=(\a,\a) e^{2\pi\i \<\a,
h\>}.\endgathered$$

If we choose an orthonormal basis $x_l, l=1,\dots, r$ in $\h$ and
denote by $c_l$ the corresponding coordinates in $\h$: $h=\sum c_l
x_l, c_l\in \C$ then these operators can be rewritten as follows: 

$$\gathered
\d_{x_l}=\frac{1}{2\pi\i} \frac{\d}{\d c_l}\\
\Delta_\h=-\frac{1}{4\pi^2} \sum \frac{\d^2}{\d
c_l^2}.\endgathered\tag 8.1.6$$

Now we can write the
differential operators $\hat L, \Mhat$ defined in Section~7.1 
as usual differential operators in $Y$ with analytic
coefficients:

\proclaim{Proposition 8.1.4} The coefficients of the differential
operators $\hat L, \Mhat$ defined by \rom{(7.1.1), (7.1.5)}
respectively are analytic functions in the region 
$\bigl\{0<\text{\rm Im  } \<\a, h\><\text{\rm Im } \tau$ for all 
$\a\in R^+\bigr\}$, and in this
region the following formulas hold: if we consider them as
operators acting on functions of level $K$ then   

$$\gathered
\hat L_k=\Delta_\h -2Kp\frac{\d}{\d p} + \frac{k(k-1)}{4\pi^2} 
\sum_{\a\in R^+} (\wp(\<\a, h\>)+c(\tau)),\\
\Mhat= \Delta_\h -2(K+kh\v)p\frac{\d}{\d p} -2k \sum_{\a\in R^+}
   \sigma(\<\a, h\>)\d_\a \\
-2kK\sum_{n\ge 1} \left(r\frac{1}{1-p^{-n}}n +
\sum_{\a\in R} \frac{1}{1-p^{-n}e^\a}n\right),\endgathered\tag
8.1.7$$
where the functions $\theta_1(x), \wp(x), \sigma(x)$ are the usual
theta-function, Weierstrass $p$-function and sigma-function,
respectively \rom{(}see Appendix\rom{)},
 and the constant $c(\tau)$ is defined in the Appendix 
\rom{(}formula \rom{A4)}.

\endproclaim
    
\demo{Proof} Explicit calculation.\enddemo

\proclaim{Theorem 8.1.5} For every $\lhat\in \Phat^+$, 
the affine Jacobi polynomial $\hat J_\lhat$ converges in $Y$. 
\endproclaim

    \demo{Proof} Assume $\lhat\in \Phat_K$. Due to
Theorem~6.2.2, we can  write $\hat J_\lhat$ as linear combination of
orbitsums:  $\hat J_\lhat =\sum_{\l\in
P^+_K} f_\l(p) m_{\l+K\eps}$ for some $f_\l \in \C((p))$. Let us
substitute this in the defining equation 

$$\Mhat \hat J_\lhat= (\lhat, \lhat+2k\rhat)\hat J_\lhat.$$

Since
$\hat\Delta(m_\lhat)=(\lhat, \lhat)m_\lhat$, it is easy to deduce that
 $f_\l(p)$ satisfy the following system of  first-order 
differential equations:

$$-2(K+k h\v)p\frac{\d}{\d p}f_\l = \sum_{\mu\in P^+_K}f_\mu
a_{\mu\l},$$ 
where the coefficients $a_{\l\mu}\in \C((p))$ are defined by

$$\Mhat m_{\l+K\eps}=\sum_{\mu\in P^+_K}a_{\l\mu} m_{\mu+K\eps}.$$

It follows from Proposition~8.1.4 that $\Mhat m_{\l+K\eps}$ converges
in the region $0<\text{\rm Im  } \<\a, h\><\text{\rm Im } \tau$; this
is only possible if each  $a_{\l\mu}$ converges in $0<|p|<1$.
 It is well known that this implies that $f_\l(p)$
absolutely converges in the same region. Since $m_\lhat$ converges in
$Y$ (Theorem~8.1.2), we get the statement of the theorem. \qed
\enddemo

    \proclaim{Corollary 8.1.6} For $\g=\sln$, the
generalized characters $\varphi_\lhat$ defined by
\rom{(6.2.1)}, converge in $Y$.\endproclaim

    In fact, it is a general result, which we prove later (see
Chapter~9): if
$L$ is an integrable module over $\ghat$, $U$ -- finite-dimensional
representation of $\g$, and $\Phi:L\to L\o U(z)$ is the
$\gtilde$-intertwiner then $\chi_\Phi$ converges in $Y$.

\head 8.2 Normalized characters and modular invariance\endhead

In this section we study modular properties of affine Jacobi polynomials,
considered -- in accordance with the previous section -- as functions on
$\h\times\C\times \H$ (see 8.1.1). In this section we fix level $K\in
\Z_+$ and assume that $\varkappa=K+kh\v\ne 0$, thus
excluding the trivial case
$K=k=0$.
 Then  $f\in A^\What_K$ can be written as 
$f(h, u, \tau)=e^{2\pi\i K u} f(h, 0, \tau)$. For this reason, we can
consider $f$ as function of only $h, \tau$ without losing information. 

To make our functions modular invariant we need to
introduce some factors of the form $p^t, t\in \Bbb Q$. Thus, we need
to consider slightly more general setting than in
the previous sections. Namely, instead of the weight lattice $\Phat$ we
consider a bigger abelian group $\Pp =P\oplus\C\delta \oplus\Z\eps$.
Also, we can consider
the algebra $A'$ formed by finite sums of the form 
$\sum p^{a_i} f_i, a_i\in \C, f_i\in A$, and
the subalgebra of $\What$-invariants in $A'$ in a manner quite similar to 
that  of the previous section. All the results of Chapter~6 hold with
obvious changes. Again, we can consider elements of $A'$ as functions
of $h, \tau$, repeating all the arguments of Section~8.1; the only
difference is that elements of $A'$ define multivalued functions of
$p$; of course, this ambiguity vanishes if we consider them as
functions of $\tau$. 

Define the normalized analogues of $\rhat$ and $\hat\delta$ as
follows:
$$\gathered
\rhat'=\rhat-\frac{(\rho,\rho)}{2h\v}\delta=
\rho-\frac{(\rho,\rho)}{2h\v}\delta+h\v\eps,\\ 
\hat\delta'=e^{\rhat'} \prod_{\ahat\in\Rhat^+}
(1-e^{-\ahat}).\endgathered\tag 8.2.1$$
 This renormalization is chosen so that $\hat\Delta
\hat\delta'=0$; another reason for this renormalization is that so
defined $\hat\delta'$ possesses nice modular properties. Indeed, if we
consider $\hat\delta'$ as a function then (8.1.4)  implies

$$\hat\delta'(h,u, \tau)= e^{2\pi\i h\v u}
\i  ^{|R^+|}\eta (p)^{r-|R^+|} 
\prod_{\a\in R^+} \theta_1(\<\a, h\>), \tag 8.2.2$$
since it is known that 

$$\frac{(\rho, \rho)}{2h\v}=\frac{\text{dim }\g}{24}$$
(``strange formula'' of Freudenthal--de Vries, see \cite{Ka1, Chapter 12}).

Now, let us define the renormalized operator 

$$\Mhat'=\hat\delta^{\prime -k}\hat L \hat\delta^{\prime k}= 
\Mhat -\frac{Kk(\rho,\rho)}{h\v}.\tag 8.2.3$$

It is easy to see that the definition of affine Jacobi polynomials 
$\hat J_\lhat$
which was given for $\lhat\in \Phat^+$ can be easily extended to
$\lhat=\l+K\eps+a\delta\in \Pp, K\in \Z_+,\l\in P^+_K, a$-- arbitrary. 
In particular, the following choice of $a$ is of special interest for
us: 

\definition{Definition 8.2.1} Let  $K\in\Z_+, \lambda\in P^+_K$.
Define the normalized affine Jacobi polynomial  $J_{\lambda, K}$  by 

$$J_{\lambda, K}=\hat J_{\l+\a\delta+K\eps}, $$
where 
$$a=\frac{k(\rho,\rho)}{2h\v}
-\frac{(\lambda+k\rho,\lambda+k\rho)}{2(K+kh\v)}.\tag 8.2.4$$ 
\enddefinition

Note that for $k=1$ they are precisely the (usual) characters of
integrable highest-weight modules, and the normalization coincides
with that in \cite{Ka1, Chapter 13}.

Note that in the case $\g=\sln$ it follows from Theorem~7.2.4 that 
these polynomials can be defined as
follows: 

$$J_{\l, K}(h, \tau)=\frac{1}{\hat\delta^{\prime k-1}}
\Tr_{L_{\lk, K'}}\biggl(\Phi p^{-d-\frac{c}{24}}e^{2\pi\i h}\biggr), 
\tag 8.2.5$$
where, as before, $\lk=\l+(k-1)\rho, K'=K+(k-1)h\v$, $L_{\lk,
K'}=L_{\lk+K'\eps+a\delta}$, where 
$$a=-\frac{(\lk, \lk+2\rho)}{2(K'+h\v)},$$
and $c$ is the central charge for the action of Virasoro algebra in
$L_{\lk, K'}$: $c=\frac{K'\text{dim }\g}{K'+h\v}$ (see \cite{Ka1}).
This form is quite usual in the conformal field theory (WZW model on
the torus).

It follows from the definition of affine Jacobi polynomials (see
Definition~7.1.5) that the normalized affine Jacobi
polynomials 
satisfy the following differential  equation:

$$\Mhat' J_{\lambda,K}=0.\tag 8.2.6$$

\proclaim{Theorem 8.2.2} The space of solutions of the equation
$\Mhat' f=0$ in $A^{\prime\What}_K$ is finite-dimensional, and the
basis over $\C$ in the  space of solutions is given by the 
normalized affine Jacobi polynomials $J_{\lambda,K}$.
\endproclaim

\demo{Proof} Let $f\in A^{\prime\What}_K$ be such that $\Mhat' f=0$.
Due to Theorem 7.1.6, we can write 
$f=\sum_{\l\in P^+_K} f_\l(p)J_{\l, K}$. Substituting it in
$\Mhat' f=0$, we get $p\frac{\d}{\d p}f_\l =0$, which proves the
theorem.\qed\enddemo

It follows from Theorem~8.1.5 that $J_{\l, K}$ converges on
$\h\times \H$. Thus, we can formulate the analytic version of the
theorem above: 

\proclaim{Theorem 8.2.3} 
For any $K\in\Z_+$, the normalized   affine Jacobi polynomials
$\{J_{\lambda,K}\}_{\lambda\in P_K^+}$ 
form a basis over $\C$ of solutions of the equation $\Mhat'f=0$ in the 
space of $W$-symmetric theta-functions of level $K$ \rom{(}see
conditions \rom{(8.1.5))}. 
\endproclaim

Now let us recall some facts about the modular group and its action. 
Recall that the modular group $\Gamma=SL_2(\Z)$ is generated
by the elements 
$$S=\pmatrix 0&-1\\ 1 &0\endpmatrix, T=\pmatrix 1&1\\ 0
&1\endpmatrix$$ 
satisfying the defining relations $(ST)^3=S^2, S^2T=TS^2, S^4=1$. 
This group acts in a natural way on $Y$ as follows: 

$$\pmatrix a&b\\c&d\endpmatrix (h,u,\tau)=\left(\frac{h}{c\tau+d},
u-\frac{c(h,h)}{2(c\tau+d)}, \frac{a\tau+b}{c\tau+d}\right).\tag 8.2.7$$ 

In particular, 

$$\align
T(h,u,\tau)&=(h, u,\tau+1),\\
S(h,u,\tau)&=\left(\frac h \tau , u-\frac{(h,h)}{2\tau} ,
-\frac{1}{\tau}\right).\endalign$$

Also, for any $j\in \C$ we will  define a right action of $\Gamma$
on functions on $Y$ as follows: if $\alpha=\pmatrix
a&b\\c&d\endpmatrix$ then let

$$(f[\alpha]_j)(h,u,\tau)=(c\tau+d)^{-j} f(\alpha(h,u,\tau)).$$

In fact, this is a projective action, which is related to the
ambiguity in the choice of $(c\tau+d)^{-j}$ for non-integer $j$; to
make it a true  action one must consider a central extension of 
$SL_2(\Z)$; we are not going into details here, only
mentioning that the corresponding cocycle takes values in
the unit circle in $\C$. We will call this
action ``an action of weight~$j$''. This action obviously commutes
with the action of $W$. Moreover, it is well known (and is easy to
check) that for every $K\in \Z_+$, this action preserves the space of
theta-functions of level $K$. 
 
Our main goal will be to find the behavior of the (normalized) affine
Jacobi polynomials under modular transformations. The first result in
this direction is

\proclaim{Theorem 8.2.4} Fix $\varkappa \in \N$.
 Then the space of all solutions
of the equation $\hat L f=0$ in the space of 
functions of level $\varkappa$ on $Y$ \rom{(}i.e., functions $f$
satisfying $f(h,u,\tau)=e^{2\pi\i \varkappa u}f(h,0,\tau)$
\rom{)} is invariant under the action  of $\Gamma$ of
weight $j=\frac r 2\left(1+\frac{k(k-1)h\v}{\varkappa}\right)$. \endproclaim

\demo{Proof} It is easy to see that the operator $\hat L$ is invariant
under the action of $T\in \Gamma$. Thus, to prove the theorem, it
suffices to prove the following formula: 

$$\hat L (f[S]_j)=\tau^{-2} ((\hat L g)[S]_j)-
\frac{1}{2\pi\i\tau}(\varkappa(r-2j)+k(k-1)rh\v) f[S]_j.$$

This is based on formula (8.1.7) for $\hat L$, which we rewrite in the
following form: 

$$\hat L_k=\Delta_\h -2\varkappa p\frac{\d}{\d p}-k(k-1)\sum_{\a\in
R^+} (\a,\a)\varphi(\<\a,h\>, \tau),$$
where 

$$\gather 
\varphi(x, \tau)=\sum_{n\in \Z}\frac{p^ne^{2\pi\i x}}{(1-p^ne^{2\pi\i x})^2}
=\frac{1}{4\pi^2}\d_x^2\log \theta_1(x)\text{\phantom{aasdfhjkk}}\\
\text{\phantom{asdfghjkl}}=-\frac{1}{4\pi^2}\wp(x)+c(\tau)\text{ for
some constant }c(\tau).\endgather$$
(see formula (A4) in the Appendix). 

 Using modular
properties of the theta-function (see, for example, \cite{Mu}), we can
show that $\varphi(\frac x\tau, -\frac 1\tau)= \tau^2\varphi(x,
\tau)+\frac{\i}{2\pi}\tau$. 

Also, it is not too difficult to check that 

$$\hat\Delta \left( f[S]_j\right) = \tau ^{-2}( (\hat\Delta f)[S]_j)
- \frac{\varkappa}{2\pi\i\tau}(r-2j)(f[S]_j),$$
using the expression for $\hat\Delta$ given in Section~8.1

Since $\sum_{\a\in R^+}(\a,\a)=rh\v$ (see proof of Lemma~7.1.3), we get
the desired formula. \qed\enddemo

\proclaim{Corollary 8.2.5} The space of solutions of the equation 
$\Mhat'f=0$, where $f$ is a theta-function of level $K$, is invariant
under the action of modular group with the weight 
$j=-\frac{K(k-1)r}{2(K+kh\v)}$.\endproclaim
\demo{Proof} First recall that the space of $W$-invariant
theta-functions of level $K$ is preserved by this action.
Now the statement of the theorem 
follows from the previous theorem and the following
well known fact (see \cite{Ka1, Chapter 13}):  
$$\hat\delta'[\a]_{r/2}=l(\a)\hat\delta',$$
for some  function (not a character) 
$l:\Gamma\to \C^\times$ such that $l^{24}=1$. 
\enddemo

Since we know that the basis of $W$-invariant solutions is given by
the normalized 
affine Jacobi polynomials, we can rewrite Corollary~8.2.5 in the following
form:

\proclaim{Theorem~8.2.6} Let $K\in \Z_+$. Denote 
$$V_K=\bigoplus_{\l\in P^+_K} \C J_{\l, K}.\tag 8.2.8$$

We consider elements of $V_K$ as functions on $Y$. Then $V_K$ is
preserved by the action of $\Gamma$ with weight
$j=-\frac{K(k-1)r}{2(K+kh\v)}$. In particular, this means that $V_K$ is
naturally endowed with a structure of a projective representation of
$\Gamma$, and  the corresponding cocycle takes values
in the unit circle $S^1=\{z\in\C| |z|=1\}$. 
\endproclaim


\vfill\newpage

\vbox{\vskip 1in}
\specialhead\chapter{9} 
\centerline{CORRELATION FUNCTIONS ON THE TORUS}
\centerline{ AND ELLIPTIC KZ EQUATIONS}
\endspecialhead

\vbox{\vskip 0.5in}

In this chapter we consider generalized characters for affine Lie
algebras with values in a
tensor product of evaluation representations. These characters have
appeared in the (twisted) Wess-Zumino-Witten (WZW) model of conformal field
theory as correlation functions on the torus, and the differential
equation for the generalized characters deduced in Chapter~6 can be
interpreted as describing their dependence on the modular parameter
$\tau$ of the torus.  It turns out that their dependence on the
parameters $z_i$ of evaluation representations can also be described
by a simple differential equation, which was first derived by Bernard
(\cite{Be}). We derive rigorously these equations and study their
monodromies, which has not been done before. We call them elliptic
Knizhnik-Zamolodchikov (KZ) equations,
 the reason being that in the limit $\tau\to\i\infty$ these
generalized characters become the usual correlation functions on the
sphere, i.e. certain matrix coefficients of product of intertwining
operators, and the equations reduce to well known
Knizhnik-Zamolodchikov equations (see \cite{KZ, TK, FR}). 

\head 9.1 Intertwiners and currents\endhead

In this section we briefly review the facts about intertwining
operators, commutation relations and (usual) KZ equation we are going
to use. We refer the reader to \cite{TK, FR} for the proofs.
We keep all the notations of Chapters~6,7. Recall that we have defined 
for every $\lhat\in \hhat^*$  Verma module $M_\lhat$ and irreducible
highest-weight module $L_\lhat$. In this chapter we always assume that
we have fixed level $K$ so that we only consider weights of the form
$\lhat=\l+a\delta+K\eps$. Moreover, we will always assume that
$\l, K$ satisfy the following condition:

$$\l\in P^+, \qquad K\notin \Bbb Q.\tag 9.1.1$$

This condition ensures that $L_\lhat$ is equal to the module over
$\ghat$ obtained by induction from the module $L_\l$ over $\g$.
 Note that  this condition implies that $K+h\v\ne 0$. We discuss 
possible generalizations later.

It is known that the structure of the module $L_{\l+a\delta+K\eps}$ does not
depend on the choice of $a$. It will be convenient for us to fix $a$
as follows: we let $a=-\Delta(\l)$, where

$$\Delta(\l)=\frac{(\l, \l+2\rho)}{2(K+h\v)}.\tag 9.1.2$$

We denote the corresponding irreducible module by $L_{\l, K}$ and
Verma module by $M_{\l, K}$.  The reason for such a choice is that in
this case the action of $d$ can be written in terms of action of
$\gtilde$ as follows (see \cite{Ka1}): 

$$d= -\frac{1}{2(K+h\v)}
\sum_{n\in \Z} 
\left(\sum_{\a\in R^+}(:e_\a[n] f_\a[-n]: +
:f_\a[n]e_\a[-n]:)  + 
\sum :x_l[n]x_l[-n]:\right),\tag 9.1.3$$
(compare with (6.3.3)), where normal ordered product is defined by (6.3.4).

Also, we have defined evaluation representations $V(z)$. It will be
convenient to change this definition, considering $z$ as a formal
variable rather than a complex number, as follows:

\definition{Definition 9.1.1} Let $V$ be a module over $\g$. The
evaluation module $V(z)$ is the representation of $\ghat$ in the space
$V\o z^{-\Delta}\C[z, z^{-1}]$, where $\Delta$ is an arbitrary complex
number,  with the action defined as follows:

$$\gathered
\pi_{V(z)}(x[n])=z^n \pi_V(x),\\
 \pi_{V(z)}(c)=0,
\pi_{V(z)}(d)=z\frac{d}{dz}.\endgathered \tag 9.1.4$$
\enddefinition

Note that so defined evaluation representations are indeed modules
over the full affine Lie algebra $\ghat$, not only $\gtilde$, and they
have weight decomposition: $V(z)=\bigoplus_{\lhat\in \hhat^*}
V(z)[\lhat]$ if $V$ has a weight decomposition as $\g$-module.

We will be interested in the  intertwining operators of the form 

$$\Phi: L_{\l, K}\to L_{\nu, K}\hat \otimes L_\mu(z).\tag 9.1.5$$ 
where the completed tensor product is defined as follows: if $V, W$
are $\ghat$-modules then 

$$(V\hat \o W)[\lhat]=
\biggl\{\sum a_i\o b_i, a_i\in V[\lhat_i], b_i\in W[\lhat'_i]\biggm|
\lhat_i+\lhat_i'= \lhat, (\lhat_i, \rhat) \to  -\infty\biggr\}.\tag 9.1.6$$

\proclaim{Lemma 9.1.2}
\roster \item Let $\Phi$ be a $\gtilde$ intertwining operator of the
form \rom{(9.1.5)}. Then it can be uniquely written in the form 

$$\Phi=\sum_{n\in \Z} \Phi[n]z^{-n-\Delta},\tag 9.1.7$$ 
where $\Delta$ is the complex number used in the definition of
evaluation representation, and $\Phi[n]$ are  $\g$-intertwiners
$L_{\l, K}\to L_{\nu, K}\o L_\mu$ such that 

$$\Phi[0](L_\l)\subset  L_\nu\o L_\mu,$$
where $L_\l$ is the $\g$-submodule in $L_{\l, K}$ generated by the
highest-weight vector, and similarly for $\nu$. 

\item The map $\Phi\mapsto \Phi[0]$ establishes isomorphism 
between the space of $\gtilde$ \rom{(}not $\ghat$\rom{!)} intertwiners of the 
form \rom{(9.1.5)} and the space of $\g$-intertwiners $L_\l\to L_\nu\o
L_\mu$.

\item A $\gtilde$-intertwiner $\Phi$ defined above is a $\ghat$
intertwiner iff $\Delta= \Delta(\l)-\Delta(\nu)$, where $\Delta(\l)$
is defined by \rom{(9.1.2)}.
\endroster\endproclaim

From now on, whenever we consider an intertwining operator of the form
(9.1.5) we assume that $\Delta$ is chosen to be
$\Delta(\l)-\Delta(\nu)$, so that $\Phi$ is a $\ghat$-intertwiner.
We will treat $\Phi$ as a formal power series in $z$. Note that the
condition that $\Phi$ commutes with the action of $d$ can be rewritten
as follows: 

$$\Phi d- (d\o 1)\Phi= z\frac{d}{dz} \Phi.\tag 9.1.8$$

We introduce the following series, quite usual in physical
literature. Let $x\in \g$ be a homogeneous element with respect to the
root decomposition, and let $z$ be a formal variable. Define

$$\tilde J_x(z)=\sum_{n\in \Z} x[n]z^{-n}.\tag 9.1.9$$

Similarly, we define

$$\gathered
\tJ^+_x(z)=\sum_{n<0} x[n]z^{-n}+\cases 0, \quad x\in \frak n^+,\\
					\half x, \quad x\in \h,\\
					x, \quad x\in \frak n^-,
				\endcases\\	
\tJ^-_x(z)=-\sum_{n>0} x[n]z^{-n}-\cases 0, \quad x\in \frak n^-,\\
					\half x, \quad x\in \h,\\
					x, \quad x\in \frak n^+.
				\endcases
\endgathered \tag 9.1.10$$

Thus, $\tJ_x(z)=\tJ^+_x(z)-\tJ^-_x(z)$.  We extend these definitions
from homogeneous $x$ to arbitrary $x\in \g$ by linearity. 
We consider $\tJ(z)$ as a series in $z$
with coefficients from $U\ghat$; similar agreement will be applied to
other series we introduce. 

Usually in the  literature a slightly different form of these currents
and polarization is used, namely: 

$$\gathered
J_x(z)=\sum x[n]z^{-n-1}=J^+_x(z)-J^-_x(z),\\
J^+_x(z)=\sum_{n<0}x[n]z^{-n-1}.\endgathered$$
It is easy to express $\tJ$ in terms of $J$ and vice versa; however,
the form we use will be more convenient in our calculations. 

Note that for every $\ghat$-module $V$ from category $\Cal O$,
$\tJ^\pm_x(z)$ can be considered as an operator $V\to \overline{V}\o
\C[z,z^{-1}]$,
where $\overline{V}$ is the completion defined similarly to
(9.1.6). Moreover, $\tJ^-_x(z)$ is in fact a well-defined operator 
 $V\to V\o\C[z,z^{-1}]$.

We can rewrite the definition of intertwining operator in terms
of commutation relations with currents as follows. Let $\Phi(z)$ be an
intertwiner of the form (9.1.5), and $u\in L_\mu^*$. Define  
$\Phi_u(z):L_{\l, K}\to L_{\nu, K}$ by 
$\Phi_u(z)v= \<\Phi(z)v, u\>$. Then the definition of intertwiner
operator takes the form

$$[x[n], \Phi_u(z)]=z^n\Phi_{x u}(z),\tag 9.1.11$$
where $xu$ denotes the usual action of $x\in \g$ in the dual
representation $L_\mu^*$. Note that 

$$\Phi_{xu}(z)=-\<(1\o x)\Phi(z), u\>.$$

\proclaim{Lemma 9.1.3} Let $\Phi$ be a $\ghat$-intertwining operator
of the form \rom{(9.1.5)}. Then 
we have the following identities of power
series in $z,w$: 

$$[\tJ_x^\pm(w), \Phi_u(z)]=\cases 
		    \frac{w}{z-w}\Phi_{xu}(z), x\in \frak n^+,\\
	\half	    \frac{z+w}{z-w}\Phi_{xu}(z), x\in \h,\\
		    \frac{z}{z-w}\Phi_{xu}(z), x\in \frak n^-.
		\endcases\tag 9.1.12$$

Here  by definition, we let 

$$\frac{1}{z-w}=z^{-1}\sum_{n\ge 0} \biggl(\frac{w}{z}\biggr)^n\tag 9.1.13$$ 
in the formula involving $\tJ^+$ and 

$$\frac{1}{z-w}=-w^{-1}\sum_{n\ge 0} \biggl(\frac{z}{w}\biggr)^n\tag 9.1.14$$ 
in the formula involving $\tJ^-$. 

Often we will write $\frac{1}{z-w}, |z|>|w|$ for the expansion
\rom{(9.1.13)}, meaning by this that we are expanding in positive
powers of $w/z$; similarly, we write  $\frac{1}{z-w}, |w|>|z|$ for
\rom{(9.1.14)}. 
\endproclaim

The proof is straightforward use of (9.1.11) and can be found, for
example, in \cite{FR}.

\proclaim{Theorem 9.1.4}{\rm (\cite{FR})}
 In the notations of Lemma~\rom{9.1.3}, we
have

$$(K+h\v)z\frac{d}{dz}\Phi_u(z)= \sum_a \vdots \tJ_a (z) \Phi_{au} (z)
\vdots +\Phi_{h_\rho u}(z),\tag 9.1.15$$
where $a$ runs over an orthonormal basis in $\g$, and 

$$\vdots \tJ_a(z) \Phi_{au}(z)\vdots
=\tJ_a^+(z)\Phi_{au}(z)-\Phi_{au}(z)\tJ_a^-(z).\tag 9.1.16$$ 
\endproclaim

\demo{Proof} The simplest way to prove it is to substitute in the
commutation relation (9.1.8) between  $d$ and $\Phi_u$ the expression
(9.1.3) for $d$. After some trivial though boring computations, we
get the desired formula (9.1.15). \enddemo

\head 9.2 The correlation functions on the torus and elliptic KZ 
equations. \endhead

In this section we  define the correlation function on the torus. We
keep all notations and conventions of the previous sections. Let
$L_{\l_i, K}, i=0\ldots n$ be a collection of irreducible
highest-weight modules such that
$\l_0=\l_n=\l$, and let $\Phi^i(z_i):L_{\l_i,K}\to
L_{\l_{i-1}, K}\otimes L_{\mu_i}(z_i)$ be intertwining operators. 
Then we can consider
the following ``correlation function on the torus'':

$$\chi_{\l,K} (z_1\ldots z_n,p,h)= \Tr|_{L_{\l, K}}
\left(
   \Phi^1(z_1)\ldots \Phi^n(z_n) p^{-d} e^{2\pi\i h}
\right),
\tag 9.2.1$$
where $h\in \h$. This is nothing but the generalized character
$\chi_{\Psi}$, where
$$\gathered
\Psi=\Phi^1(z_1)\dots \Phi^n(z_n): L_{\l, K}\to L_{\l, K}\o V,\\
V=L_{\mu_1}\o\dots\o L_{\mu_n}.\endgathered\tag 9.2.2$$

So far, we consider $\chi$  as a formal series in $p, z_i$ with values in 
the finite-dimensional space $V$; in fact, it follows from weight
considerations that $\chi$ takes values in the zero-weight subspace
$V[0]$. Later we will study the questions of its convergence.
Note also  that since  $\Psi$ commutes with $d$ we have

$$\sum_i z_i\frac{\d}{\d z_i}\chi =0,$$
thus $\chi$ only depends on the ratios $z_i/z_j$. 
We will write $\pi_i(x)$ for action of $x\in \g$ in
the $i$-th factor $L_{\mu_i}$ in the product $V=L_{\mu_1}\o\dots\o
L_{\mu_n}$. 

These functions are the main objects of study in this chapter. We call
them ``correlation functions on the torus'': readers familiar with
Wess-Zumino-Witten model of conformal field theory immediately
recognize that if we let $h=0$ in (9.2.1) then we get what is known as
$n$-point correlation function on the torus (up to a factor of
$p^{-c/24}$, where $c$ is the central charge of Virasoro).

The following theorem  follows immediately from Lemma~9.1.2(2).

\proclaim{Theorem 9.2.1} 
The space of $\ghat$-intertwining operators $\Psi$ of the form
\rom{(9.2.2)} is isomorphic  to the space of all $\g$-intertwiners 
$L_\l\to L_\l\o V$.\endproclaim

Let us consider the limit of $\chi$ as $p\to 0$ of this correlation function.

\proclaim{Lemma 9.2.2}\roster
\item In the limit $p\to 0$, the trace $\chi$ defined by \rom{(9.2.1)}
has the following asymptotic:

$$\chi\sim p^{\Delta(\l)} \Tr_{L_{\l}}
\left(
   \Phi^1(z_1)\ldots \Phi^n(z_n) e^{2\pi\i h}
\right).\tag 9.2.3$$   

Note that the last trace is taken over the finite-dimensional 
 module $L_\l\subset L_{\l,
K}$, not over the whole module $L_{\l, K}$, and is independent of
$p$. 

\item In the limit $p\to 0, h\to +\i\rho\infty$ so that $pe^{2\pi\i
\<h,\a\>}\to 0$ for all $\a\in R$, $e^{-2\pi\i \<h,\a\>}\to 0$ for all
$\a\in R^+$, the trace $\chi$ has the following asymptotics: 

$$\chi\sim p^{\Delta(\l)} e^{2\pi\i \<h, \l\>}
f(z_1, \dots, z_n), \tag 9.2.4$$
where
$$f(z_1, \dots, z_n)=
\<v^*_{\l, K},
   \Phi^1(z_1)\ldots \Phi^n(z_n)v_{\l, K}\> \tag 9.2.5$$ is a function
of $z_1, \dots, z_n$ with values in $V$.

The limits should be understood in the formal sense \rom{(}as limits of
formal power series\rom{)}. 

\item The function $f$ defined by \rom{(9.2.5)} converges to an
analytic function in the region $|z_1|>|z_2|>\dots>|z_n|$ and satisfies
there  the
following system of differential  equations \rom{(}trigonometric form
of Knizhnik-Zamolodchikov equations\rom{):}

$$(K+h\v)z_i\frac{\d}{\d z_i}f =\biggl(\sum_{j\ne i}
\frac{z_j \Omega^+_{ij}+z_i\Omega^-_{ij}}{z_i-z_j}-\pi_i(h_{\l+\rho})
\biggr)f, \tag 9.2.6$$ 
where 

$$\gathered
\Omega^+=\sum_{\a\in R^+}e_\a\o f_\a + \half \sum_l x_l\o x_l,\\
\Omega^-=\sum_{\a\in R^+}f_\a\o e_\a + \half \sum_l x_l\o
x_l\endgathered$$ 
and 
$\Omega^\pm _{ij}=\pi_i\o \pi_j(\Omega^\pm )$.\endroster
\endproclaim

The proof of the first two statements is quite trivial; as for the
last one, we refer the reader to \cite{FR}.

So far, we considered $\chi_{\l, K}$ as a formal series in $z_i$,
$p$. Let us prove that in fact it defines an analytic function in a
certain region. 
Let $\Cal D$ be the set of all $(z_1, \dots, z_n, h, p), z_i, p\in \C,
h\in \h$ satisfying the following conditions:

$$\aligned
&\qquad  0<|p|<1\\
&\qquad |pe^{2\pi\i \<h,\a\>}|<1 \quad \text{for all }\a \in R\\
&\qquad |e^{-2\pi\i \<h,\a\>}|<1 \quad \text{for all }\a\in R^+\\
&\qquad |z_1|>|z_2|>\ldots >|z_n|>|pz_1|\endaligned \tag 9.2.7$$

Recall   the notations $\d_x, x\in \h, \Delta_\h$, introduced in
Section~2.2 (see also formula (8.1.6)). 

\proclaim{Theorem 9.2.3} \roster\item 
The trace $\chi_{\l, K}(z_i, h, p)$ defined by \rom{(9.2.1)} converges
to an analytical function in $\Cal D$. 

\item This trace satisfies the following differential equation in
$\Cal D$: 

$$(\Delta_\h -2(K+h\v)p\frac{\d}{\d p}
+A)(\chi\dhat)=(\rho,\rho)\chi\dhat,\tag 9.2.8$$ 
where
$$\aligned
A=&-\sum_{i,j=1}^n \pi_i\otimes \pi_j
	\biggl[  2\sum \Sb \a\in R^+\\ m\in \Z\endSb 
		f_\a\otimes e_\a   (z_i/z_j)^m	
	\frac{e^{2\pi\i\<\a,h\>}p^m}{(1-e^{2\pi\i \<\a,h\>}p^m)^2}\\
	&\hbox{\hskip 3cm} +\sum\Sb l=1\ldots r\\ m\ne 0\endSb
		x_l\otimes x_l (z_i/z_j)^m
	\frac{p^m}{(1-p^m)^2}\biggr]\\
=& -\sum_{i,j=1}^n \pi_i\otimes \pi_j
	\biggl[ 2\sum_{\a\in R^+}\phi(-\<\a,h\>,\zeta_i-\zeta_j)
			f_\a\otimes e_\a
	+\sum_l \phi_0(\zeta_i-\zeta_j)x_l\otimes x_l\biggr],\endgathered
\tag 9.2.9$$
where $\zeta_i$ are defined by $z_i=e^{2\pi\i \zeta_i}$, 
the functions $\phi, \phi_0$ are defined in the Appendix and 
$\dhat$ is the affine Weyl denominator given by \rom{(6.2.4) (}see also
\rom{(8.1.4))}. We use the following convention: $\pi_i\o \pi_i(a\o
b)= \pi_i(ab)$. 
\endroster
\endproclaim

\demo{Proof} The proof goes as follows. First,  formula
(9.2.8), interpreted as a an equality of formal series, follows from
more general formula (6.3.8). Next, it is easy to see that in fact all
coefficients of the differential equation (9.2.8) converge in $\Cal
D$. It follows from the Lemma~9.2.2 and well-known facts about
convergence of matrix elements of products of intertwining operators
(see \cite{TK})  that in the limit $p\to 0$ the trace $\chi$
converges for $|z_1|>\dots>|z_n|$.
 These observations along with the standard fact
from the theory of ordinary differential equations (namely, that if a
formal series satisfies a differential equation with analytic
coefficients then this series converges to an analytic solution) prove
that $\chi\dhat$ is an analytic function in $\Cal D$; slight
modification of the above arguments, using (6.3.7), allows to prove
that in fact  $\chi$ converges in $\Cal D$. \enddemo

As we have mentioned before, in the case $n=1$ the trace $\chi$ is independent
of $z$. However, in general case it is not so, and we can study the
dependence of this correlation function on $z_i$, which should
generalize the Knizhnik-Zamolodchikov equation (9.2.6).

\proclaim{Theorem 9.2.4} The product $\chi\dhat$,
where $\chi$ is given by \rom{(9.2.1)} and $\dhat$ is the affine 
Weyl denominator 
\rom{(6.2.4)}, satisfies the following system of equations: 

$$(K+h\v)\frac{\d}{\d \zeta_i} \chi\dhat = 
\left(\sum_{j\ne i} r_{ij}(\zeta_i-\zeta_j) 
	-2\pi\i \sum_l \pi_i(x_l)\d_{x_l} \right)\chi\dhat,\tag 9.2.10$$
where 

$$r(\zeta)=2\pi\i\left[\sum_{\a\in R^+} \bigl(e_\a\otimes f_\a \cdot
g(\zeta, \<\a, h\>)+
f_\a\otimes e_\a \cdot g(\zeta, -\<\a,h\>)\bigr) - \sum x_l\otimes x_l
\cdot \sigma(\zeta)\right],\tag 9.2.11$$
with the elliptic functions $g,\sigma$ defined in the Appendix. Here the
parameters $\zeta_i$ are related with $z_i$  by
$z_i=e^{2\pi\i\zeta_i}$.  As
usual, we let $r_{ij}=\pi_i\otimes \pi_j(r)$. \endproclaim

We will call the equations (9.2.10) {\it elliptic KZ equations}; they
are also known under the name Knizhnik-Zamolodchikov-Bernard
equations, since they first appeared in the papers of Bernard (see
\cite{Be}).

The remaining part of this section is devoted to the proof of this
theorem. 

\demo{Proof} The idea of proof is quite simple, though the
computations are a bit boring. Let us fix $u_1\in L_{\mu_1}^*, \dots,
u_n \in L_{\mu_n}^*$ and consider $\chi(u)=\<\chi, u_1\o u_2\o \dots\o
u_n\>\in \C$. Then, 
using the expression for the derivative
of $\Phi$, given in Theorem~9.1.4, we can
 write:

$$(K+h\v)z_i\frac{\d}{\d z_i}\chi(u) = \sum_a
\Tr \left(\Phi^1\dots 
\vdots\tJ_a(z_i) \Phi_{a u_i}(z_i)\vdots
 \dots \Phi^n e^{2\pi\i h}p^{-d}\right)-(\pi_i(h_\rho)\chi)(u),$$
where we write for brevity $\Phi^k$ instead of
$\Phi^k_{u_k}(z_k)$, and $a$ is any orthonormal basis in
$\g$. Since $\sum a\o a=\sum_{\a\in R^+} (e_\a\o f_\a+f_\a\o e_\a) +
\sum_l x_l\o x_l$, where $x_l$ is an orthonormal basis in $\h$, we get
the following expression: 

$$(K+h\v)z_i\frac{\d}{\d z_i}\chi(u)= 
\sum_{\a\in R^+} (a_\a^+ + a_\a^-+b_\a^+ +b_\a ^-) +\sum_l
(c_l^+ +c_l^-) - (\pi_i(h_\rho)\chi)(u),\tag 9.2.11$$
where

$$\aligned
a_\a^+&=\Tr \left(\Phi^1\dots \tJ^+_{e_\a}(z_i)\Phi^i_{f_\a u_i}(z_i) \dots
\Phi^n e^{2\pi\i h}p^{-d}\right),\\
a_\a^-&=-\Tr \left(\Phi^1\dots \Phi^i_{f_\a u_i}(z_i)\tJ^-_{e_\a}(z_i) \dots
\Phi^n e^{2\pi\i h}p^{-d}\right),\\
b_\a^+&=\Tr \left(\Phi^1\dots \tJ^+_{f_\a}(z_i)\Phi^i_{e_\a u_i}(z_i) \dots
\Phi^n e^{2\pi\i h}p^{-d}\right),\\
b_\a^-&=-\Tr \left(\Phi^1\dots \Phi^i_{e_\a u_i}(z_i)\tJ^-_{f_\a}(z_i) \dots
\Phi^n e^{2\pi\i h}p^{-d}\right),\\
c_l^+&=\Tr \left(\Phi^1\dots \tJ^+_{x_l}(z_i)\Phi^i_{x_l u_i}(z_i) \dots
\Phi^n e^{2\pi\i h}p^{-d}\right),\\
c_l^-&=-\Tr \left(\Phi^1\dots \Phi^i_{x_l u_i}(z_i)\tJ^-_{x_l}(z_i) \dots
\Phi^n e^{2\pi\i h}p^{-d}\right).\endaligned\tag 9.2.12$$

Let us calculate $a_\a^+$. Using the commutation relation 

$$\Phi_u(z_j)\tJ^+_{e_\a}(z_i)=\tJ^+_{e_\a}(z_i)\Phi_u(z_j)
-\frac{z_i}{z_j-z_i}\Phi_{e_\a u}(z_j)$$
for $|z_j|>|z_i|$ (see Lemma~9.1.3), we get 

$$ \multline
a_\a^+=\sum_{j<i}-\frac{z_i}{z_j-z_i} \chi(u_1, \dots, e_\a u_j,
\dots, f_\a u_i,\dots)\\
+\Tr \left(\tJ^+_{e_\a}(z_i)\Phi^1\dots \Phi^i_{f_\a u_i}(z_i) \dots
\Phi^n e^{2\pi\i h}p^{-d}\right).\endmultline$$

Using the cyclic property of the trace and the obvious relation 
$$e^{2\pi\i h}p^{-d}\tJ^+_{e_\a}(z_i)=e^{2\pi\i \<h, \a\>}\tJ^+_{e_\a}(pz_i)
e^{2\pi\i h}p^{-d},$$
we can rewrite it as follows: 

$$ \multline
a_\a^+=\sum_{j<i}-\frac{z_i}{z_j-z_i} \chi(u_1, \dots, e_\a u_j,
\dots, f_\a u_i,\dots)\\
+e^{2\pi\i \<h, \a\>}\Tr \left(\Phi^1\dots \Phi^i_{f_\a u_i}(z_i) \dots
\Phi^n \tJ^+_{e_\a}(pz_i)e^{2\pi\i h}p^{-d}\right).\endmultline$$

Now we again can commute $\tJ^+_{e_\a}$ with all $\Phi^j$, and so
on. It follows from the definition of the domain $\Cal D$, see
(9.2.7), that 
$\lim_{m\to \infty}e^{2\pi\i m\<h, \a\>}\tJ^+_{e_\a}(p^m z_i)=0$.  
Thus,  we get the following expression: 

$$\multline
a_\a^+=\biggl(\sum_{m\ge 0} \sum_{j<i}
\frac{p^m e^{2\pi\i m\<h, \a\>}z_i}{p^m z_i-z_j} \pi_i\o \pi_j(f_\a\o
e_\a) \\
+  \sum_{m> 0} \sum_{j\ge i}
\frac{p^m e^{2\pi\i m\<h, \a\>}z_i}{p^m z_i-z_j} \pi_i\o \pi_j(f_\a\o
e_\a)\biggr)\chi (u),\endmultline\tag 9.2.13$$
where for $j=i$ we let, by definition, $\pi_i\o \pi_j(f_\a\o e_\a)
=\pi_i(f_\a e_\a)$. Repeating similar arguments for $a_\a^-$ (but
moving $\tJ^-$ in the opposite direction), we get the following
result: 

$$
a_\a^+ + a_\a^-=\biggl(\sum_{m, j}
\frac{p^m e^{2\pi\i m\<h, \a\>}z_i}{p^m z_i-z_j} \pi_i\o \pi_j(f_\a\o
e_\a)\biggr) \chi(u),\tag 9.2.14$$
where the sum is taken over all pairs $m\in \Z, j=1,\dots, n$ except
the pair $(m=0, j=i)$, and we have the same convention as before for
the term with $j=i$. 

In the same way we get the following expression: 

$$
b_\a^+ + b_\a^-=\biggl(\sum_{m, j}
\frac{z_j e^{-2\pi\i m\<h, \a\>}}{p^m z_i-z_j} \pi_i\o \pi_j(e_\a\o
f_\a)\biggr) \chi(u),\tag 9.2.15$$
with the same conventions as before. 

As for the terms $c_l^\pm$, we must be more careful,
since $\lim_{m\to \infty}\tJ^\pm_{x_l}(p^{\pm m}z_i)=\pm \half
x_l$. This gives us

$$\aligned 
c_l^+ + c_l^-=&\biggl(\sum_{m, j}\half
\frac{p^m z_i+z_j}{p^m z_i-z_j} \pi_i\o \pi_j(x_l\o x_l)\biggr)
\chi(u)\\
&\qquad +\Tr \bigl(\Phi^1\dots \Phi^i_{x_l u_i}(z_i)\dots \Phi^n x_l e^{2\pi\i
h}p^{-d}\bigr)\\
=&\biggl(\sum_{m, j}\half
\frac{p^m z_i+z_j}{p^m z_i-z_j} \pi_i\o \pi_j(x_l\o
x_l)-\pi_i(x_l)\d_{x_l}\biggr) \chi (u), \endaligned \tag 9.2.16$$
where in addition to the previous conventions we must also specify
that the summation is done in the following order: $\sum_{m,j}=\sum_m
\sum_j$. This can be replaced by $\sum_{m,j}=\sum_j \text{V.P.} \sum
_m$, where the principal value of a  series is given by 

$$\text{V.P.} \sum_{m\in \Z}a_m= \lim_{N\to +\infty} \sum_{m=-N}^N
a_m.\tag 9.2.17$$

Summarizing the previous calculations, we get the following answer:

$$(K+h\v)z_i\frac{\d}{\d z_i} \chi =
\biggl(\sum_{j, m} r^m_{ij}(z_i/z_j)-\sum
\pi_i(x_l)\d_{x_l}-\pi_i(h_\rho)\biggr) \chi,\tag 9.2.18$$
where the sum is taken over  all pairs $m\in \Z, j=1,\dots, n$ except
the pair $(m=0, j=i)$, the order of summation is given by 
$ \sum_{m,j}=\sum_j \text{V.P.} \sum
_m$,  
and 

$$\aligned
r^m_{ij}(z_i/z_j)=&\sum_{\a\in R^+}
 \biggl(\frac{p^m e^{2\pi\i m\<h, \a\>}z_i}{p^m z_i-z_j} \pi_i\o \pi_j(f_\a\o
e_\a)\\
&\qquad \qquad +\frac{z_j e^{-2\pi\i m\<h, \a\>}}{p^m z_i-z_j} 
		\pi_i\o \pi_j(e_\a\o f_\a)\biggr)\\
&+\half\sum_l \frac{p^m z_i+z_j}{p^m z_i-z_j} \pi_i\o \pi_j(x_l\o
x_l).\endaligned \tag 9.2.19$$

Recall that by definition  for $j=i$ we let $\pi_i\o \pi_j(a\o
b)=\pi_i(ab)$.

We can simplify the expression above. Namely, note that the terms with
$j=i$ in (9.2.18) have the form: 
$$\aligned
\sum_{m\ne 0} r^m_{ii}(1)=&
\pi_i\left(\sum_{\a\in R^+} \sum_{m\ne 0}
	\biggl(\frac{p^m e^{2\pi\i m \<h, \a\>}}{p^m-1} f_\a e_\a  
		+\frac{ e^{-2\pi\i m \<h, \a\>}}{p^m-1} e_\a f_\a
	\biggr)
      \right)\\
&\quad + \pi_i\biggl(\text{V.P.} \sum_l \half \frac{p^m+1}{p^m-1}x_l^2
	      \biggr)\\
=&\pi_i\biggl(\sum_{\a\in R^+} \sum_{m\ne 0}
		\frac{ e^{-2\pi\i m \<h, \a\>}}{p^m-1} [e_\a, f_\a]
	\biggr)\\
=&\pi_i\biggl(
\sum_{\a\in R^+} -\biggl(\sigma(-\<\a, h\>)-\half\biggr)h_\a\biggr)\\
=& \pi_i\biggl( h_\rho +\sum_{\a\in R^+} \sigma(\<\a, h\>)h_\a\biggr),
\endaligned\tag 9.2.20$$
where $\sigma(x)$ is the logarithmic derivative of theta-function
defined in the Appendix (see formulas (A2), (A7)).  

It turns out that if we write the differential equation for the
product $\chi\dhat$ then these terms cancel. Indeed, since $\dhat$
does not depend on $z$, we have 

$$\aligned
(K+h\v)z_i&\frac{\d }{\d z_i} (\chi\dhat)= 
\biggl((K+h\v)z_i\frac{\d }{\d z_i} \chi\biggr) \dhat\\
&=\dhat \biggl(\sum r_{ij}^m(z_i/z_j)-\sum
\pi_i(x_l)\d_{x_l}-\pi_i(h_\rho)\biggr) \chi \\
&=\biggl(\sum r_{ij}^m(z_i/z_j)-\sum
\pi_i(x_l)\d_{x_l}-\pi_i(h_\rho)\biggr) (\chi\dhat) 
+\sum_l \frac{\d_{x_l} \dhat}{\dhat}
\pi_i(x_l)(\chi\dhat).\endaligned\tag 9.2.21$$

Since $\dhat=a(\tau) \prod_{\a\in R^+} \theta_1(\<\a, h\>)$ (see
(8.1.4)), we have:

$$\frac{\d_{x_l}\dhat}{\dhat}=\frac{1}{2\pi\i} \sum_{\a\in R^+} 
\frac{\theta'_1(\<\a, h\>)}{\theta_1(\<\a, h\>)} \<\a, x_l\>=
-\sum_{\a\in R^+}\sigma(\<\a, h\>)\<\a, x_l\>. $$

Therefore, 

$$\sum_l \frac{\d_{x_l} \dhat}{\dhat}
\pi_i(x_l)(\chi\dhat)=-\sum_{\a\in R^+} \sigma(\<\a, h\>)\pi_i(h_\a)
(\chi\dhat).$$

Substituting this in (9.2.21) we see that this term cancels the terms
$r^m_{ii}$. Introducing $\zeta_i$ such that $z_i=e^{2\pi\i \zeta_i}$
and using the formulas from the Appendix, we can rewrite
$r_{ij}, j\ne  i$ in terms of elliptic functions $g, \sigma$, which
gives us the statement of the Theorem. 
\qed
\enddemo

\remark{Remark}
In fact, all the constructions of this section can be
generalized to the case when $\l$ is not necessarily dominant; 
we could replace the conditions 
$\l\in P^+, K\notin\Bbb Q$ by the following conditions:

$$\gathered
(\lhat+\rhat, \ahat)-\frac{N}{2} (\ahat, \ahat)=0\text { with }N\in
\N,
\ahat\in \Rhat^+\\
\text{is only possible if $\ahat\in R^+$}.\endgathered$$

These conditions ensure that $L_{\l, K}$ is induced from the module
$L_\l$ over $\g$, and all the arguments of Sections~9.1, 9.2 can be
repeated literally. 

We can also include in the consideration integrable modules; in this
case, of course, the differential equations are the same, but it is
not true that the space of intertwiners $\Phi$ of form (9.1.5) is
isomorphic to the space of $\g$-intertwiners (see, for example,
\cite{TK}, where this is discussed in detail for
$\g=\sltwo$). \endremark

\head 9.3 Monodromies of elliptic KZ equations\endhead

In this section we study some properties of the elliptic KZ equations
(9.2.10). Let us consider the trace $\chi=\chi(z_1, \dots, z_n, h, p)$
defined by (9.2.1) as a function of $z_1, \dots, z_n, p$ with values
in the space 

$$\V= (L_{\mu_1}\o \dots \o L_{\mu_n})[0] \o C^\infty(\h).\tag 9.3.1$$

Then the elliptic KZ equations can be rewritten in the following form:

$$(K+h\v)\frac{\d}{\d \zeta_i}\varphi=A_i\varphi , \tag 9.3.2$$
where $\varphi =\chi\dhat$, and the operators $A_i:\V\to \V$ are
defined by  

$$A_i(\zeta_1, \dots, \zeta_n)=\sum_{j\ne i}r_{ij}(\zeta_i-\zeta_j)
-2\pi\i \sum_l \pi_i(x_l)\d_{x_l},\tag 9.3.3$$
where all the notations are as in Theorem~9.2.5; in particular,
$r(\zeta)$ is defined by (9.2.11). 

\proclaim{Lemma 9.3.1} The elliptic KZ system is consistent, i.e., 
$$[(K+h\v)\frac{\d}{\d \zeta_i}-A_i,(K+h\v) \frac{\d}{\d \zeta_j}-A_j]=0
	\tag 9.3.4$$
as operators in $\V$.\endproclaim

\demo{Proof} The proof is based on the fact that we have sufficiently
many solutions of this system given by the correlation functions.
  Indeed, let us consider the commutator (9.3.4); denote it by $D_{ij}$. It is
easy to see that it does not contain derivatives in $\zeta_i$ and thus is
some operator in $\V$. More precisely, it is some first order differential
operator in $\h$ tensored by some operator in $(L_{\mu_1}\o \dots\o
L_{\mu_n})[0]$.  It follows from Theorem~9.2.4 that
it annihilates all the correlation functions
$\varphi=\chi\dhat$. 

Now, let us take $\l$ such that $\<\l, \a_i\v\>$ are large enough. It
is known that we can do it in such a way that the space of all
$\g$-intertwiners $\Psi: L_\l\to L_\l\o V$
(recall that $V=L_{\mu_1}\o \dots \o L_{\mu_n}$) is
isomorphic to the zero-weight subspace $V[0]$. 

In this case it is easy to show, using Theorem~9.2.1 and some
properties of correlation functions on the sphere, that Lemma~9.2.2
can be rewritten as follows:

\proclaim{Lemma 9.3.2} Let us fix $z_1, \dots, z_n$ such that
$|z_1|>\dots >|z_n|$. Then in the limit $\tau\to +\i\infty$ the trace
$\chi$ has the following asymptotics: 

$$\chi(z_1, \dots, z_n,p, h-\frac{\tau\rho}{h\v})\simeq p^{\Delta(\l)
-\frac{(\l, \rho)}{h\v}} e^{2\pi\i \<h, \l\>} v_\chi$$
for some vector $v_\chi\in V[0]$. Moreover, for generic $z_i$ the
mapping $\chi\mapsto v_\chi$ is an isomorphism of the space of all
traces of the form \rom{(9.2.1)} for all possible choices of $\l_i$
and the space $V[0]$. \endproclaim

Therefore, we see that since the commutator $D_{ij}$ annihilates
every trace $\chi\delta$, it also annihilates every function of the
form  $\dhat e^{2\pi\i \<h, \l\>} v_\chi$ for every sufficiently large
$\l$ and small enough $p$. But this is only possible if it is identically
zero. \qed\enddemo

Note that straightforward proof of the consistency would require use
of Riemann identities for theta-functions, and would necessarily be
very tiresome.

Now, let us study some properties of the elliptic KZ system. First of
all, we want to know the elliptic properties, i.e. the behavior of
the coefficients under the translations $\zeta_i\mapsto \zeta_i+1,
\zeta_i\mapsto \zeta_i+\tau$. Explicit calculation gives the following
result: 

\proclaim{Theorem 9.3.3} Let $r(\zeta)$ be defined by
\rom{(9.2.11)}. Then 
\roster 
\item $r_{12}(\zeta)+r_{21}(-\zeta)=0$, where $r_{12}=r, r_{21}=P(r),
P(a\o b)=b\o a$.  

\item $r$ has a pole of first order at $\zeta=0$, and 
$\text{Res}_{\zeta=0} r(\zeta)=\Omega$, where $\Omega$ is the canonical
$\g$-invariant element in $\g\otimes \g$ defined by \rom{(1.1.3)}.  

\item 
$$\aligned
r(\zeta+1)&=r(\zeta),\\
r(\zeta+\tau)&=(e^{2\pi\i h}\otimes 1)\circ r(\zeta)
	 \circ (e^{-2\pi\i h}\otimes 1) -2\pi\i \sum_l x_l\otimes x_l\\
	&= (1\otimes e^{-2\pi\i h}) \circ r(\zeta) 
	\circ (1\otimes e^{2\pi\i h})-2\pi\i \sum x_l\otimes x_l.
\endaligned\tag 9.3.5$$
\endroster
\endproclaim

\proclaim{Corollary 9.3.4} Let $A_i$ be the operator in $\V$
defined by \rom{(9.3.3)}. Then for any $j=1, \dots, n$ \rom{(}including 
$j=i$\rom{)}, we have 
$$\gathered
A_i(\zeta_1,\ldots, \zeta_j+1,\ldots \zeta_n)=A_i(\zeta_1,\ldots,
\zeta_n), \\ 
A_i(\zeta_1,\ldots, \zeta_j+\tau,\ldots \zeta_n)=
\pi_j(e^{2\pi\i h}) A_i(\zeta_1,\ldots,\zeta_n)\pi_j(e^{-2\pi\i h}).
\endgathered\tag 9.3.6$$
\endproclaim

Thus, the elliptic KZ system is ``almost'' invariant under the elliptic
transformations. Let us modify it so that it becomes really invariant.

Introduce the following operator in $\V$: 

$$\psi(\zeta_1,\ldots,\zeta_n,h)=\prod\limits_{i<j}
\frac{\theta_1\biggl(\zeta_i-\zeta_j+\frac{1}{n}(\pi_i(h)-\pi_j(h))
\biggr)}{\theta_1(\zeta_i-\zeta_j)};\tag 9.3.7$$
then

$$\gathered
\psi(\zeta_1,\ldots \zeta_i+1,\ldots,
\zeta_n)=\psi(\zeta_1,\ldots,\zeta_n),\\
\psi(\zeta_1,\ldots \zeta_i+\tau,\ldots,
\zeta_n)=\pi_i(e^{-2\pi\i
h})\psi(\zeta_1,\ldots,\zeta_n).\endgathered\tag 9.3.8$$

Then we have the following theorem: 

\proclaim{Theorem 9.3.5} 
\roster\item 
 If $\varphi$ is a solution of the elliptic KZ system \rom{(9.3.2)} then
$\tilde \varphi=\psi \varphi$, where $\psi$ is given by \rom{(9.3.7)}, is a 
solution of the system
$$(K+h\v)\frac{\d }{\d \zeta_i}\tilde \varphi =\tilde
A_i(\zeta_1,\ldots,\zeta_n) \tilde \varphi,\tag 9.3.9$$
where 
$$\tilde A_i= \psi A_i \psi^{-1} +(K+h\v)\frac{\d }{\d \zeta_i}
\log \psi.\tag 9.3.10$$

\item This new system is invariant under the translations by $1,\tau$:
for every $j\in 1, \dots, n$ 

$$\gather
\tilde A_i(\zeta_1,\ldots \zeta_j+1,\ldots,
\zeta_n)=\tilde A_i(\zeta_1,\ldots,\zeta_n),\\
\tilde A_i(\zeta_1,\ldots \zeta_j+\tau,\ldots,
\zeta_n)=\tilde A_i(\zeta_1,\ldots,\zeta_n).\endgather$$
\endroster
\endproclaim

This means that we can consider (9.3.9) as a local system on the torus
$X=E^n\setminus\{\zeta_i=\zeta_j\}$, where 
$E=\C/(\Z+\Z\tau)$. Therefore, it makes sense to consider monodromies of
this local system. It will give us some representation of the
fundamental group of $X$ in the space $\V$.
This fundamental group, called the pure
braid group of the torus, is well known and can be described in terms
of generators and relations. It is convenient to consider a larger
group: the braid group of the torus $BT_n=\pi_1(X/S_n)$, where the
symmetric group $S_n$ acts in $X$ by permutation of coordinates. Then 
$\pi_1(X)=\text{Ker }\sigma$, where $\sigma$ is the natural
homomorphism $\sigma:BT_n\to S_n$. The local system above defines a
representation of $BT_n$ in the following sense: to each $\gamma\in
BT_n$ we assign a holonomy operator
 
$$M(\gamma):\bigl(L_{\mu_1}\o\dots\o L_{\mu_n}\bigr) [0]
\o C^{\infty}(\h)\to
 \bigl(L_{\mu_{\sigma(1)}}\o\dots\o L_{\mu_{\sigma(n)}}\bigr)
[0] \otimes C^{\infty}(\h),$$
where $\sigma=\sigma(\gamma)$. 

To describe this representation, recall that the braid group $BT_n$ is
generated by the elements $T_i, 1\le i\le n-1, X_j, Y_j, 1\le j \le
n$ (in fact, it suffices to take only $T_i, X_1, Y_1$). Their
geometric meaning is as follows: if we take a base point
$\zeta=(\zeta_1,\ldots,\zeta_n)$ such that $\zeta_i=x_i+ y_i \tau,\,  
x_i,y_i\in\R, 0<x_1<\ldots <x_n<1, 0<y_1<\ldots<y_n<1$ then $T_i$
corresponds to a transposition
of $\zeta_i$ and $\zeta_{i+1}$ such that $\zeta_i$ passes 
$\zeta_{i+1}$ from the left and $X_i,Y_i$ correspond to $\zeta_i$
going  around the $1$ and $\tau$ cycles  in negative direction:

\psfig{figure=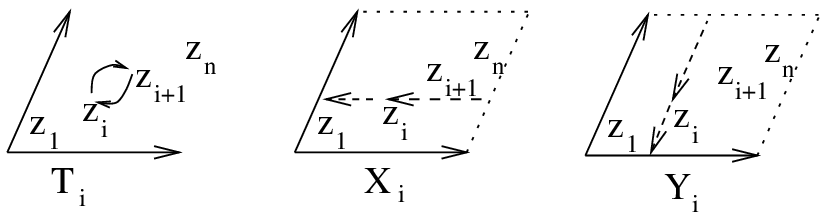} 

To calculate the monodromies of the local system above we use the
fact that we have a large number of solutions given by the 
functions $\tilde \varphi=\psi\dhat\chi$, where $\chi$ is
the correlation function on the torus (9.2.1). 

Recall the known results about "commutation" of vertex operators. 
Denote by $H_\l^{\nu\mu}$ the space of $\g$-intertwiners $L_\l\to
L_\nu\o L_\mu$:

$$H_\l^{\nu\mu}=\Hom_\g (L_\l, L_\nu\o L_\mu).$$

 Then, due to Lemma~9.1.2, for every $g\in
H_\l^{\nu\mu}$ we have a unique $\ghat$-intertwiner $\Phi^g(z): L_{\l,
K}\to L_{\nu, K}\o L_\mu(z)$. Similarly, if we have weights $\l_0, \l_1,
\l_2$ and $\mu_1, \mu_2$ and 
$g_i\in H_{\l_i}^{\l_{i-1}\mu_i}$ then we can define a
$\ghat$-intertwiner

$$\Phi^{g_1}(z_1)\Phi^{g_2}(z_2): L_{\l_2, K}\to L_{\l_0, K}\o
L_{\mu_1}(z_1) \o L_{\mu_2}(z_2).\tag 9.3.11$$

Summing over all possible $\l_1$, we can assign to every intertwiner
$g$ from 

$$\bigoplus_{\l_1} H_{\l_1}^{\l_0 \mu_1}\o H_{\l_2}^{\l_1 \mu_2}=
\Hom_\g (L_{\l_2}, L_{\l_0}\o L_{\mu_1}\o L_{\mu_2})$$ 
(note that this sum is in fact finite) a corresponding
$\ghat$-intertwiner $\Phi^{g^{(1)}} (z_1)\Phi^{g^{(2)}} (z_2)= 
\sum \Phi^{g^{(1)}_i} (z_1)\Phi^{g^{(2)}_i} (z_2)$ if 
$g=\sum g^{(1)}_i\o g^{(2)}_i$. 

Note that  the product $\Phi^{g^{(1)}} (z_1)\Phi^{g^{(2)}} (z_2)$ is a
well-defined operator for $|z_1|>|z_2|, z_i\notin \R_-$ in the sense 
that all its
matrix coefficients are well-defined analytic functions in this
region. It turns out that we can analytically continue this product as
a multivalued function in the whole region $z_1\ne z_2, z_i\ne
0$. Denote by $A^+$ the analytic continuation along the path when
$z_1$ passes $z_2$ from the left:
$$\psfig{figure=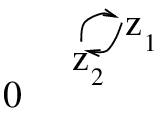}$$

. The following theorem is well-known
(see, for example, \cite{MS}). 

\proclaim{Theorem 9.3.6}
 As before, let us fix the weights $\l_0,
 \l_2, \mu_1, \mu_2$ and define 

$$
 H^{\mu_1\mu_2} =\bigoplus_{\l_1} H_{\l_1}^{\l_0 \mu_1}\o
H_{\l_2}^{\l_1 \mu_2}=\Hom_\g (L_{\l_2}, L_{\l_0}\o L_{\mu_1}\o
L_{\mu_2}). \tag 9.3.12$$

Then there exists some linear operator $R: H^{\mu_1\mu_2}\to
H^{\mu_2\mu_1}$  such that for every
$g\in H^{\mu_1\mu_2}$ we have the following identity in the region
$|z_2|>|z_1|$: 
 
$$A^+(\Phi^{g^{(1)}}(z_1)\Phi^{g^{(2)}} (z_2))=
\Phi^{(Rg)^{(1)}}(z_2) \Phi^{(Rg)^{(2)}}(z_1).\tag 9.3.13$$
\endproclaim

\remark{Remark} In fact, this matrix $R$ coincides (in a suitably
chosen basis) with the matrix $\check R=PR$, where $R$ is the $R$-matrix
for the corresponding quantum group. This is probably one of the
deepest results in conformal field theory; see its discussion, for
example, in the book of Kassel \cite{Kas} or in \cite{KL}. However, we
won't need it. 
 \endremark

This statement allows us to find the monodromies of the elliptic KZ
system. Let us recall (see Lemma~9.1.2) that for every collection of
$\g$-intertwiners  
$g^i: L_{\l_i}\to L_{\l_{i-1}}\o L_{\mu_i}$ we can define the
corresponding $\ghat$-intertwiners $\Phi^{g^i}(z_i): L_{\l_i, K}\to
L_{\l_{i-1}, K}\o L_{\mu_i}(z_i)$, and the corresponding trace 
$\chi^{g^1\o \dots g^n}$, defined by (9.2.1). We can extend this by
linearity to define $\chi^g$ for arbitrary intertwiner

$$\aligned 
g\in H^{\mu_1\dots \mu_n} = &\Hom_\g\bigl(L_\l,  L_\l\o
L_{\mu_1}\o\dots \o L_{\mu_n}\bigr)\\ 
& =\bigoplus_{\l_1, \dots, \l_{n-1}} H_{\l_1}^{\l
\mu_1}\o H_{\l_2}^{\l_1\mu_2}\o  \dots \o
H_{\l}^{\l_{n-1}\mu_n}.\endaligned \tag 9.3.14 $$

Let us fix the base point $\zeta^0$ as before
and  define the the subspace $\V^0\subset \V$
which is spanned by the expressions of the form 

$$\tilde \varphi^g=\chi^g_{\l, K}(\zeta^0, p, h) \dhat(h)\psi(\zeta^0,
p,h),\tag 9.3.15$$
for all possible choices of $\l, g$. In some sense, 
we could say that such expressions
form a basis in $\V$; however, since this space is infinite-dimensional
this requires more accurate approach, and we prefer just to speak of the
subspace $\V^0$ in $\V$ formed by finite linear combinations of the
functions of the form (9.3.15). 

\proclaim{Theorem 9.3.7} The  representation $M: BT_n\to
\operatorname{End } \V$, given by the monodromies of the local system
\rom{(9.3.9)}, preserves the subspace $\V^0$ and in this subspace 
is given by 

$$\aligned
T_i&:\tilde \varphi^g \mapsto \tilde \varphi^{R_{i, i+1}(g)}, \\
X_i&: \tilde \varphi^g \mapsto  e^{2\pi\i
(\Delta_{i}-\Delta_{i-1})}\tilde \varphi^g, \\
Y_i&: \tilde \varphi^g \mapsto \tilde \varphi^{g'},\\
&\qquad g'= R_{i+1, i}R_{i+2, i}\dots R_{n, i}R_{1,i}\dots
R_{i-1,i}g.
\endaligned\tag 9.3.16
$$

Here $R_{i i+1}$ is defined as follows: 

$$R_{i i+1}=1\o\dots R\o \dots\o 1:H^{\mu_1\dots \mu_n} \to
H^{\mu_1\dots\mu_{i+1} \mu_i\dots \mu_n},\tag 9.3.17$$  
where $H^{\mu_1\dots \mu_n}$ is defined by \rom{(9.3.14)} and the
operator $R$ was defined in Theorem~\rom{9.3.6}. 

The numbers $\Delta_i$ are defined as follows: 
if $g\in H^{\mu_1\dots \mu_n}$ lies
in the component  $H_{\l_1}^{\l
\mu_1} \o \dots \o H_{\l}^{\l_{n-1}\mu_n}$ for some $\l_1, \dots, \l_{n-1}$
then 
$\Delta_i=\Delta(\l_i)$ \rom{(}see \rom{(9.1.2))}.
\endproclaim

\demo{Proof} Let us assume that $g$ is homogeneous, i.e. 
$g\in H_{\l_1}^{\l
\mu_1} \o \dots\o  H_{\l}^{\l_{n-1}\mu_n}$ for some $\l_1, \dots, \l_{n-1}$.
 Let $\tilde \varphi^g(\zeta^0, h,p)\in \V^0$ be a function of the
form (9.3.15). 
Consider the corresponding function $\tilde \varphi^g(\zeta, h,p)$,
for arbitrary $\zeta$.  Then we know that $\tilde
\varphi^g(\zeta, h,p)$   is a solution of the
elliptic KZ system (9.3.9) (see Theorem~9.3.5). Thus, calculation, for
example, of 
$X_i \tilde \varphi^g$ is equivalent to finding the result of the
analytic continuation of $\tilde \varphi^g$ along the path $X_i$, i.e. the path
which connects $\zeta^0=(\zeta^0_1, \dots, \zeta^0_n)$ and $(\zeta^0_1,
\dots, \zeta^0_i-1, \dots, \zeta^0_n)$. Since the dependence of $\tilde
\varphi$ on $\zeta_i$ has the form 

$$\tilde\varphi =e^{-2\pi\i\zeta_i
(\Delta(\l_i)-\Delta(\l_{i-1}))}\cdot\text{\it some periodic function of
$\zeta_i$},$$
 we get the formula for action of $X_i$. 

The action of $T_i$ can be calculated quite similarly; it is essentially
the statement on the commutation of intertwiners (Theorem~9.3.6) 

To prove the statement about $Y_i$,
we should commute the operator $\Phi^{g^{(i)}}(z_i)$ with all
$\Phi^{g^{(k)}}(z_k), 
1\le k\le i-1$; then use the cyclic property of trace and the
relation 
$p^{-d}e^{2\pi\i h}\Phi(z)=\pi(e^{2\pi\i h})\Phi(pz)p^{-d}e^{2\pi\i
h}$. The factor $\pi(e^{2\pi\i h})$ is canceled by the factor
appearing in the transformation law for $\psi$ (see (9.3.8)).\qed
\enddemo


\vfill\newpage


\vbox{\vskip 1in}
\specialhead\chapter{10} 
\centerline{PERSPECTIVES AND OPEN QUESTIONS} 
\endspecialhead

\vbox{\vskip 0.5in}

In this chapter we briefly review the results and list some open
questions. So far, we have defined generalized characters for simple
finite-dimensional Lie algebras, corresponding quantum groups and affine
Lie algebras. We also have derived a number of interesting properties
of these characters and have shown that for $\g=\sln$ and some special
choice of $U$ they coincide (after renormalization) with well-known
special functions -- Jack polynomials; similarly, for $U_q\sln$ we get
$q$-analogue of Jack polynomials -- Macdonald's polynomials of type
$A_{n-1}$. This construction allowed us to reprove certain 
 identities for Macdonald's polynomials. We also generalized a large
part of this theory to the affine case. 

However, a number of interesting questions remain open. Here are some
of them. 

\subhead 1. Relation with affine Hecke algebras \endsubhead
In a recent series of papers (\cite{C1, C2}), Ivan Cherednik has shown
that Macdonald's polynomials naturally appear in representations of
double affine Hecke algebras. Using this technique, he proved Macdonald's
inner product conjectures for arbitrary root systems. It would be very
interesting to establish a relation of our construction with his work.
Morally, it should be something like the duality between representations
of $GL_n$ and $S_m$; so far, no results in this direction are known. 

\subhead 2. Critical level limit\endsubhead 
So far, we have assumed in the affine case that the level $K$  of
highest-weight modules is not equal to the dual Coxeter number $h\v$.
However, from many points of view the case $K=-h\v$ (critical level) 
 is very interesting.
For example, it was shown by E.~Frenkel and B.~Feigin that (a certain
completion of) the algebra $U\ghat$ has a large center at the critical
level. Using this result, Etingof has shown that the elliptic
Calogero-Sutherland operator can be included in a large commutative
family of differential operators (see \cite{E1}). 

However, it is easy to check that one can not define the generalized
characters at the critical level in the same manner as we have done
before, since in this case there are no intertwining operators of the
form we need. One way to study the critical level is to  study the
asymptotics of generalized characters as $K\to -h\v$. For $\g=\sltwo$ it
was done in the paper \cite{EK3}, using explicit integral formulas for
the generalized characters. This allowed us to recover the classical
formulas for Lam\'e functions. In general case, no integral formulas are
known so far; even if one writes such formulas (which can be done),
finding their asymptotics is technically rather difficult. 

It seems plausible that there is a direct way to get the generalized
characters at the critical level, using the language of coinvariants
rather than intertwining operators. Such a construction should be
closely related with geometric Langlands correspondence as defined in
recent papers of Beilinson and Drinfeld.

\subhead 3. Inner product in the affine case \endsubhead
In the finite-dimensional case, we had two ways to characterize
generalized characters: in terms of differential equations or in terms
of inner product (due to the orthogonality theorem). In the affine case, 
the differential equations work just fine, but defining inner product
is quite complicated: any straightforward approach fails (the formulas
have singularities). However, some arguments from conformal field theory
suggest that there exists a natural inner product on the generalized
characters such that we have an analogue of the orthogonality theorem.
Also, this inner product should be modular invariant. This question 
is the subject of the author's present work (\cite{K2}).


\vfill\newpage

\vbox{\vskip 0.5in}
\specialhead \chapter\nofrills{Appendix}
\centerline{LIST OF FORMULAS RELATED TO ELLIPTIC FUNCTIONS} 
\endspecialhead

\vbox{\vskip 0.5in}

In this Appendix, we briefly list the main definitions and formulas for
elliptic (theta-, sigma-) functions we are using. Many of them can be
found in \cite{WW, Mu}; others can be proved in the usual way. 

Throughout this section, we assume that we are given a number 
$p=e^{2\pi\i\tau}\in \C,|p|<1$, and all our functions depend on it.

In some formulas, we use the notion of the principal values of a
series, which is defined as follows: 

$$\text{V.P.}\sum_{n\in \Z} a_n =\lim_{N\to \infty} \sum_{n=-N}^N
a_n.$$

\head Functions and identities\endhead

$$\gather
\theta_1(x)=2p^{1/8}\sin \pi x
	\prod_{n\ge 1} (1-e^{2\pi \i x}p^n)(1-e^{-2\pi\i x}p^n)(1-p^n)
	\tag A1\\
\sigma(x)=-\frac{1}{2\pi\i}\frac{\theta'_1(x)}{\theta_1(x)}=
\sum_{n\ge 0} \frac{p^n e^{2\pi\i x}}{1- p^n e^{2\pi\i x}}+
\sum_{n< 0} \frac{1}{1- p^n e^{2\pi\i x}}+\frac{1}{2}
\tag A2\\	
\wp(x)=\frac{1}{x^2} + \sum_{(m,n)\in \Z^2\setminus\{(0,0)\} } 
\left(\frac{1}{(x-m-n\tau)^2}-\frac{1}{(m+n\tau)^2}\right)
\tag A3\\
\sum_{n\in \Z}\frac{p^ne^{2\pi\i x}}{(1-p^ne^{2\pi\i x})^2}
=\frac{1}{4\pi^2}\d_x^2\log \theta_1(x)=\frac{1}{2\pi\i}\sigma'(x)\\
\text{\phantom{asdfghjkl}}=-\frac{1}{4\pi^2}\wp(x)+c(\tau)\text{ for some constant }c(\tau)
\tag A4\\
\eta(\tau)=p^{\frac{1}{24}}\prod_{n\ge 1} (1-p^n)
	\tag A5\endgather$$

$$\gather
g(x,\zeta)=\sum_{m\in\Z}\frac{e^{2\pi\i m\zeta}}{1-p^me^{-2\pi\i x}}
	=-\frac{1}{2\pi\i } \frac{\theta_1(x-\zeta)\theta'_1(0)}
			    {\theta_1( x)\theta_1(\zeta)}
	\tag A6\\
\sum_{m\ne 0}\frac{e^{2\pi\i m\zeta}} {1-p^m} 
	=\sigma(\zeta)-\frac{1}{2}\tag A7\\
\frac{1}{2} \text{V.P.} \sum_{m\in \Z}
	\frac{1+p^me^{2\pi\i\zeta}}{1-p^me^{2\pi\i \zeta} } 
	=\sigma(\zeta)\tag A8\\
\phi(x,\zeta)=\sum_{m\in\Z} \frac{e^{-2\pi\i x}p^me^{2\pi\i m\zeta} }
				{(1-p^m e^{-2\pi\i x})^2}
	=-\frac{1}{2\pi\i}\frac{\d}{\d x} g(x,\zeta)\tag A9\\
\phi_0(\zeta)=\sum_{m\ne 0} \frac{p^me^{2\pi\i m\zeta}} {(1-p^m)^2}=
	\frac{\theta''_1(\zeta)}{\pi^2\theta_1(\zeta)} 
	-\frac{\theta'''_1(0)}{24\pi^3\theta_1(0)}+\frac{1}{12}\tag A10
\endgather$$

\head Symmetries and elliptic properties\endhead
\nopagebreak
$$\gathered
\theta_1(-x)=-\theta_1(x),\\
\theta_1(x+1)=-\theta_1(x),\qquad
\theta_1(x+\tau)=-p^{-1/2}e^{-2\pi \i x}
			\theta_1(x), \endgathered\tag A11$$

$$\gathered
\wp(-x)=\wp(x), \\
\wp(x+1)=\wp(x+\tau)=\wp(x),\endgathered\tag A12$$

$$\gathered
\sigma(-\zeta)=-\sigma(\zeta),\\
\sigma(\zeta+1)=\sigma(\zeta),\qquad \sigma(\zeta+\tau)=\sigma(\zeta)+1,
\endgathered\tag A13$$

$$\gathered
g(x,\zeta)=-g(\zeta,x)=-g(-x,-\zeta),\\
g(x+1,\zeta)=g(x,\zeta),\qquad g(x,\zeta+1)=g(x,\zeta),\\
g(x+\tau,\zeta)=e^{2\pi\i\zeta}g(x,\zeta),\qquad 
g(x,\zeta+\tau)=e^{2\pi\i x}g(x,\zeta).\endgathered\tag A14$$


\vfill\newpage


\enddocument
\end